\documentclass[journal,onecolumn,compsoc]{IEEEtran}

\usepackage{amsmath,amsfonts,amssymb}
\usepackage{amsthm}
\usepackage{mathrsfs} 
\usepackage{algorithm}
\usepackage{array}
\usepackage{textcomp}
\usepackage{stfloats}
\usepackage{url}
\usepackage{verbatim}
\usepackage{graphicx}
\usepackage{bm}
\usepackage{bbm}
\usepackage{mathtools}
\usepackage{subcaption} 
\usepackage{multicol} 
\usepackage{orcidlink} 
\usepackage{stmaryrd} 
\usepackage{empheq} 
\usepackage{enumitem} 
\usepackage{optidef}
\usepackage{algpseudocode} 

\newcommand{\nocontentsline}[3]{}
\newcommand{\tocless}[2]{\bgroup\let\addcontentsline=\nocontentsline#1{#2}\egroup}

\hypersetup{
    colorlinks,
    linkcolor={red!50!black},
    citecolor={blue!50!black},
    urlcolor={blue!80!black}
}

\newtheorem{theorem}{Theorem}
\newtheorem{lemma}{Lemma}
\newtheorem{corollary}{Corollary}
\newtheorem{proposition}{Proposition}

\newtheorem{definition}{Definition}

\newtheorem{assumption}{Assumption}
\newtheorem{remark}{Remark}
\newtheorem{fact}{Fact}

\newtheorem{procedure}{Procedure}

\algnewcommand{\LineComment}[1]{\State \(\triangleright\) #1} 

\algdef{SE}[DOWHILE]{Do}{doWhile}{\algorithmicdo}[1]{\algorithmicwhile\ #1}%

\begin{document}
\newcommand{\parens}[1]{\left( #1 \right )}
\newcommand{\brackets}[1]{\left[ #1 \right]}
\newcommand{\braces}[1]{\left\{ #1 \right\}}

\newcommand{\Real}{\mathbb{R}}
\newcommand{\Natural}{\mathbb{N}}

\newcommand{\Dcal}{\mathcal{D}}
\newcommand{\Ecal}{\mathcal{E}}
\newcommand{\Gcal}{\mathcal{G}}
\newcommand{\Ical}{\mathcal{I}}
\newcommand{\Jcal}{\mathcal{J}}
\newcommand{\Kcal}{\mathcal{K}}
\newcommand{\Lcal}{\mathcal{L}}
\newcommand{\Mcal}{\mathcal{M}}
\newcommand{\Pcal}{\mathcal{P}}
\newcommand{\Scal}{\mathcal{S}}
\newcommand{\Vcal}{\mathcal{V}}
\newcommand{\Ocal}{\mathcal{O}}

\newcommand{\Pscr}{\mathscr{P}}

\newcommand{\av}{\bm{a}}
\newcommand{\bv}{\bm{b}}
\newcommand{\cv}{\bm{c}}
\newcommand{\dv}{\bm{d}}
\newcommand{\pv}{\bm{p}}
\newcommand{\qv}{\bm{q}}
\newcommand{\rv}{\bm{r}}
\newcommand{\sv}{\bm{s}}
\newcommand{\tv}{\bm{t}}
\newcommand{\uv}{\bm{u}}
\newcommand{\vv}{\bm{v}}
\newcommand{\wv}{\bm{w}}
\newcommand{\xv}{\bm{x}}
\newcommand{\yv}{\bm{y}}
\newcommand{\zv}{\bm{z}}

\newcommand{\betav}{\bm{\beta}}
\newcommand{\thetav}{\bm{\theta}}
\newcommand{\zetav}{\bm{\zeta}}
\newcommand{\xiv}{\bm{\xi}}
\newcommand{\etav}{\bm{\eta}}
\newcommand{\muv}{\bm{\mu}}
\newcommand{\nuv}{\bm{\nu}}

\newcommand{\Imat}{\bm{I}}
\newcommand{\Am}{\bm{A}}
\newcommand{\Bm}{\bm{B}}
\newcommand{\Cm}{\bm{C}}
\newcommand{\Dm}{\bm{D}}
\newcommand{\Em}{\bm{E}}
\newcommand{\Fm}{\bm{F}}
\newcommand{\Gm}{\bm{G}}
\newcommand{\Lm}{\bm{L}}
\newcommand{\Mm}{\bm{M}}
\newcommand{\Pm}{\bm{P}}
\newcommand{\Qm}{\bm{Q}}
\newcommand{\Rm}{\bm{R}}
\newcommand{\Sm}{\bm{S}}
\newcommand{\Tm}{\bm{T}}
\newcommand{\Xm}{\bm{X}}
\newcommand{\Ym}{\bm{Y}}
\newcommand{\Zm}{\bm{Z}}

\newcommand{\Xim}{\bm{\Xi}}
\newcommand{\Lmdm}{\bm{\Lambda}}

\newcommand{\abs}[1]{\left\lvert #1 \right\rvert}
\newcommand{\norm}[1]{\left\lVert #1 \right\rVert}
\newcommand{\innerp}[2]{\langle #1, #2 \rangle}
\newcommand{\sign}{\mathrm{sign}}
\newcommand{\aff}{\mathrm{aff}}
\newcommand{\supp}{\mathrm{supp}}
\newcommand{\interior}{\mathrm{int}}
\newcommand{\closure}{\mathrm{cl}}

\newcommand{\trace}{\mathrm{tr}}
\newcommand{\rank}{\mathrm{rank}}
\newcommand{\diag}{\mathrm{diag}}
\newcommand{\blkdiag}{\mathrm{blkdiag}}
\newcommand{\Tr}{\top}

\newcommand{\Null}{\mathrm{Null}}
\newcommand{\Col}{\mathrm{Col}}

\newcommand{\dom}{\mathrm{dom\,}}
\newcommand{\levelset}{\mathrm{lev}}
\newcommand{\argmin}{\mathrm{arg\,min}}
\newcommand{\Prox}{\mathrm{Prox}}

\newcommand{\paragroup}{(\Am,\rho,\bv,\lambda)}
\newcommand{\paraspace}{\Real^{m\times n}\times [0,1)\times\Real^{2m}\times\Real_{++}}

\newcommand{\sGMCpara}{\text{sGMC}\paragroup}
\newcommand{\sGMClmd}{\text{sGMC}(\lambda)}

\newcommand{\SLA}{\Scal_{\text{LA}}}
\newcommand{\betaLA}{\betav_{\text{LA}}}
\newcommand{\gammaLA}{\gamma_{\text{LA}}}

\newcommand{\Sp}{\Scal_{\text{p}}}
\newcommand{\betap}{\betav_{\text{p}}}
\newcommand{\gammap}{\gamma_{\text{p}}}

\newcommand{\Sd}{\Scal_{\text{d}}}
\newcommand{\betad}{\betav_{\text{d}}}
\newcommand{\gammad}{\gamma_{\text{d}}}

\newcommand{\Se}{\Scal_{\text{e}}}
\newcommand{\betae}{\betav_{\text{e}}}
\newcommand{\gammae}{\gamma_{\text{e}}}

\newcommand{\wstar}{\wv_{\star}}
\newcommand{\xstar}{\xv_{\star}}
\newcommand{\zstar}{\zv_{\star}}

\newcommand{\Ee}{\Ecal_{\text{e}}}
\newcommand{\sOPT}{\hat{\sv}_{\text{OPT}}}
\newcommand{\se}{\mathrm{csign}}
\newcommand{\equipair}{(\Ee(\lambda),\se(\lambda))}

\newcommand{\Estar}{\Ecal_{\star}}
\newcommand{\sstar}{\sv_{\star}}
\newcommand{\Ebar}{\bar{\Ecal}}
\newcommand{\sbar}{\bar{\sv}}

\newcommand{\Spoly}{\hat{\Scal}_{\text{EQ-NQ}}}
\newcommand{\wEQ}{\hat{\wv}_{\text{EQ}}}

\newcommand{\Deltab}{\bm{\Delta}_{b}}
\newcommand{\Deltalmd}{{\Delta}_{\lambda}}

\newcommand{\lselect}{l_{\text{EN}}}
\newcommand{\rselect}{r_{\text{EN}}}
\newcommand{\Iselect}{\hat{I}_{\text{EN}}}

\newcommand{\Vo}{\Vcal_{\star}}
\newcommand{\Eo}{\Ecal^{\star}}
\newcommand{\so}{\sv^{\star}}

\newcommand{\Io}{I_{\star}}
\newcommand{\barIo}{\bar{I}_{\star}}

\newcommand{\Ijoin}{\Ical_{\text{insert}}}
\newcommand{\Idel}{\Ical_{\text{delete}}}

\newcommand{\la}{l^{\text{a}}}
\newcommand{\lb}{l^{\text{b}}}

\newcommand{\ra}{r^{\text{a}}}
\newcommand{\rb}{r^{\text{b}}}

\newcommand{\Ia}{I^{\text{a}}}
\newcommand{\Ib}{I^{\text{b}}}

\newcommand{\sadj}{\sv_{\text{adj}}}
\newcommand{\tadj}{t_{\text{adj}}}
\newcommand{\Tzone}{\hat{T}_{\text{EN}}}

\newcommand{\txtred}[1]{\textcolor{red}{#1}}

\title{\huge Piecewise Linearity of Min-Norm Solution Map of a Nonconvexly Regularized Convex Sparse Model}

\author{Yi Zhang$^{\orcidlink{0000-0003-4870-7884}}$, \IEEEmembership{Student Member, IEEE} and Isao Yamada$^{\orcidlink{0000-0002-6563-7526}}$, \IEEEmembership{Fellow, IEEE}
\thanks{This work was supported in part by JSPS Grants-in-Aid (21J22393), JSPS Grants-in-Aid (19H04134) and by JST SICORP (JPMJSC20C6). A shorter version of this paper has been accepted by ICASSP 2024. \textit{(Corresponding author: Yi Zhang)}}
\thanks{The authors are with the Department of Information and Communications Engineering, Tokyo Institute of Technology, Meguro-ku, Tokyo 152-8550, Japan (email: yizhang.ch.2015@gmail.com; isao@sp.ce.titech.ac.jp)}}

\maketitle

\begin{abstract}
It is well known that the minimum $\ell_2$-norm solution of the convex LASSO model, say $\xstar$, is a continuous piecewise linear function of the regularization parameter $\lambda$, and its signed sparsity pattern is constant within each linear piece. The current study is an extension of this classic result, proving that the aforementioned properties extend to the min-norm solution map $\xstar(\yv,\lambda)$, where $\yv$ is the observed signal, for a generalization of LASSO termed the scaled generalized minimax concave (sGMC) model. The sGMC model adopts a nonconvex debiased variant of the $\ell_1$-norm as sparse regularizer, but its objective function is overall-convex. Based on the geometric properties of $\xstar(\yv,\lambda)$, we propose an extension of the least angle regression (LARS) algorithm, which iteratively computes the closed-form expression of $\xstar(\yv,\lambda)$ in each linear zone. Under suitable conditions, the proposed algorithm provably obtains the whole solution map $\xstar(\yv,\lambda)$ within finite iterations. Notably, our proof techniques for establishing continuity and piecewise linearity of $\xstar(\yv,\lambda)$ are novel, and they lead to two side contributions: (a) our proofs establish continuity of the sGMC solution set as a set-valued mapping of $(\yv,\lambda)$; (b) to prove piecewise linearity and piecewise constant sparsity pattern of $\xstar(\yv,\lambda)$, we do not require any assumption that previous work relies on (whereas to prove some additional properties of $\xstar(\yv,\lambda)$, we use a different set of assumptions from previous work).
\end{abstract}


\tableofcontents

\section{Introduction}
\label{sec:intro}
In the past two decades, the sparse least-squares problems have attracted significant attention in signal processing and machine learning \cite{candes2006,donoho2006,eldar2012,qaisar2013,hastie2015},
which are usually formulated as the following optimization program:
\begin{equation}\label{eq:sparse_least_squares}
\underset{\xv\in\Real^n}{\text{minimize}}\;\;J(\xv)\coloneqq \frac{1}{2}\norm{\yv-\Am\xv}^2_2+\lambda\, \Psi(\xv),
\end{equation}
where $\yv\in\Real^m$ is the measured signal, $\Am\in\Real^{m\times n}$ is the sensing matrix, $\Psi:\Real^n\to\Real$ is a sparseness-promoting penalty which approximates the $\ell_0$ pseudo-norm (i.e., the cardinality of nonzero components in $\xv$), and $\lambda> 0$ is the regularization parameter which trades off between sparsity and data fidelity.

A prominent choice for $\Psi$ is the $\ell_1$-norm, which is known as the tightest convex approximation of the $\ell_0$ pseudo-norm. The resultant convex program (\ref{eq:sparse_least_squares}), also known as the LASSO \cite{tibshirani1996} (or basis
pursuit denoising \cite{chen2001}) problem, has been widely used \cite{eldar2012,hastie2015} for several favorable properties: on one hand, the LASSO cost function is convex, which allows us to develop efficient and scalable algorithms \cite{combettes2005,beck2009} that are provably convergent to a global minimizer of (\ref{eq:sparse_least_squares}); on the other hand, LASSO has many elegant geometric properties that facilitates its analysis, computation and application \cite{tibshirani1996,osborne2000,efron2004,hastie2007,zou2007,mairal2012,tibshirani2012,tibshirani2013}. These properties mainly include:
\begin{enumerate}[label=(\alph*)]
\item \label{it:property:unique} the LASSO solution is almost surely unique and has at most $\min\braces{m,n}$ nonzero components under a mild condition, which leads to a desired sparse solution in a high-dimensional setting (i.e., $m\ll n$),

\item \label{it:property:path} the minimum $\ell_2$-norm LASSO regularization path\footnote{A regularization path (or solution path) of a sparse least-squares problem (\ref{eq:sparse_least_squares}) is the solution of (\ref{eq:sparse_least_squares}) as a function of the regularization parameter $\lambda$.} is continuous and piecewise linear in $\lambda$ with its signed sparsity pattern being constant within each linear piece, which provides important clues for analyzing its variable selection behavior,

\item \label{it:property:LARS} the whole minimum $\ell_2$-norm LASSO regularization path can be computed by the least angle regression (LARS \cite{efron2004}) algorithm within finite steps, which eases tuning of the regularization parameter $\lambda$.
\end{enumerate}
We refer the interested readers to the seminal work of R. J. Tibshirani \cite{tibshirani2013} for a brilliant review of these properties.

Nevertheless, since the $\ell_1$-norm is coercive (i.e., $\norm{\xv}_1$ goes to infinity as $\norm{\xv}_2\to+\infty$), the LASSO solution tends to underestimate large-amplitude components, which leads to estimation bias. Although nonconvex regularizers such as SCAD \cite{fan2001} and MC \cite{zhang2010} penalties have been proposed to yield less biased sparse regularization, they usually sacrifice the convexity of (\ref{eq:sparse_least_squares}) and lose the favorable properties of LASSO.

To resolve this dilemma, the generalized minimax concave (GMC \cite{selesnick2017}) penalty has been proposed to yield less biased convex sparse regularization. Its formulation is as follows:
\begin{equation}\label{eq:GMC_penalty}
	\Psi_{\text{GMC}}(\xv;\Bm)\coloneqq \norm{\xv}_1-\min_{\zv\in\Real^n}\parens{\norm{\zv}_1+ \frac{1}{2}\norm{\Bm(\xv-\zv)}^2_2},
\end{equation}
where the subtrahend part is a debiasing function which cancels out the overpenalization effect caused by the $\ell_1$-norm, and the steering matrix $\Bm\in\Real^{q\times n}$ therein is a shape controlling parameter which adjusts the shape of $\Psi_{\textrm{GMC}}$ to achieve overall-convexity of (\ref{eq:sparse_least_squares}). Especially, if $\Bm^{\Tr}\Bm=\mathbf{O}_{n\times n}$, then the GMC penalty reduces to the $\ell_1$-norm; if $\Bm^{\Tr}\Bm$ is diagonal, then the GMC penalty is separable and reduces to a weighted sum of the MC penalties. Hence the GMC penalty is a nonseparable multidimensional extension of the MC penalty (and the LASSO penalty), which accounts for its name.

Notably, while the GMC penalty itself is nonconvex in general, it has been proven in \cite[Thm. 1]{selesnick2017} that as long as $\Bm$ satisfies the following linear matrix inequality
\begin{equation}\label{eq:overall_convexity_condition}
	\Am^{\Tr}\Am\succeq \lambda\Bm^{\Tr}\Bm,
\end{equation}
then the cost function of the resultant regularization problem (\ref{eq:sparse_least_squares}) is convex, which makes an essential difference from earlier nonconvex sparse regularizers \cite{fan2001,zhang2010}. Moreover, in contrast to the conventional LASSO model which solely convexifies the regularizer, the GMC model permits the regularizer to be nonconvex and convexifies the cost function as a whole. Hence the GMC model is a better convex approximation of the $\ell_0$-regularization model than LASSO and can possibly serve as a less biased surrogate of LASSO in certain applications.

Considerable efforts have been made to facilitate the application of the GMC model. For example, several extensions of GMC have been proposed to handle more complicated scenarios than sparse regression \cite{lanza2019,abe2020,alshabili2021,lanza2021,yata2022,liu2023,zhang2023}. Moreover, exploiting the overall-convexity, a number of efficient and scalable algorithms have been developed for finding a global minimizer of the GMC model \cite{selesnick2017} or its extensions \cite{lanza2019,abe2020,alshabili2021,zhang2021,zhang2022,yata2022,liu2023,zhang2023}. Nevertheless, up to now, little is known about the solution-set geometry and regularization path of the GMC model, which constitutes the motivation of the current study. 

In this paper, we consider the GMC model equipped with a specially designed\footnote{The steering matrix design (\ref{eq:B_of_sGMC}) was firstly proposed in the pioneering work of Selesnick \cite[Eq. (48)]{selesnick2017} that proposes the GMC model, and has been widely used in previous work \cite{alshabili2021,lanza2021,zhang2023,liu2023}. We note that although \cite[Eq. (48)]{selesnick2017} permits $\rho$ to be 1, Prop. \ref{prop:properties_of_bar_f} (Appx. \ref{proof:geometry_and_uniqueness}) in the present paper requires $\rho$ to be smaller than 1, which is a preliminary proposition for proving many subsequent results. Hence we assume $\rho<1$ in this study.} $\Bm$ that can automatically satisfy (\ref{eq:overall_convexity_condition}):
\begin{equation}\label{eq:B_of_sGMC}
	\Bm\coloneqq \sqrt{\rho/\lambda}\Am,\;\; \text{where }\rho\in[0,1).
\end{equation}
We formally refer to the GMC model using the steering matrix design (\ref{eq:B_of_sGMC}) as the scaled GMC (sGMC) model, and prove that: 
\begin{enumerate}[label=\arabic*)]
	\item despite the presence of a nonconvex debiasing term in the sGMC penalty, the sGMC model preserves the aforementioned properties \ref{it:property:unique} \ref{it:property:path} \ref{it:property:LARS} of LASSO, hence can reduce the estimation bias of LASSO without losing its geometric advantages,
	
	\item the piecewise linearity and piecewise constant sparsity pattern of the sGMC regularization path indeed extends to the min-norm sGMC solution map with respect to $(\yv,\lambda)$ in (\ref{eq:sparse_least_squares}), and the whole solution map can be computed by an extension of the LARS algorithm within finite iterations.
\end{enumerate}
Interestingly, our study implies that for a fixed linear system, the min-norm sGMC/LASSO solution map has an exact neural network representation, which draws a connection to deep learning methods for sparse recovery (e.g., \cite{gregor2010}).

\subsection{Our Contributions}\label{sec:contributions}

Our main contribution is threefold:
\begin{enumerate}[label=\arabic*)]
	\item We conduct a min-max analysis of a reformulation of the sGMC model, whereby we show that the sGMC model has a structurally similar optimality condition to LASSO (i.e., (\ref{eq:OPT}) in Sec. \ref{sec:optimality_condition}) and enjoys similar solution-set geometry (i.e., the property \ref{it:property:unique} of LASSO),
	
	\item We prove that continuity, piecewise linearity and piecewise constant sparsity pattern of the LASSO regularization path extend to the min-norm sGMC solution map, say $\xstar$, with respect to $(\yv,\lambda)$. This extension is achieved by two novel proof techniques\footnote{Previous work \cite{efron2004,tibshirani2013} proved conditional correctness of the LARS algorithm by induction, whereby the aforementioned property \ref{it:property:path} of LASSO can be implied by properties of the LARS iterative step. In contrast, our proof for extending the property \ref{it:property:path} to the min-norm sGMC solution map is a purely geometric analysis. Notably, piecewise linearity of the solution map naturally implies continuity within each linear zone, but continuity on the boundary is troublesome and requires an independent proof, which motivates us to study set-valued continuity of the sGMC model.}:
	\begin{itemize}
		\item to prove continuity, we study continuity of the sGMC solution set as set-valued mapping of $(\yv,\lambda)$,
		
		\item to prove piecewise linearity, we study the conditions that lead to an equality between $\xstar(\yv,\lambda)$ and elements from a finite set of affine functions of $(\yv,\lambda)$.
	\end{itemize}
	
	\item We propose an extension of the LARS algorithm, termed E-LARS, for iteratively computing the closed-form expression of $\xstar(\yv,\lambda)$ within each linear zone. Under suitable conditions, we prove that E-LARS obtains the whole solution map of $\xstar(\yv,\lambda)$ within finite iterations.
\end{enumerate}

It is worth noting that since the optimality condition (\ref{eq:OPT}) of the sGMC model structurally resembles the KKT condition of LASSO, with some efforts, it is possible to adapt the analysis in previous work \cite{efron2004,tibshirani2013} to the sGMC model, whereby one can extend the properties \ref{it:property:path} \ref{it:property:LARS} of LASSO to sGMC under certain assumptions (cf. Sec. \ref{sec:comments_on_tib13}). Nevertheless, such adaptation cannot extend these results to the min-norm sGMC solution map $\xstar(\yv,\lambda)$. In addition, our novel proofs for the second contribution above yield two side contributions:
\begin{enumerate}[label=\arabic*)]
	\item Our proof for continuity of $\xstar(\yv,\lambda)$ establishes set-valued continuity (i.e., \textit{outer} and \textit{inner semicontinuity} \cite{rockafellar2009}) of the sGMC/LASSO solution set with respect to $(\yv,\lambda)\in\Real^m\times (0,+\infty)$ in (\ref{eq:sparse_least_squares}).
	
	\item Our proof for piecewise linearity and piecewise constant sparsity pattern of $\xstar(\yv,\lambda)$ does not require any assumption that previous work relies on. (However, to prove some additional property of $\xstar(\yv,\lambda)$ and correctness of E-LARS, we require a different set of assumptions from previous work; see Thm. \ref{thm:correspondence} \ref{it:corres:GP} and Thm. \ref{thm:E_LARS} for details).
\end{enumerate}
 
\subsection{Related Works}
We first review some notable previous studies on the piecewise linear regularization path of LASSO or its variants. Osborne et al. \cite{osborne2000LARS} firstly proposed a homotopy method to compute the LASSO regularization path, which is essentially the same idea as LARS \cite{efron2004}, but the proof of correctness of the algorithm was not given explicitly. Efron et al. \cite{efron2004} developed the well known LARS algorithm, where continuity and piecewise linearity of the regularization path are established as a by-product of the conditional correctness of LARS. However, the proof in \cite{efron2004} requires two assumptions\footnote{The original proof in \cite{efron2004} assumed that the sensing matrix has full column rank (cf. \cite[p. 413, second paragraph]{efron2004}) to ensure solution uniqueness, which can be relaxed to the mild ``general positioning" condition \cite{dossal2012}. The ``one-at-a-time" assumption is not well-understood theoretically, but it holds with high probability in random numerical experiments. It is not known whether these two assumptions imply each other.}: solution uniqueness of LASSO and the ``one-at-a-time" assumption \cite[Thm. 1]{efron2004}. Rosset and Zhu \cite{rosset2007} extended the result of LARS to a more general $\ell_1$-regularization model with the family of ``almost quadratic" loss functions. They developed a piecewise linear solution path algorithm, but the derivation implicitly used the same assumptions as \cite{efron2004}. R. J. Tibshirani and Taylor \cite{tibshirani2011,tibshirani2011phd} extended many results of LARS to the generalized LASSO model. R. J. Tibshirani also improved the proof of the conventional LARS-LASSO algorithm and provided many new results along with a brilliant review of previous studies in \cite{tibshirani2013}. The proofs in \cite{tibshirani2011,tibshirani2011phd,tibshirani2013} do not require the solution to be unique, but implicitly used the ``one-at-a-time" assumption and another assumption (cf. Sec. \ref{sec:comments_on_tib13} for details). Recently, Berk et al. \cite{berk2023} proved the local Lipschitz continuity of the LASSO solution with respect to $(\yv,\lambda)$, assuming the solution to be unique.

There also have been some studies on the regularization path of more complicated extensions of LASSO, which are not piecewise linear. However, in this case, even proving the continuity of the regularization path is difficult. Xiao et al. \cite{xiao2015} proposed a regularization path algorithm for the group LASSO model, which involves solving an ODE system in every iteration. They established the continuity of the unique regularization path assuming the loss function to be strictly convex and twice differentiable. Yukawa and Amari \cite{yukawa2016} studied the path of critical points of the nonconvex $\ell_p$-regularized ($0<p<1$) least squares problem and confirmed its discontinuity. Mishkin and Pilanci \cite{mishkin2022} established almost everywhere continuity of the regularization path for group LASSO, which can be improved to everywhere continuity under the solution uniqueness assumption. The same authors \cite{mishkin2023} also studied solution-set geometry and regularization path of a constrained group LASSO problem. It should be noted that the tools therein for studying solution uniqueness and set-valued continuity are the same as that in this paper (but the proof techniques are not), and outer-semicontinuity of the solution set is proved.

Another notable work that is closely related to this paper is conducted by Liu et al. \cite{liu2023}, where the authors considered group GMC, an extension of the GMC model for group sparse linear regression, and proved the following results: 
\begin{itemize}
	\item if the solution of the group GMC model is unique, then its regularization path is continuous\footnote{In \cite[Thm. 2.1]{liu2023}, the authors originally claimed that solution uniqueness holds, but one can verify that their proof is flawed, by noticing LASSO as a special instance of group GMC whose solution is not unique in general \cite{tibshirani2013}.},
	
	\item under the restricted eigenvalue condition (which is commonly used in statistical analysis of LASSO \cite{hastie2015}), the group GMC estimator has a similar error bound to that of group LASSO, where the steering matrix $\Bm$ only affects the constant factor of the bound.
\end{itemize}
However, we should note that since the constant factor involving $\Bm$ does not have an explicit expression, it is difficult to quantitatively study the effect of $\Bm$ using the error bound in \cite{liu2023}. Therefore, the statistical performance gain of the GMC model over LASSO remains to be an open problem.

\subsection{Paper Organization}

The remainder of this paper is organized as follows. Sec. \ref{sec:preliminaries} introduces the notations used in this paper. In Sec. \ref{sec:static_geometry}, we present some useful reformulations of the sGMC model, its optimality condition and static geometric properties (where problem parameters are fixed). In Sec. \ref{sec:min_norm_solution_map_analysis}, we consider $(\yv,\lambda)$ as varying parameters, and present results regarding continuity, piecewise linearity and piecewise constant sparsity pattern of the min-norm sGMC solution map $\xstar(\yv,\lambda)$. In Sec. \ref{sec:E_LARS}, we propose the E-LARS algorithm, and we analyze its conditional correctness and complexity. Sec. \ref{sec:proof:continuity} and Sec. \ref{sec:proof:correspondence} respectively present our novel proofs for continuity and for piecewise linearity/piecewise constant sparsity pattern of $\xstar(\yv,\lambda)$. Sec. \ref{sec:conclusion} closes this paper with some concluding remarks.

\section{Notations}
\label{sec:preliminaries}

\subsection{Numbers, Sets and Functions}
Let $\Natural,\Real,\Real_+,\Real_{++}$ respectively be the sets of nonnegative integers, real numbers, nonnegative real numbers and positive real numbers. For a number set $I\subset\Real$, we denote $\sup(I)$ and $\inf(I)$ respectively as its supremum and infimum; especially, when $I=\emptyset$, we adopt the convention $\sup(\emptyset)\coloneqq -\infty$ and $\inf(\emptyset)\coloneqq +\infty$. For the $n$-dimensional Euclidean space $\Real^n$, $\norm{\cdot}_{p}\;(p\geq 1)$ denotes the $\ell_p$-norm in $\Real^n$. For a set $S\subset\Real^n$, $\interior(S)$ and $\closure(S)$ are respectively its interior and closure.

For an extended real-valued function $f:\Real^n\to \Real\cup\{+\infty\}$, the domain of $f$ is $\dom{f}\coloneqq\braces{\xv\in\Real^n\mid f(\xv)<+\infty}$. For a convex function $f:\Real^n\to\Real\cup\braces{+\infty}$, the subdifferential of $f$ at $\xv\in\Real^n$ is defined as
\begin{equation*}
	\partial f(\xv)\coloneqq\{\uv\in\Real^n\mid(\forall \zv\in\Real^n)\;\uv^\Tr(\zv-\xv)+f(\xv)\leq f(\zv)\},
\end{equation*}
if $\uv\in\partial f(\xv)$, then $\uv$ is called a subgradient of $f$ at $\xv$. Especially, if $f$ is differentiable at $\xv\in\Real^n$, we denote $\nabla f(\xv)$ as the gradient of $f$ at $\xv$, and the subdifferential of $f$ is given by $\partial f(\xv)=\braces{\nabla f(\xv)}$. We denote the set of all proper, lower semicontinuous\footnote{For $f:\Real^n\to\Real\cup\{+\infty\}$, $f$ is referred to as \textit{proper} if $\dom{f}\neq\emptyset$, and \textit{lower semicontinuous} if every lower level set of it is closed.} convex functions from $\Real^n$ to $\Real\cup\{+\infty\}$ by $\Gamma_0(\Real^n)$. For two sets $X$ and $Y$, we adopt the notation $\Mcal:X\rightrightarrows Y$ to denote a set-valued mapping (a.k.a. multifunction) from $X$ to the power set of $Y$, i.e., $\Mcal$ maps a point in $X$ to a (possibly empty) subset of $Y$.

\subsection{Vectors and Matrices}

For a vector $\xv\in\Real^n$, we use $x_i$ or $[\xv]_i$ to denote its $i$th component, and respectively define its support as $\supp(\xv)\coloneqq \braces{i\in\braces{1,\dots,n}\;\middle\vert\; x_i\neq 0}$ and define its signs vector as $\sign(\xv)\coloneqq [\sign(x_1),\sign(x_2),\dots,\sign(x_n)]^{\Tr}\in\braces{+1,0,-1}^n$. For a matrix $\Am\in\Real^{m\times n}$, $\Am^{\Tr}\in\Real^{n\times m}$ is its transpose, $\Am^{\dagger}\in\Real^{n\times m}$ is its Moore-Penrose pseudoinverse, $\Null(\Am)$ and $\Col(\Am)$ respectively denote its null space and
column space. For $\xv,\yv\in\Real^n$, $\xv\geq \yv$ means that $x_i\geq y_i$ for every $i\in\{1,\dots,n\}$. For $\Am,\Bm\in\Real^{n\times n}$, $\Am\succeq \Bm$ means that $\Am-\Bm$ is positive semidefinite, $\Am\odot\Bm$ denotes their Hadamard product (a.k.a. entrywise product). The symbols $\mathbf{0}_n$, $\mathbf{1}_n$, $\mathbf{O}_{m\times n}$ and $\Imat_n$ are respectively the $n\times 1$ zero vector, the $n\times 1$ all-ones vector, the $m\times n$ zero matrix and the $n\times n$ identity matrix (when the ambient dimension is known, we may omit the subscripts). We denote the diagonal matrix with diagonal entries $d_1,\dots,d_n$ by $\mathrm{diag}\parens{d_1,\dots,d_n}$, and denote the block diagonal matrix with diagonal blocks $\Am_1\in\Real^{m_1\times n_1}$, $\dots$, $\Am_k\in\Real^{m_k\times n_k}$ by $\blkdiag(\Am_1,\dots,\Am_k)$. 

For a vector $\xv\coloneqq [x_1,\dots,x_n]^{\Tr}\in\Real^n$, a matrix $\Am\coloneqq [\av_1,\dots,\av_n]\in\Real^{m\times n}$ and an index set $\Ical=\braces{i_1,\dots,i_k}\subset \braces{1,\dots,n}$ with $i_1<i_2<\cdots<i_k$, we define $\xv_{\Ical}\coloneqq[\xv]_{\Ical}\coloneqq [x_{i_1},\dots, x_{i_k}]^{\Tr}$, $\Am_{\Ical}\coloneqq [\av_{i_1},\dots,\av_{i_k}]$ and $\neg\Ical\coloneqq \braces{1,\dots,n}\setminus\Ical$, respectively. In particular, we define $\Am_{\emptyset}\coloneqq \mathbf{0}_{m}$ and $\xv_{\emptyset}\coloneqq 0$ such that the slicing operation for the empty set does not lead to a void value, which avoids possible ambiguity in computations involving $\Am_{\emptyset}$ and $\xv_{\emptyset}$. The slicing operation will be repeatedly used in proving piecewise linearity of the minimum $\ell_2$-norm sGMC solution map (Sec. \ref{sec:min_norm_solution_map_analysis} and \ref{sec:proof:correspondence}).

Finally, we introduce the general-positioning condition \cite{dossal2012}, which is a mild condition widely used for studying solution uniqueness of LASSO and its variants.
\begin{definition}[{\cite[Def. 1.1]{dossal2012}}]\label{def:general_position}
	We say $\Am\coloneqq [\av_1,\dots,\av_n]\in\Real^{m\times n}$ has columns in \textit{general position} if for every $\Ical\subset\braces{1,\dots,n}$ satisfying $|\Ical|\leq m$ and $(s_j)_{j\in\Ical}\subset\braces{\pm 1}$, we have
	\begin{equation}\label{eq:general_position}
		(\forall i\in \Ical)\;\; s_i\av_i\not\in\textrm{aff}\braces{s_j\av_j}_{j\in\Ical\setminus \braces{i}},
	\end{equation}
	where for $\braces{\bv_1,\dots,\bv_n}\subset\Real^n$, $\textrm{aff}\braces{\bv_1,\dots,\bv_n}$ denotes the affine hull of $\braces{\bv_1,\dots,\bv_n}$.
\end{definition}
Notice that the logical complement of (\ref{eq:general_position}) leads to a set of $\Am$ with Lebesgue measure zero, one can deduce that if the entries of $\Am$ follow a continuous probability distribution on $\Real^{m\times n}$, then $\Am$ has columns in general position with probability one \cite[the paragraph after Def. 1.1]{dossal2012}.

\section{Static Solution-Set Geometry}
\label{sec:static_geometry}
In this section, we introduce reformulations of the sGMC model, its optimality condition and its static solution-set geometry (i.e., we consider the scenario where $(\yv,\lambda)$ is fixed). This lays the foundation for our subsequent analysis, and addresses the first major contribution stated in Sec. \ref{sec:contributions}.

\subsection{Three Equivalent Formulations of the sGMC Model and the Optimality Condition}\label{sec:optimality_condition}

To begin with, we introduce a reformulation of the sGMC model as a min-max problem and its associated primal and dual minimization problems.

Substituting (\ref{eq:B_of_sGMC}), (\ref{eq:GMC_penalty}) into (\ref{eq:sparse_least_squares}) and conducting some basic transformations, we can rewrite the sGMC model as the following min-max problem:
\begin{equation}\label{eq:sGMC_model}
	\min_{\xv\in\Real^n}\max_{\zv\in\Real^n}\;\; \frac{1}{2}\norm{\yv-\Am\xv}^2_2+\lambda\norm{\xv}_1-\lambda\norm{\zv}_1-\frac{\rho}{2}\norm{\Am\xv-\Am\zv}^2_2.
\end{equation}
For convenience of subsequent analysis, here we introduce an additional term\footnote{This additional term does not have a concrete meaning; it is only used for easing the proof of piecewise linearity of the min-norm sGMC solution map with respect to $(\yv,\lambda)$.} $\sqrt{\rho}\rv^\Tr\Am \zv$ in the objective function of (\ref{eq:sGMC_model}), i.e., we consider the following function $G(\xv,\zv)$:
\begin{equation}\label{eq:def_of_G}
	\begin{aligned}
		G(\xv,\zv) \coloneqq  \frac{1}{2}\norm{\yv-\Am\xv}^2_2+\lambda\norm{\xv}_1-\lambda\norm{\zv}_1  -\frac{\rho}{2}\norm{\Am\xv-\Am\zv}^2_2 +\sqrt{\rho}\rv^\Tr\Am \zv.
	\end{aligned}	
\end{equation}
where $\rv\in\Real^m$ is an auxiliary problem parameter. Throughout this paper, we refer to the convex-concave min-max problem
\begin{equation}\label{eq:sGMC_as_minmax}
	\min_{\xv\in\Real^n}\max_{\zv\in\Real^n}\;\; G(\xv,\zv)
\end{equation}
as the min-max form of the sGMC model, and refer to the convex program
\begin{equation}\label{eq:primal_sGMC}
	\underset{\zv\in\Real^n}{\text{minimize}}\;\; {J}_{\text{sGMC}}(\zv)\coloneqq {\max_{\zv\in\Real^n}G(\xv,\zv)}
\end{equation}
as the primal form of the sGMC model.

Moreover, since $G(\cdot,\zv)$ (resp. $-G(\xv,\cdot)$) is convex, continuous and coercive for every $\zv\in\Real^n$ (resp. $\xv\in\Real^n$), we have the following minimax equality from an extension of Sion's minimax theorem by Hartung \cite[Thm. 1]{hartung1982}:
\begin{equation}\label{eq:minimax_equality}
	\min_{\xv\in\Real^n} \max_{\zv\in\Real^n} G(\xv,\zv)=\max_{\zv\in\Real^n} \min_{\xv\in\Real^n} G(\xv,\zv),
\end{equation}
which implies the existence of a saddle point of (\ref{eq:sGMC_as_minmax}). Then from the RHS of the minimax equality (\ref{eq:minimax_equality}), we can further define the dual form of the sGMC model:
\begin{equation}\label{eq:dual_sGMC}
	\underset{\zv\in\Real^n}{\text{minimize}}\;\; \bar{J}_{\text{sGMC}}(\zv)\coloneqq \parens{-\min_{\xv\in\Real^n}G(\xv,\zv)}.
\end{equation}

Next we introduce some notations for denoting the solution sets of the three formulations of the sGMC model.

\begin{definition}\label{def:solution_sets}
	We respectively define the sets of minimizers of (\ref{eq:primal_sGMC}) and (\ref{eq:dual_sGMC}) as the \textit{(primal) solution set} (denoted by $\Sp$) and \textit{dual solution set} (denoted by $\Sd$) of the sGMC model. In addition, we define the set of vectors $\wv\coloneqq \begin{bmatrix}
		\xv^\Tr & \zv^\Tr
	\end{bmatrix}^\Tr$ with $(\xv,\zv)\in\Real^n\times\Real^n$ being a saddle point of (\ref{eq:sGMC_as_minmax}) as the \textit{extended solution set} (denoted by $\Se$) of the sGMC model. 
\end{definition}

Let $\bv\coloneqq \begin{bmatrix}
\yv^\Tr & \rv^\Tr
\end{bmatrix}^\Tr\in\Real^{2m}$. Then $\Sp$, $\Sd$ and $\Se$ are all set-valued mappings of the following problem parameters:
\begin{equation}
	\paragroup\in\paraspace=:\Pscr,
\end{equation}
where $\Pscr$ is the space of suitable problem parameters. Hereafter, when we want to emphasize such dependence on certain problem parameter (e.g., on $(\bv,\lambda)$), we may explicitly write $\Sp,\Sd,\Se$ or any associated quantity of the sGMC model as a mapping of the parameters of interest (e.g., $\Sp(\bv,\lambda)$).

It should be noted that if $\rho=0$ and $\rv=\mathbf{0}_m$, then the primal form (\ref{eq:primal_sGMC}) of sGMC model reduces to the LASSO model:
\begin{equation}\label{eq:LASSO}
	\underset{\xv\in\Real^n}{\text{minimize}}\;\;J_{\text{LASSO}}(\xv)\coloneqq \frac{1}{2}\norm{\yv-\Am\xv}^2_2+\lambda \norm{\xv}_1.
\end{equation}
Hence every theoretical result of the sGMC model introduced in this paper naturally applies to the LASSO model. We denote the solution set of (\ref{eq:LASSO}) by $\SLA$, or $\SLA(\Am,\yv,\lambda)$ to emphasize its dependence on problem parameters.

Finally, we present the optimality condition of an extended solution, which will be a main subject of study in the subsequent min-norm solution map analysis.

\begin{lemma}\label{lemma:OPT}
	Given $\paragroup\in\Pscr$, $\wv\in\Real^{2n}$ is an extended solution in $\Se$ (i.e., $\wv$ is the concatenation of some primal solution $\xv\in\Sp$ and dual solution $\zv\in\Sd$) if and only if the following inclusion holds for every $i\in\{1,2,\dots,2n\}$:
	\begin{equation}
		\tag{OPT} \cv_i^{\Tr}(\bv-\Dm\Cm\wv)\in\begin{cases}
			\braces{\lambda\,\sign([\wv]_i)}, & \textrm{if }[\wv]_i\neq 0,\\
			[-\lambda,\lambda], & \textrm{if }[\wv]_i=0,
		\end{cases}\label{eq:OPT}
	\end{equation}
	where $\Cm \coloneqq \blkdiag(\Am,\sqrt{\rho} \Am)$,
	\begin{align}
		\Dm &\coloneqq \begin{bmatrix}(1-\rho)\Imat_m & \sqrt{\rho}\Imat_m\\ -\sqrt{\rho}\Imat_m & \Imat_m\end{bmatrix},\label{eq:matrix_D}
	\end{align}
	and $\cv_i$ is the $i$th column vector of $\Cm$.
\end{lemma}
\begin{IEEEproof}
	(\ref{eq:OPT}) can be derived by transforming the saddle point condition of (\ref{eq:sGMC_as_minmax}). See Appendix \ref{app:derive_OPT} for details.
\end{IEEEproof}

One may realize that (\ref{eq:OPT}) is structurally similar to the KKT condition of the LASSO model (\ref{eq:LASSO}):
\begin{equation}\label{eq:KKT_LASSO}
	\tag{KKT} \av_i^\Tr(\yv-\Am\xv)\in\begin{cases}
		\braces{\lambda\,\sign([\xv]_i)}, & \textrm{if }[\xv]_i\neq 0,\\
		[-\lambda,\lambda], & \textrm{if }[\xv]_i=0,
	\end{cases}
\end{equation}
where $i\in\braces{1,2,\dots,n}$, $\av_i$ is the $i$th column vector of $\Am$. 

Particularly, if $\Dm$ in (\ref{eq:OPT}) is an arbitrary positive semidefinite matrix (without considering the constraint (\ref{eq:matrix_D})), then (\ref{eq:OPT}) can be rewritten into the same form of (\ref{eq:KKT_LASSO}) by some simple transformations. However, due to the constraint (\ref{eq:matrix_D}), $\Dm$ is not positive semidefinite for every $0<\rho<1$. Indeed, while the KKT conditions of constrained convex programs are also derived from saddle point conditions (of their Lagrangians), the saddle point condition (\ref{eq:OPT}) cannot be written as the KKT condition of a simple convex program\footnote{This can be shown by noticing that a Lagrangian is always linear in its dual variable (Lagrangian multipliers), whereas the objective function $G(\xv,\zv)$ of (\ref{eq:sGMC_as_minmax}) is not. We want to thank Professor R. J. Tibshirani for a helpful discussion that inspires the first author to notice this difference.}.

\subsection{Shape, Uniqueness and Sparseness of the Solution Sets}

In the sequel, we introduce geometric properties of $\Sp$, $\Sd$ and $\Se$ with a fixed group of problem parameters $\paragroup$, which extend the property \ref{it:property:unique} of LASSO introduced in Section \ref{sec:intro} to the sGMC model, and serve as preliminary results for the subsequent min-norm solution map analysis.

\begin{theorem}\label{thm:shape_uniqueness_sparseness}
	Given $\paragroup\in\Pscr$, then the solution sets $\Sp$, $\Sd$ and $\Se$ of the sGMC model are all nonempty, closed, convex and bounded sets. Moreover, the following holds:
	\begin{enumerate}[label=(\alph*)]
		\item \label{it:shape:linear_fit} Every $\xv\in\Sp$ (resp. $\zv\in\Sd$) gives the same linear fit (i.e., multiplication with $\Am$) and the same $\ell_1$-norm, i.e.,
	\begin{align*}
		\parens{\forall \xv_1,\xv_2\in\Sp}\;\; &\Am\xv_1=\Am\xv_2,\;\;  \norm{\xv_1}_1=\norm{\xv_2}_1,\\
		\parens{\forall \zv_1,\zv_2\in\Sd}\;\; &\Am\zv_1=\Am\zv_2,\;\;  \norm{\zv_1}_1=\norm{\zv_2}_1.
	\end{align*}
		
		\item \label{it:shape:unique} If $\Am$ has columns in general position, then the primal solution $\xv\in\Sp$ and the dual solution $\zv\in\Sd$ are both unique with at most $\min\braces{m,n}$ nonzero components.
		
		\item \label{it:shape:bound} We have the following upper bound for the common $\ell_1$-norm of every $\xv\in\Sp$ and every $\zv\in\Sd$:
		\begin{equation}\label{eq:l1_bound}
			\max\braces{\norm{\xv}_1,\norm{\zv}_1}\leq \frac{\norm{\yv}_2^2}{2\lambda(1-\rho)}+\frac{\norm{\rv}^2_2}{2\lambda}.
		\end{equation}
		
		\item \label{it:shape:Se} The following relation holds among $\Sp$, $\Sd$, $\Se$ and the LASSO solution set mapping $\SLA(\cdot,\cdot,\cdot)$:
		\begin{align}
			\Se =& \Sp\times\Sd \label{eq:Se_SpSd}\\ =&\SLA\parens{\Am,\frac{\yv-\rho\betad}{1-\rho},\frac{\lambda}{1-\rho}}\times \SLA(\sqrt{\rho}\Am,\rv+\sqrt{\rho}\betap,\lambda),\label{eq:Se_variant_expression}
		\end{align}
		where $\betap$ and $\betad$ are respectively the common linear fits of every $\xv\in\Sp$ and $\zv\in\Sd$. Note that $\betap$ and $\betad$ are vector-valued mappings of $\paragroup\in\Pscr$.
	\end{enumerate}
\end{theorem}

We remark that the results \ref{it:shape:linear_fit} and \ref{it:shape:unique} are known for the LASSO case, whilst \ref{it:shape:bound} and \ref{it:shape:Se} are new results. In the following, we present some consequences of these results.

\begin{remark}[Consequences of Theorem \ref{thm:shape_uniqueness_sparseness}]
	\begin{enumerate}[label=\arabic*)]
		\item Applying the result \ref{it:shape:linear_fit} into (\ref{eq:sGMC_model}), one can verify that every $\xv\in\Sp$ gives the same data fidelity and the same sGMC penalty value, hence all sGMC solutions are equally good with respect to goodness of fit and the sGMC penalty.
		
		\item From the mildness of general positioning condition, the result \ref{it:shape:unique} yields almost sure uniqueness of the sGMC solution in a probabilistic setting, and further ensures its sparseness in the high dimensional setting where $m\ll n$.
		
		\item The result \ref{it:shape:Se} indicates that the interplay between $\Sp$ and $\Sd$ is simple: when constructing a saddle point $(\xv,\zv)$ of (\ref{eq:sGMC_as_minmax}), the selection of $\xv\in\Sp$ is independent of that of $\zv\in\Sd$. Moreover, (\ref{eq:Se_variant_expression}) elucidates an interesting fact: the sGMC solution set $\Sp$ can be rewritten as a nonlinearly parameterized LASSO solution set (albeit $\betad$ in the parameterization mapping is an implicit function).
	\end{enumerate}
\end{remark}

It is worth noting that in previous work \cite{tibshirani2013}, the result \ref{it:shape:linear_fit} for LASSO \cite[Lemma 1]{tibshirani2013} is not obtained by analyzing (\ref{eq:KKT_LASSO}), hence one cannot easily adapt the analysis therein to sGMC by exploiting the similarity between (\ref{eq:OPT}) and (\ref{eq:KKT_LASSO}), which further hinders the adaptation of the proof of result \ref{it:shape:unique}. In this paper, we circumvent this difficulty by another reformulation of the sGMC model as a $\ell_1$-regularization model with nonquadratic data-fidelity term; see the following proof sketch. The complete proof can be found in Appendix \ref{proof:geometry_and_uniqueness}.

\begin{remark}[Proof sketch of Theorem \ref{thm:shape_uniqueness_sparseness}]
	 The proof of Theorem \ref{thm:shape_uniqueness_sparseness} is developed by recognizing another reformation of (\ref{eq:primal_sGMC}) as an $\ell_1$-regularization model with an involved data fidelity term. The key point is to show that the involved data fidelity therein is differentiable, strictly convex and bounded from below as a function of $\Am\xv$, whereby we can apply results in \cite[Sec. 2.3]{tibshirani2013} to derive properties \ref{it:shape:linear_fit} and \ref{it:shape:unique} for $\Sp$. The properties \ref{it:shape:linear_fit} and \ref{it:shape:unique} for $\Sd$ are obtained via its interplay with $\Sp$ in (\ref{eq:OPT}). The results \ref{it:shape:bound} and \ref{it:shape:Se} can be derived from intermediate results of previous analysis.
\end{remark}

Finally, before proceeding to the analysis of the minimum $\ell_2$-norm sGMC solution map, we present a corollary ensuring the existence and uniqueness of this min-norm solution.

\begin{corollary}\label{cor:wstar}
	Given $\paragroup\in\Pscr$, then there exists uniquely a minimum $\ell_2$-norm element in $\Se$. Moreover, let this min-norm element be $\wstar\coloneqq \begin{bmatrix}\xstar^\Tr & \zstar^\Tr\end{bmatrix}^\Tr$ with $\xstar,\zstar\in\Real^n$, then $\xstar$ and $\zstar$ are respectively the unique minimum $\ell_2$-norm elements in $\Sp
	$ and $\Sd$. Notice that $\wstar$, $\xstar$ and $\zstar$ are all vector-valued mappings of $\paragroup\in\Pscr$.
\end{corollary}
\begin{IEEEproof}
	Since $\Se$ (resp., $\Sp$, $\Sd$) is nonempty, closed and convex, the projection of $\mathbf{0}$ onto $\Se$ (resp., $\Sp$, $\Sd$) exists and is unique from \cite[Thm. 3.16]{bauschke2017}, which guarantees the existence and uniqueness of $\wstar$ (resp. $\xstar$, $\zstar$). From (\ref{eq:Se_SpSd}), if $\wstar$ achieves minimum $\ell_2$-norm in $\Se$, its first and second half parts must achieve the minimum $\ell_2$-norm respectively in $\Sp$ and $\Sd$. Hence $\wstar\coloneqq \begin{bmatrix}\xstar^\Tr & \zstar^\Tr\end{bmatrix}^\Tr$.
\end{IEEEproof}

\section{Min-Norm Solution Map Analysis}
\label{sec:min_norm_solution_map_analysis}

In this section, we consider $(\bv,\lambda)$ as varying parameters, and establish continuity, piecewise linearity and piecewise constant sparsity pattern of the min-norm extended solution map $\wstar(\bv,\lambda)$. This section contains the most important theoretical findings and proof techniques of the current study, which addresses the second major contribution and the two side contributions stated in Sec. \ref{sec:contributions}.

\subsection{Continuity Properties of the Solution Sets, Min-Norm Solution Maps and Associated Quantities}

This subsection presents continuity properties of the sGMC model with respect to $(\bv,\lambda)$. All results here are obtained via a set-valued variational analysis on the extended solution set mapping $\Se(\bv,\lambda)$. Before presenting these results, we briefly introduce necessary notions for describing continuity of set-valued mappings\footnote{We note that there are multiple ways of defining set-valued continuity. The notions of outer and inner semicontinuity that we are using here are also known as \textit{closedness} and \textit{openness} of a set-valued mapping in \cite{hogan1973}. Please see \cite[Commentary in Sec. 5]{rockafellar2009} for other related notions.}. We start with the ``outer limit" and ``inner limit" of a sequence of sets, which are respectively analogies of the upper and lower limits of a real sequence. 

\begin{definition}[{\cite[Def. 4.1 and the paragraph after it]{rockafellar2009}}]\label{def:limit_of_sets}
	For a sequence of sets $(S_k)_{k\in\Natural}$ with $S_k\subset\Real^n$ for every $k\in\Natural$:
	\begin{enumerate}[label=(\alph*)]
		\item Its \textit{outer limit} $\limsup_k S_k$ is the set of all possible cluster points of every $(\xv_k)_{k\in\Natural}$ satisfying $\xv_k\in S_k$ for all $k$,
		
		\item Its \textit{inner limit} $\liminf_k S_k$ is the set of all possible limit points of every $(\xv_k)_{k\in\Natural}$ satisfying $\xv_k\in S_k$ for all $k$.
	\end{enumerate}
	
	From this definition, we naturally have $\liminf_{k}S_k\subset\limsup_{k}S_k$ for every $(S_k)_{k\in\Natural}$.
\end{definition}

Recall that a real-valued function $f$ is continuous at $x\in\Real$ if and only if $f(x)=\limsup_{z\to x}f(z)=\liminf_{z\to x}f(z)$, the continuity of set-valued mappings can be defined in a similar way; see the following definition.

\begin{definition}[{\cite[Def. 5.4]{rockafellar2009}}]\label{def:continuity_of_setvalued_mapping}
	For $D\subset\Real^m$, a point $\bar{\thetav}\in D$ and a set-valued mapping $\Mcal:D\rightrightarrows{\Real^n}$, we say:
	\begin{itemize}
		\item $\Mcal$ is \textit{outer semicontinuous} (osc) at $\bar{\thetav}$ relative to $D$ if for every $(\thetav_k)_{k\in\Natural}\subset D$ with $\lim_{k\to\infty}\thetav_k=\bar{\thetav}$, we have $\limsup_{k\to\infty} \Mcal(\thetav_k)\subset\Mcal(\bar{\thetav})$.
		\item $\Mcal$ is \textit{inner semicontinuous} (isc) at $\bar{\thetav}$ relative to $D$ if for every $(\thetav_k)_{k\in\Natural}\subset D$ with $\lim_{k\to\infty}\thetav_k=\bar{\thetav}$, we have $\liminf_{k\to\infty} \Mcal(\thetav_k)\supset\Mcal(\bar{\thetav})$.
		
		\item $\Mcal$ is called \textit{continuous} at $\bar{\thetav}$ relative to $D$ if it is both osc and isc at $\bar{\thetav}$ relative to $D$, i.e., if for every $(\thetav_k)_{k\in\Natural}\subset D$ with $\lim_{k\to\infty}\thetav_k=\bar{\thetav}$, we have $\limsup_{k\to\infty} \Mcal(\thetav_k)=\liminf_{k\to\infty}\Mcal(\thetav_k)=\Mcal(\bar{\thetav})$.
	\end{itemize}
	If $\Mcal$ is osc/isc/continuous at every point in $D$ relative to $D$, we say $\Mcal$ is osc/isc/continuous relative to $D$.
\end{definition}

Notice that single-valued mappings can be viewed as special set-valued mappings whose values are singletons. Then one can verify that in the context of single-valued mappings, Definition \ref{def:continuity_of_setvalued_mapping} coincides with the conventional definition of continuity of functionals (cf. \cite[Def. Thm. 4.2 and 4.6]{rudin1976}).

The following theorem summarizes continuity properties of the solution sets, min-norm solutions and other associated quantities of the sGMC model. For simplicity of presentation, all results are presented with respect to $\Se$ or quantities defined over $\Se$. However, each continuity property stated here naturally extends to $\Sp$, $\Sd$ and their associated quantities; one can readily transfer these properties via (\ref{eq:Se_SpSd}).

\begin{theorem}\label{thm:continuity}
	Let $\Se$ be the extended solution set (Def. \ref{def:solution_sets}) of the sGMC model. Define the following quantities:
	\begin{itemize}
		\item $\betae$ is the common ``extended" linear fit of every $\wv\in\Se$ (Thm. \ref{thm:shape_uniqueness_sparseness} \ref{it:shape:linear_fit} and (\ref{eq:Se_SpSd})), i.e., $\betae\coloneqq \Cm\wv$ with $\wv\in\Se$,
		
		\item $\gammae$ is the common $\ell_1$-norm of every $\wv\in\Se$ (Thm. \ref{thm:shape_uniqueness_sparseness} \ref{it:shape:linear_fit} and (\ref{eq:Se_SpSd})), i.e., $\gammae\coloneqq \norm{\wv}_1$ with $\wv\in\Se$,
		
		\item for given $\bar{\wv}\in\Real^{2n}$, $P_{\Se}\parens{\bar{\wv}}$ is the unique projection of $\bar{\wv}$ onto $\Se$. Notice that $P_{\Se}(\mathbf{0}_{2n})$ reproduces the minimum $\ell_2$-norm extended solution $\wstar$ (Cor. \ref{cor:wstar}).
	\end{itemize}
	
	Then we have the following continuity properties:
	\begin{enumerate}[label=(\alph*)]
		\item \label{it:continuity:para} Consider $\Se$ as a set-valued mapping and consider $\betae$, $\gammae$ as vector/scalar-valued mappings of $\paragroup$, then
		\begin{itemize}
			\item $\Se$ is osc relative to $\Pscr$,
			
			\item $\betae$ and $\gammae$ are continuous relative to $\Pscr$.
		\end{itemize}
		
		\item \label{it:continuity:lmd} Given $(\Am,\rho)\in\Real^{m\times n}\times [0,1)$. Consider $\Se$ as a set-valued mapping and consider $P_{\Se}(\bar{\wv})$ as a vector-valued mapping of $(\bv,\lambda)$, then
		\begin{itemize}
			\item $\Se$ is continuous relative to $\Real^{2m}\times\Real_{++}$,
			
			\item $P_{\Se}(\bar{\wv})$ is continuous relative to $\Real^{2m}\times\Real_{++}$; in particular, $\wstar$ is continuous relative to $\Real^{2m}\times\Real_{++}$.
		\end{itemize}
	\end{enumerate} 
\end{theorem}

We remark that continuity of the minimum $\ell_2$-norm LASSO regularization path $\xstar(\lambda)$, i.e., a special case of the second statement in Theorem \ref{thm:continuity} \ref{it:continuity:lmd}, is known under the ``one-at-a-time" assumption in previous work \cite{tibshirani2013}. Except this, all results in Theorem \ref{thm:continuity} is new.

Notably, to prove Theorem \ref{thm:continuity}, it is sufficient to prove outer semicontinuity of $\Se$ in $\paragroup$ and its inner-semicontinuity in $(\bv,\lambda)$, whilst vector-valued continuity of $\betae$, $\gammae$ and $\wstar$ can be implied from semicontinuity of $\Se$, which demonstrates the power of set-valued analysis. It is also worth noting that outer-semicontinuity of a set-valued mapping can usually be proven in a more direct way\footnote{As pointed out by one of the reviewers, for minimization problems, one can prove outer-semicontinuity of the solution set mapping fairly immediately from continuity of the objective function \cite[Thm. 8]{hogan1973}. In Theorem \ref{thm:continuity} \ref{it:continuity:para} of this study, since we are studying saddle point set of a min-max problem, \cite[Thm. 8]{hogan1973} does not apply here, but we can still prove outer-semicontinuity of $\Se$ quite directly using (\ref{eq:OPT}) and Theorem \ref{thm:shape_uniqueness_sparseness}.}, whilst it is typically harder to verify inner semicontinuity \cite[p. 155, para. 2]{rockafellar2009}. In this study, we address the difficulty of establishing inner-semicontinuity using a novel repeated-decomposition technique; see the following remark. The complete proof of Thm. \ref{thm:continuity} will be presented in Sec. \ref{sec:proof:continuity}.

\begin{remark}[Proof sketch of inner semicontinuity of $\Se$]
	The key point of our proof for inner semicontinuity lies in Lemma \ref{lemma:continuity_of_composition} (see Sec. \ref{sec:proof:continuity}), which indicates that the composition of an osc/isc set-valued mapping and a continuous single-valued mapping is also osc/isc. Hence we can decompose a complicated set-valued mapping into the composition of simpler mappings, and transfer continuity properties from simpler ones to the original complicated one. In the context of the sGMC model, we first use (\ref{eq:Se_variant_expression}) to decompose $\Se$ into the composition of the LASSO solution set mapping $\SLA$ and a continuous parameterization function, then further decompose $\SLA$ as a parameterized solution set mapping of a linear program. Inner semicontinuity of the solution set of linear programs can be derived from a previous work \cite{mangasarian1987} which uses different notions of set-valued continuity, whereby we can obtain inner semicontinuity of $\Se$.
\end{remark}

\subsection{Preparation for Studying Piecewise Simple Behavior of the Min-Norm Solution Map}
\label{sec:preparation_for_corres}
In the next two subsections, we introduce piecewise linearity and piecewise constant sparsity pattern of the min-norm solution map $\wstar(\bv,\lambda)$. The core of our results lies in an interesting one-to-one correspondence between continuous and discrete objects associated with $\wstar(\bv,\lambda)$, i.e.,
\begin{itemize}
	\item continuous objects: linear zones of $\wstar(\bv,\lambda)$,
	\item discrete objects: constant signed sparsity patterns $\sv\in\braces{+1,0,-1}^{2n}$ of $\wstar(\bv,\lambda)$ within each linear zone.
\end{itemize}
Such continuous-discrete correspondence is well-known in the study of the min-norm LASSO regularization path \cite{efron2004,tibshirani2013}, where the signed sparsity pattern $\sv$ and its support $\Ecal\coloneqq \supp(\sv)$ are respectively known as \textit{equicorrelation signs and set} in \cite{tibshirani2013}. In this paper, when characterizing this correspondence for $\wstar(\bv,\lambda)$ of the sGMC model, we do not explicitly use the support of $\sv$. In addition, we refer to $\sv$ by a shorter name: \textit{indicator}, emphasizing the fact that it embraces all information about a linear zone.

\begin{definition}\label{def:indicator_space} For every $\sv\in\braces{+1,0,-1}^{2n}$, we call $\sv$ a \textit{candidate indicator}. If $\sv$ equals to $\sign(\wstar(\bv,\lambda))$ within the interior some linear zone of $\wstar(\cdot,\cdot)$, then we further call it a \textit{zone indicator} or the indicator of that linear zone.
\end{definition}

In the sequel, we introduce necessary tools for connecting the continuous and discrete objects, which mainly include
\begin{itemize}
	\item an ``encoding function" $\sOPT(\wv;\bv,\lambda)$ that can map a vector $\wv\in\Real^{2n}$ into a candidate indicator $
	\sv$ at $(\bv,\lambda)$,
	
	\item a ``decoding function" $\wEQ(\sv;\bv,\lambda)$ that can map a candidate indicator $\sv$ into a linear function of $(\bv,\lambda)$.
\end{itemize}
These tools will be used in the next subsection to characterize the continuous-discrete correspondence of $\wstar(\bv,\lambda)$, which further addresses its piecewise simple behavior.

We first present the ``encoding function" $\sOPT(\wv;\bv,\lambda)$. Recall that given $(\bv,\lambda)\in\Real^m\times\Real_{++}$, $\wv\in\Se(\bv,\lambda)$ satisfies the following inclusion problem for every $i\in\braces{1,2,\dots,2n}$:
\begin{equation}
\cv_i^{\Tr}(\bv-\Dm\Cm\wv)\in\begin{cases}
\braces{\lambda\,\sign([\wv]_i)}, & \textrm{if }[\wv]_i\neq 0,\\
[-\lambda,\lambda], & \textrm{if }[\wv]_i=0.
\end{cases}\tag{\ref{eq:OPT}}
\end{equation}
Hence we can estimate the value of $\sign(\wv)$ (which appears in the RHS of (\ref{eq:OPT})) by $\xi_i(\wv)\coloneqq \cv_i^\Tr(\bv-\Dm\Cm\wv)$ (i.e., the LHS of (\ref{eq:OPT})). Extending this intuition from an extended solution $\wv\in\Se(\yv,\lambda)$ to an arbitrary vector $\wv\in\Real^{2n}$ leads to the definition of $\sOPT(\wv;\yv,\lambda)$.

\begin{definition}\label{def:Es_mapping}
	For $(\bv,\lambda)\in\Real^{2m}\times\Real$ and $\wv\in\Real^{2n}$,
	\begin{enumerate}[label=\arabic*)]
		\item we define $\sOPT(\wv;\bv,\lambda)\in\braces{+1,0,-1}^{2n}$ as:
		\begin{align*}
			\brackets{\sOPT(\wv;\bv,\lambda)}_i&\coloneqq \begin{cases}
				\sign\left(\xi_i(\wv)\right), & \textrm{if }\abs{\xi_i(\wv)}=\lambda,\\
				0, & \textrm{otherwise},
			\end{cases}
		\end{align*}
		where $i\in\braces{1,2,\dots,2n}$, $\xi_i(\wv)$ is the LHS of (\ref{eq:OPT}):
		\begin{equation*}
			\xi_i(\wv)\coloneqq \cv_i^\Tr(\bv-\Dm\Cm\wv).
		\end{equation*}
		One can verify that $\sOPT(\cdot;\bv,\lambda)$ maps an arbitrary vector $\wv\in\Real^{2n}$ to a candidate indicator $\sv$.
		
		\item Notice that for $\wv\in\Se(\bv,\lambda)$, the value of $\sOPT(\wv;\bv,\lambda)$ is only determined by $\Cm\wv$, which only depends on $(\bv,\lambda)$ (Thm. \ref{thm:shape_uniqueness_sparseness} \ref{it:shape:linear_fit}). By this observation, we further define \textit{equicorrelation signs} of the sGMC model as
		\begin{equation}\label{eq:equi_pair}
			\se(\bv,\lambda)\coloneqq \sOPT(\wv;\bv,\lambda), \text{ where }\wv\in\Se(\bv,\lambda),
		\end{equation}
		where ``c" in $\se(\cdot,\cdot)$ stands for ``correlation", emphasizing the fact that $\se(\bv,\lambda)$ is computed by a correlation vector $\xiv(\wv)\coloneqq \Cm^\Tr(\bv-\Dm\Cm\wv)$.
	\end{enumerate}
\end{definition}

One can verify that in the LASSO case, $\se(\bv,\lambda)$ reduces to the conventional definition of equicorrelation signs\footnote{In \cite{tibshirani2013}, equicorrelation signs of LASSO is originally defined as a discrete-valued mapping of $\lambda$, which is later proven to be piecewise constant under certain assumptions. Hence the terminology naturally extends to the associated constant values within linear pieces. In contrast, in the present paper we define the candidate indicators as constant discrete objects from the very beginning. This is a huge difference in viewing these quantities, and we will explain its consequence in Remark \ref{remark:compare_proof_corres}.} \cite[Eq. (4) and (5)]{tibshirani2013}. Hence the ``encoding function" $\sOPT(\cdot;\bv,\lambda)$ is a slight extension of the classic equicorrelation signs.

Next we introduce the ``decoding function" $\wEQ(\cdot;\bv,\lambda)$ along with some auxiliary quantities. Suppose that $\sv$ is the signed sparsity pattern of the unknown min-norm solution map $\wstar(\bv,\lambda)$, and we want to estimate the value of $\wstar(\bv,\lambda)$ by some vector $\wv\in\Real^{2n}$. Then $\wv$ must satisfy $\sign(\wv)=\sv$, which is equivalent to
\begin{align}
	&(\forall i\in\Ecal) &&[\sv]_i[\wv]_i>0, \label{temp:EQNQ1} &\\
	&(\forall i\in\neg\Ecal)&& [\wv]_{i}=0, \label{temp:EQNQ2}&
\end{align}
where $\Ecal\coloneqq\supp(\sv)$. In addition, notice that $\wv$ must satisfy (\ref{eq:OPT}) for every $i\in\braces{1,2,\dots,2n}$, we have
\begin{align}
	&(\forall i\in\Ecal) &&\cv_i^\Tr(\bv-\Dm\Cm\wv)=\lambda[\sv]_i, \label{temp:EQNQ3} &\\
	&(\forall i\in\neg\Ecal)&& \abs{\cv_i^\Tr(\bv-\Dm\Cm\wv)}\leq \lambda. \label{temp:EQNQ4}&
\end{align}
Let us divide the four conditions (\ref{temp:EQNQ1})-(\ref{temp:EQNQ4}) into two groups based on whether they are equalities or inequalities. Then we can obtain the following system of linear equations:
\begin{subequations}
	\makeatletter
	\def\@currentlabel{EQ}
	\makeatother
	\label{eq:EQ}
	\renewcommand{\theequation}{EQ-\alph{equation}}
	\begin{empheq}[left=\empheqlbrace]{align}
		&(\forall i\in\Ecal) &&\cv_{i}^{\Tr}(\bv-\Dm\Cm\wv)=\lambda[\sv]_{i}, \label{eq:EQ-E} &\\
		&(\forall i\in\neg\Ecal)&& [\wv]_{i}=0, \label{eq:EQ-nE}&
	\end{empheq}
\end{subequations}
as well as the following system of linear inequations:
\begin{subequations}
	\makeatletter
	\def\@currentlabel{NQ}
	\makeatother
	\label{eq:NQ}
	\renewcommand{\theequation}{NQ-\alph{equation}}
	\begin{empheq}[left=\empheqlbrace]{align}
		&(\forall i\in\Ecal) && [\sv]_i[\wv]_i\geq 0, \label{eq:NQ-E} &\\
		&(\forall i\in\neg\Ecal) &&  |\cv_i^{\Tr}(\bv-\Dm\Cm\wv)|\leq \lambda. &\label{eq:NQ-nE}
	\end{empheq}
\end{subequations}
Notice that the strict inequality in (\ref{temp:EQNQ1}) is relaxed in (\ref{eq:NQ-E}). The names of (\ref{eq:EQ}) and (\ref{eq:NQ}) respectively stand for ``equalities" and ``inequalities". Based on these two groups of conditions, we can define the ``decoding function" $\wEQ(\sv;\bv,\lambda)$ and some other useful quantities.

\begin{definition}\label{def:EQ_NQ_quantities}
	Given $(\bv,\lambda)\in\Real^{2m}\times\Real$, for a candidate indicator $\sv\in\braces{+1,0,-1}^{2n}$, let $\Ecal\coloneqq\supp(\sv)$. If $\wv\in\Real^{2n}$ satisfies (\ref{eq:EQ}) (resp. (\ref{eq:NQ})), we say $\wv$ satisfies $\sv$-(\ref{eq:EQ}) (resp. (\ref{eq:NQ})) at $(\bv,\lambda)$ to emphasize the associated parameters.
	\begin{enumerate}[label=\arabic*)]
		\item We define $\wEQ(\sv;\bv,\lambda)$ as the unique least-squares fitting\footnote{Notice that $\sv$-(\ref{eq:EQ}) does not necessarily have a solution at $(\bv,\lambda)$.} of $\sv$-(\ref{eq:EQ}) at $(\bv,\lambda)$, and term it the \textit{candidate solution map} of $\sv$ at $(\bv,\lambda)$.
		
		\item We define $\Spoly(\sv;\bv,\lambda)$ as the set of $\wv\in\Real^{2n}$ that simultaneously satisfies $\sv$-(\ref{eq:EQ}) and $\sv$-(\ref{eq:NQ}) at $(\bv,\lambda)$, and term it the \textit{candidate solution set} of $\sv$ at $(\bv,\lambda)$.
		
		\item We define $\Iselect(\sv)$ as the set of $(\bv,\lambda)\in\Real^{m}\times\Real_{++}$ such that $\wEQ(\sv;\bv,\lambda)\in\Spoly(\sv;\bv,\lambda)$, and term it the \textit{candidate zone} of $\sv$.
	\end{enumerate}
\end{definition}
It should be noted that in the LASSO case, the candidate solution map $\wEQ(\sv;\bv,\cdot)$ reduces to the function used in the LARS iteration \cite[Eq. (15)]{tibshirani2013} for computing the min-norm LASSO regularization path, whereas (\ref{eq:NQ}) in the definition of $\Iselect(\sv)$ reduces to the conditions for computing ``joining time" and ``crossing time" in the LARS iteration \cite[Alg. 1]{tibshirani2013}.

In the following, we interpret the meanings of all quantities defined above and make some remarks on their computation.

\begin{remark}[Interpretation and computation of Definition \ref{def:EQ_NQ_quantities}]\label{remark:computation_EQ_NQ}
	For $\sv\in\braces{+1,0,-1}^{2n}$ and $(\bv,\lambda)\in\Real^{2m}\times\Real_{++}$,
	\begin{enumerate}[label=\arabic*)]
		\item $\wEQ(\sv;\bv,\lambda)$, $\Spoly(\sv;\bv,\lambda)$ and $\Iselect(\sv)$ respectively can be interpreted as approximations of
		\begin{itemize}
			\item the min-norm solution map $\wstar(\bv,\lambda)$,
			\item the extended solution set $\Se(\bv,\lambda)$,
			\item the range of the linear zone corresponding to $\sv$,
		\end{itemize}
		which are obtained by imagining that $(\bv,\lambda)$ lies in the linear zone of $\sv$; this accounts for their names.
		
		\item Given a candidate indicator $\sv$, then
		\begin{itemize}
			\item $\wEQ(\sv;\bv,\lambda)$ is a linear function of $(\bv,\lambda)\in\Real^{2m}\times\Real$ whose slope can be computed in closed-form,
			
			\item $\Iselect(\sv)$ is a convex cone in $\Real^{2m}\times\Real_{++}$ which can be computed by solving a linear program. Moreover, for any straight line $L\subset\Real^{2m}\times\Real$, $L\cap \Iselect(\sv)$ can be computed in closed-form.
		\end{itemize}
		See Appx. \ref{app:computation} for details of their computation.
	\end{enumerate}
\end{remark}

Finally, we present two lemmas, highlighting some useful properties of the three quantities in Def. \ref{def:EQ_NQ_quantities}. The proofs of these two lemmas can be found in Appx. \ref{app:auxiliary_results}.

\begin{lemma}\label{lemma:SEQNQ_property}
	Given $\sv\in\braces{+1,0,-1}^{2n}$, then we always have $\Spoly(\sv;\bv,\lambda)\subset\Se(\bv,\lambda)$ for every $(\bv,\lambda)\in\Real^{2m}\times\Real_{++}$, where $\Se(\bv,\lambda)$ is the true extended sGMC solution set.
\end{lemma}

\begin{lemma}\label{lemma:IEN_property}
	For $\sv\in\braces{+1,0,-1}^{2n}$ and $(\bv,\lambda)\in\Real^{2m}\times\Real_{++}$, $(\bv,\lambda)\in\Iselect(\sv)$ is equivalent to the following statements:
	\begin{enumerate}[label=(\alph*)]
		\item \label{it:IEN:wEQ_in_Se} $\wEQ(\sv;\bv,\lambda)\in\Se(\bv,\lambda)$,
		
		\item \label{it:IEN:wEQ_min_in_SEQNQ} $\wEQ(\sv;\bv,\lambda)$ is the unique minimum $\ell_2$-norm element in $\Spoly(\sv;\bv,\lambda)\subset\Se(\bv,\lambda)$.
	\end{enumerate}
	Moreover, for all $(\bv,\lambda)\in\interior(\Iselect(\sv))$, $\sign(\wEQ(\sv;\bv,\lambda))$ and $\sOPT(\wEQ(\sv;\bv,\lambda);\bv,\lambda)$ respectively remain constant.
\end{lemma}

\subsection{Piecewise Simple Behavior of the Min-Norm Solution Map and Continuous-Discrete Correspondence}
\label{sec:correspondence}

Exploiting the tools introduced above, in the sequel, we formally present our results about piecewise simple behavior of the min-norm solution map $\wstar(\bv,\lambda)$ and the continuous-discrete correspondence therein. 

\begin{theorem}\label{thm:correspondence}
	Given $(\Am,\rho)\in\Real^{m\times n}\times[0,1)$, then the following holds for the minimum $\ell_2$-norm extended solution map $\wstar(\bv,\lambda)$ of the sGMC model:
	\begin{enumerate}[label=(\alph*)]
		\item \label{it:corres:PL} $\wstar(\bv,\lambda)$ is piecewise linear in $(\bv,\lambda)\in\Real^{2m}\times\Real_{++}$ with finite ($\leq 3^{2n}$) linear zones, with $\sign(\wstar(\bv,\lambda))$ being constant within the interior of each linear zone.
		
		\item \label{it:corres:zone}
		The correspondence between a linear zone in $\wstar(\cdot,\cdot)$ and its zone indicator $\sv$ (cf. Def. \ref{def:indicator_space}) is one-to-one. Moreover, let $\sv$ be a zone indicator and let $I(\sv)\subset\Real^{2m}\times\Real_{++}$ be the range of its linear zone, then $I(\sv)$ is a convex cone contained in $\Iselect(\sv)$ and
		\begin{align*}
			(\forall (\bv,\lambda)\in I(\sv))&&\wstar(\bv,\lambda)&=\wEQ(\sv;\bv,
			\lambda).
		\end{align*}

		\item \label{it:corres:GP} If $\Am$ has columns in general position, then the result \ref{it:corres:zone} can be improved to $I(\sv)=\Iselect(\sv)$ and
		\begin{align*}
			(\forall (\bv,\lambda)\in I(\sv))&&\wstar(\bv,\lambda)&=\wEQ(\sv;\bv,
			\lambda),\\
			(\forall (\bv,\lambda)\in I(\sv))&& \Se(\bv,\lambda)&=\Spoly(\sv;\bv,\lambda).
		\end{align*}
		Moreover, for every $(\bv,\lambda)\in\interior(\Iselect(\sv))$, we have
		\begin{align}
			\sv \equiv \sign(\wstar(\bv,\lambda))\equiv \se(\bv,\lambda).\label{eq:s_sign_csign}
		\end{align}
	\end{enumerate}
\end{theorem}

We note that Thm. \ref{thm:correspondence} is well known\footnote{Strictly speaking, the characterization of $\Se(\bv,\lambda)$ by $\Spoly(\sv;\bv,\lambda)$ is slightly different from the previous result \cite[Eq. (9)]{tibshirani2013} in two aspects: firstly, $(\Ecal,\sv)$ in \cite[Eq. (9)]{tibshirani2013} indeed stands for $(\supp(\se(\bv,\lambda)),\se(\bv,\lambda))$ in this paper, which does not necessarily equal to $(\supp(\sv),\sv)$ for the zone indicator $\sv$ on the boundary of $\Iselect(\sv)$; secondly, the condition (\ref{eq:NQ-nE}) in the definition of $\Spoly(\sv;\bv,\lambda)$ is not included in \cite[Eq. (9)]{tibshirani2013}.} for the min-norm LASSO regularization path under certain assumptions \cite{efron2004,tibshirani2013}. Thm. \ref{thm:correspondence} improves previous work in two aspects:
\begin{enumerate}[label=\arabic*)]
	\item all properties stated in Thm. \ref{thm:correspondence} are generalized from the min-norm LASSO regularization path $\xstar(\lambda)$ to the min-norm sGMC extended solution map $\wstar(\bv,\lambda)$. Notably,
	\begin{itemize}
		\item previously, $\xstar(\lambda)$ is known to be piecewise affine in $\lambda$, where each linear zone is an interval in $\Real_{++}$,
		
		\item our study indicates that $\wstar(\bv,\lambda)$ is exactly piecewise ``linear" in $(\bv,\lambda)$, and linear zones therein are newly proven to be convex cones in $\Real^{2m}\times\Real_{++}$.
	\end{itemize}
	
	\item the proofs of Thm. \ref{thm:correspondence} \ref{it:corres:PL} \ref{it:corres:zone} do not require any assumption that previous work \cite{tibshirani2013} relies on (whereas the proof of Thm. \ref{thm:correspondence} \ref{it:corres:GP} uses a different assumption from \cite{tibshirani2013}, see Sec. \ref{sec:comments_on_tib13} for the assumptions implicitly used in \cite{tibshirani2013}),
\end{enumerate}

In the following, we discuss about the results of Thm. \ref{thm:correspondence}.

\begin{remark}[Discussion about Theorem \ref{thm:correspondence}]\label{remark:consequence_correspondence}
	\begin{enumerate}[label=\arabic*)]
		\item Thm. \ref{thm:correspondence} \ref{it:corres:PL} \ref{it:corres:zone} imply that for a fixed linear system $\Am$, the sGMC/LASSO solution map can be exactly recovered by a finite-depth finite-width neural network, which draws a connection with deep learning methods for sparse recovery \cite{gregor2010}. While the piecewise linear expression of $\wstar(\cdot,\cdot)$ can be computed via the E-LARS algorithm to be introduced later, we do not know how to transform it into a compact neural network. Finding this neural network representation may be an interesting future direction.
		
		\item It is well known that the min-norm LASSO regularization path has at most $3^n$ linear pieces\footnote{Interestingly, it has been proven that the worst-case linear piece number is exactly $(3^n+1)/2$. See \cite[Prop. 2]{mairal2012} for an adversarial strategy for building problem instances with such pathological regularization paths.}. However, in practice such worst-case linear piece number rarely occurs. Thm. \ref{thm:correspondence} \ref{it:corres:zone} provides an intuitive explanation for this phenomenon: for the LASSO model, we can show that the $(m+1)$-dimensional parameter space $(\yv,\lambda)\in\Real^{m}\times\Real_{++}$ is divided by at most $3^{n}$ linear zones (which are convex cones), thus as the min-norm solution $\xstar(\yv,\lambda)$ travels along a 1D line parallel to the $\lambda$-axis, the number of cones it can visit should be much smaller than $\Ocal(3^{n})$ unless these cones are very ill-positioned.
		
		\item \label{it:discuss_corres:GP_not_necessary} It should be noted that the general-positioning condition is not a necessary condition for Thm. \ref{thm:correspondence} \ref{it:corres:GP}. One can consider the following counterexample: let $\Am\coloneqq \begin{bmatrix}
			1 & 1
		\end{bmatrix}$ and $\rho=0$ in the sGMC model (\ref{eq:sGMC_model}). Then evidently $\Am$ does not have columns in general position. Moreover, by comparing the norm of $\wEQ(\sv;\yv,\lambda)$ for every candidate indicator $\sv$ within $\Iselect(\sv)$, one can verify that there exists only three zone indicators (here we omit the dual parts in the candidate indicator $\sv$ and parameter $\bv$):
		\begin{align*}
			\sv_0&\coloneqq \mathbf{0}_{2}\text{ with }\Iselect(\sv_0)=\braces{(y,\lambda)\;\middle\vert\; \abs{y}\leq \lambda,\; \lambda>0},\\
			\sv_1&\coloneqq \mathbf{1}_{2}
			\text{ with }\Iselect(\sv_1)=\braces{(y,\lambda)\;\middle\vert\; y\geq \lambda,\; \lambda>0},\\
			\sv_2&\coloneqq - \mathbf{1}_{2}
			\text{ with }\Iselect(\sv_2)=\braces{(y,\lambda)\;\middle\vert\; y\leq -\lambda,\; \lambda>0},
		\end{align*}
		and all results stated in Thm. \ref{thm:correspondence} \ref{it:corres:GP} holds. Indeed, the authors believe that Thm. \ref{thm:correspondence} \ref{it:corres:GP} should hold without requiring any condition, but currently they cannot prove this conjecture, and would like to leave it for future work.
	\end{enumerate}
\end{remark}

In the following remark, we give a comparison between our proof strategy of Theorem \ref{thm:correspondence} and that of \cite{efron2004,tibshirani2013}. More detailed comments on the proof of \cite{tibshirani2013} will be given in Sec. \ref{sec:comments_on_tib13}. See Sec. \ref{sec:proof:correspondence} for the complete proof of Theorem \ref{thm:correspondence}.

\begin{remark}[Comparison between  our proof strategy of Theorem \ref{thm:correspondence} and that of \cite{efron2004,tibshirani2013}]\label{remark:compare_proof_corres}
	As mentioned in the footnote of Definition \ref{def:Es_mapping}, the main subject of study in \cite{tibshirani2013} is the discrete-valued mapping $\se(\lambda)$ of $\lambda$, whilst that of the current study is the finite set of candidate indicators. This leads to a huge difference in the proof strategies.
	\begin{enumerate}[label=\arabic*)]
		\item The dynamic nature of $\se(\lambda)$ naturally leads to a constructive approach to establish the continuous-discrete correspondence in previous work \cite{efron2004,tibshirani2013}, which uses the LARS algorithm to track the dynamic behavior of $\se(\lambda)$ and proves its correctness by induction.
		
		\item In contrast, the static nature of candidate indicators leads to a geometric approach to study the continuous-discrete correspondence. Our proof is developed by introducing a partial order $\preceq$ among candidate indicators. It will be shown in Prop. \ref{prop:partial_order} that
		\begin{itemize}
			\item under certain conditions, there exists an inequality
			\[\sign(\wstar(\bv,\lambda))\preceq \sv\preceq\se(\bv,\lambda)\]
			that connects the three quantities in the equality (\ref{eq:s_sign_csign}),
			
			\item when this inequality attains partial equalities, we can further draw an equality between $\wstar(\bv,\lambda)$ and $\wEQ(\sv;\bv,\lambda)$, or between $\Se(\bv,\lambda)$ and $\Spoly(\sv;\bv,\lambda)$, leading to the results stated in Thm. \ref{thm:correspondence} \ref{it:corres:GP}.
		\end{itemize}		
		Studying the conditions that guarantee partial equalities of the discrete inequality above yields Theorem \ref{thm:correspondence}.
	\end{enumerate}	 
\end{remark}

Finally, before proceeding to the E-LARS algorithm for computing $\wstar(\cdot,\cdot)$, we present a corollary that will be useful for its initialization: in general, zone indicators in $\wstar(\cdot,\cdot)$ depend on the problem parameter $(\Am,\rho)$, but there is one exception, i.e., $\sv_0\coloneqq\mathbf{0}_{2n}$ is always a zone indicator.   

\begin{corollary}\label{cor:zeroth_linear_piece}
	Given $(\Am,\rho)\in\Real^{m\times n}\times[0,1)$, then $\sv_0\coloneqq\mathbf{0}_{2n}$ is a zone indicator of $\wstar(\cdot,\cdot)$ with its linear zone being
\begin{align*}
	I(\sv_0)=\Iselect(\sv_0)	=\braces{(\bv,\lambda)\in\Real^{2m}\times\Real_{++}\;\middle\vert\; \max_{i\in\braces{1,2,\dots,2n}} \abs{\cv_i^\Tr\bv}\leq\lambda}
\end{align*}
	and $\wstar(\bv,\lambda)\equiv \mathbf{0}_{2n}$ within this linear zone.
\end{corollary}
\begin{IEEEproof}
	One can verify that for every $(\bv,\lambda)\in \Iselect(\sv_0)$, $\wEQ(\sv_0;\bv,\lambda)=\mathbf{0}_{2n}$ satisfies (\ref{eq:OPT}), thus $\wstar(\bv,\lambda)= \mathbf{0}_{2n}$. Then from Thm. \ref{thm:correspondence} \ref{it:corres:zone} we can yield $I(\sv_0)=\Iselect(\sv_0)$.
\end{IEEEproof}

\section{Extended Least Angle Regression}\label{sec:E_LARS}

In this subsection, we introduce the extended least angle regression (E-LARS) algorithm for computing the whole min-norm extended solution map $\wstar(\cdot,\cdot)$ of sGMC, which addresses the third main contribution stated in Sec. \ref{sec:contributions}.

In previous work \cite{osborne2000LARS,efron2004,tibshirani2013}, the LARS algorithm for LASSO is developed as a homotopy method, and the correspondence between linear pieces and their indicators is deduced from its iterative step. In contrast, in this paper, since we have established the continuous-discrete correspondence independently, we can view the computation of $\wstar(\cdot,\cdot)$ from a different perspective, treating it as a discrete search problem for all zone indicators.

The E-LARS algorithm intends to solve this discrete search problem in an iterative manner, i.e., given the indicator $\sv$ of a linear zone in $\wstar(\cdot,\cdot)$, the E-LARS iteration estimates an indicator of its adjacent linear zone. See the following.

\subsection{Iterative Step, Initialization and Termination of E-LARS}\label{sec:step_ini_termi_ELARS}

Let the input of the E-LARS iteration be some zone indicator $\sv$, and let $\sadj$ be the indicator of an adjacent linear zone. Then without loss of generality, we can assume there exists a straight line connecting $\interior(I(\sv))$ and $\interior(I(\sadj))$. Let us move the parameter $(\bv,\lambda)$ along this straight line uniformly from $\interior(I(\sv))$ to $\interior(I(\sadj))$, i.e., we set
\begin{equation}\label{eq:line_parameter}
	\begin{aligned}
		\bv(t) & = \bv_0+\Deltab t,\\
		\lambda(t) & = \lambda_0 +\Deltalmd t,
	\end{aligned}
\end{equation}
where $t\in\Real$ is a time index, and we assume that the original point $(\bv_0,\lambda_0)$ and velocity vector $(\Deltab,\Deltalmd)$ are known. Then one can imagine that as time $t$ increases from $-\infty$ to $+\infty$, the min-norm sGMC solution
\begin{align}
	\wstar(t)\coloneqq \wstar(\yv(t),\lambda(t))
\end{align}
first enters $\interior(I(\sv))$ at some time, then leaves it and immediately enters $\interior(I(\sadj))$ at another time, say $\tadj$:
\begin{equation}\label{eq:t_plus}
	\begin{aligned}
		\tadj\coloneqq& \sup_{t\in\Real} \braces{t\;\middle\vert\; \wstar(t)\in\interior(I(\sv))}\\
		=&\inf_{t\in\Real}\braces{t\;\middle\vert\; \wstar(t)\in\interior(I(\sadj))}.
	\end{aligned}
\end{equation}
The goal of the E-LARS iteration is to generate an estimate $\sv_+$ for $\sadj$, using the hyperparameters $(\bv_0,\lambda_0)$ and $(\Deltab,\Deltalmd)$.

The basic idea of the E-LARS iteration is to estimate the difference between $\sv$ and $\sadj$ by observing the behavior of $\wstar(t)$ at the zone-switching time $\tadj$, then construct the estimate $\sv_+$ by modifying the value of $\sv$. In the sequel, we present a sketch of the E-LARS iteration. Please see Appx. \ref{app:ELARS_deri} for its derivation and see Alg. \ref{alg:LARS_sGMC} therein for pseudocode.

\begin{procedure}[E-LARS iteration]\label{procedure:LARS-sGMC}
	The E-LARS iteration computes an estimate $\sv_+$ for $\sadj$ by the following steps:
	\begin{enumerate}[label=\arabic*)]
		\item {Compute the expression of $\wstar(t)$ in the linear zone of $\sv$ and estimate the zone-switching time $\tadj$.}
		\begin{itemize}
			\item We compute $\wEQ(\sv;\yv(t),\lambda(t))$ to obtain the closed-form expression of $\wstar(t)$ in the linear zone of $\sv$,
			
			\item We generate an estimate $t_+$ of $\tadj$ by computing the intersection between $\Iselect(\sv)$ and the straight line that $(\bv(t),\lambda(t))$ moves along (cf. Lemma \ref{lemma:compuation_IEN_line} in Appx. \ref{app:computation}).
		\end{itemize}
		
		\item {Perform suitable deletion and insertion operations on the components of $\sv$ to obtain $\sv_+$. More precisely,}
		\begin{itemize}
			\item for every $i\in\Ecal$, if $\wstar(t_+)$ attains equality of the $i$th $\sv$-(\ref{eq:NQ-E}) condition at $(\bv(t_+),\lambda(t_+))$, then we delete this component $[\sv]_i$, i.e., we set $[\sv_+]_i\gets 0$,
			
			\item for every $i\in\neg\Ecal$, if $\wstar(t_+)$ attains equality of the $i$th $\sv$-(\ref{eq:NQ-nE}) condition at $(\bv(t_+),\lambda(t_+))$, then we insert a suitable nonzero sign to $[\sv]_i$, i.e., we set 
			\[[\sv_+]_i\gets \sign\parens{\cv_i^\Tr(\bv(t_+)-\Dm\Cm\wstar(t_+))}.\]
		\end{itemize}
		For other index $i$, we make no change, i.e., $[\sv_+]_i\gets[\sv]_i$.
	\end{enumerate}
\end{procedure}

Perhaps unsurprisingly, the E-LARS iteration resembles that of the conventional LARS-LASSO algorithm. In particular, if $\rho=0$, $\rv=\mathbf{0}$ in the sGMC model (\ref{eq:def_of_G}) and $\Deltab=\mathbf{0}$, then the E-LARS iteration reduces to the conventional LARS iteration.

Notice that to run the E-LARS algorithm, one needs to have an initial zone indicator as its input of the first iteration, and needs to obtain suitable $(\bv_0,\lambda_0,\Deltab,\Deltalmd)$ in (\ref{eq:line_parameter}) before the computation of each iteration. See the following remark for some possible strategies for these initialization issues.

\begin{remark}[Initialization of E-LARS]\label{remark:initial_ex_equipair}
	 Regarding the selection of the initial zone indicator, we suggest the following two strategies as pointed out in \cite[Appendix A.3]{tibshirani2013}:
	\begin{enumerate}[label=\arabic*)]
		\item From Cor. \ref{cor:zeroth_linear_piece}, we can always select $\sv_0=\mathbf{0}_{2n}$.
		
		\item For arbitrary $(\bv,\lambda)\in\Real^{2m}\times\Real_{++}$, we can compute an extended solution $\wv\in\Se(\bv,{\lambda})$ by some iterative algorithm. Then if $(\bv,\lambda)$ is not on the boundary of a linear zone, we can select $\sOPT(\wv;\bv,\lambda)$ as initial input.
	\end{enumerate}
	On the other hand, to obtain suitable $(\bv_0,\lambda_0,\Deltab,\Deltalmd)$ in (\ref{eq:line_parameter}), we can compute one face of $\Iselect(\sv)$ (see, e.g., \cite{fukuda1997}), then if Thm. \ref{thm:correspondence} \ref{it:corres:GP} holds, we can set $(\bv_0,\lambda_0)$ as a point in the interior of this face, and can set $(\Deltab,\Deltalmd)$ as a normal vector of it.
\end{remark}

Recall that there are at most $3^{2n}$ linear zones in $\wstar(\cdot,\cdot)$, in general we cannot obtain all the zone indicators with limited computational resources. Thus in practice, it is crucial to search zone indicators in a proper order and to terminate the algorithm early according to the needs of the user. See the following remark for a possible termination rule.

\begin{remark}[Termination of E-LARS]\label{remark:termi_ELARS}
	In real-world applications, we usually can assume that the observed signal $\yv\in\Real^{m}$ has a bounded norm, say $\norm{\yv}_2\leq R_y$, and the regularization parameter $\lambda$ we use is above some small number, say $\lambda\geq \delta_{\lambda}$. Then since all linear zones of $\wstar(\cdot,\cdot)$ are convex cones, we can terminate the E-LARS algorithm when the linear zones we have already obtained cover the following parameter set:
	\begin{equation*}
		\Lambda\coloneqq\braces{(\bv,\lambda)\in\Real^{2n}\times\Real_{++}\;\middle\vert\;  \norm{\yv}_2\leq R_y,\;\rv=\mathbf{0}_{2n},\; \lambda=\delta_{\lambda}}.
	\end{equation*}
	To efficiently search linear zones that intersect with $\Lambda$, we can further equip the E-LARS iteration with a breadth-first search strategy (where we can consider zone indicators with adjacent linear zones as adjacent vertices of an undirected graph).
\end{remark}

\subsection{Conditional Correctness and Complexity of E-LARS}

Subsequently, we present results about conditional correctness of E-LARS. Here we introduce two useful assumptions. The first one is the conventional ``one-at-a-time" assumption widely used in previous work \cite{efron2004,rosset2007,tibshirani2011,tibshirani2013}.

\begin{assumption}[One-at-a-time]\label{assump:one_at_a_time}
	Let $\sv$ and $\sv_+$ respectively be the input and output of the E-LARS iteration. We say the ``one-at-a-time" assumption holds if $\sv$ and $\sv_+$ only differs on one component.
\end{assumption}

In addition, we introduce a novel assumption termed ``sticky-equality" assumption.

\begin{assumption}[Sticky-equality]\label{assump:sticky_equality}
	For $t_+$ be the zone-switching time in (\ref{eq:t_plus}). We say the ``sticky-equality" assumption holds if for every $i\in\braces{1,2,\dots,2n}$, we have:
	\begin{itemize}
		\item the equality $[\wstar(t_+)]_i=0$ implies the existence of a neighborhood $U$ of $t_+$ such that 
		\[(\forall t\in U)\;\; [\wstar(t)]_i=0,\]
		
		\item the equality $\abs{\cv_i^\Tr(\bv(t_+)-\Dm\Cm\wstar(t_+))}=\lambda(t_+)$ implies the existence of a neighborhood $V$ of $t_+$ such that 
		\[(\forall t\in V)\;\; \abs{\cv_i^\Tr(\bv(t)-\Dm\Cm\wstar(t))}=\lambda(t).\]
	\end{itemize}
	In words, the ``sticky-equality" assumption requires the following properties of $\wstar(t)$: if $\wstar(t)$ attains equality of the $i$th (\ref{eq:NQ-E}) or (\ref{eq:NQ-nE}) condition at $t_+$, then this equality will last for a neighborhood of $t_+$.
\end{assumption}

Exploiting the aforementioned two assumptions, in the sequel, we present conditional correctness of E-LARS iteration.

\begin{theorem}\label{thm:E_LARS}
	Let $\sv$, $\sadj$ and $\sv_+$ be as described in the first paragraph of Sec. \ref{sec:step_ini_termi_ELARS}. Suppose that Thm. \ref{thm:correspondence} \ref{it:corres:GP} holds for the linear zones of $\sv$ and $\sadj$. Then the following holds:
	\begin{enumerate}[label=(\alph*)]
		\item \label{it:LARS:sticky} The E-LARS iteration is correct, i.e. $\sv_+=\sadj$, if and only if the sticky-equality assumption holds.
		
		\item \label{it:LARS:one} The one-at-a-time assumption implies the sticky-equality assumption, while the reverse does not hold.
	\end{enumerate}
\end{theorem}

We note that correctness of the LARS algorithm for LASSO under the ``one-at-a-time" assumption is known in previous work \cite{efron2004,tibshirani2013}. The other results of Theorem \ref{thm:E_LARS} are novel.

Theorem \ref{thm:E_LARS} can be proven quite directly from the continuity of $\wstar(t)$ (Thm. \ref{thm:continuity}) and the discrete equality (\ref{eq:s_sign_csign}) in Thm. \ref{thm:correspondence} \ref{it:corres:GP}; see Appx. \ref{proof:E_LARS} for its proof.

In the following remark, we present a comparison between the two assumptions used in Theorem \ref{thm:E_LARS}.

\begin{remark}[Comparison between Assumption \ref{assump:one_at_a_time} and \ref{assump:sticky_equality}]\label{remark:compare_assumptions}
	\begin{enumerate}[label=\arabic*)]
		\item The one-at-a-time condition assumes certain properties of the output of the E-LARS iteration, whilst the sticky-equality condition assumes certain geometric properties of the min-norm solution $\wstar$.
		
		\item Let $t_i$ be the time for $\wstar(t)$ to attain equality of the $i$th $\sv$-(\ref{eq:NQ-E}) (resp. (\ref{eq:NQ-nE})) condition for $i\in\Ecal$ (resp. $i\in\neg\Ecal$). Then one can imagine that if $(t_i)_{i\in\braces{1,2,\dots,2n}}$ follows a continuous probability distribution on $\Real^{2n}$, then the one-at-a-time assumption almost surely holds, which provides a non-rigorous interpretation for its empirical correctness in experiments. This further implies mildness of the sticky-equality assumption from Thm. \ref{thm:E_LARS} \ref{it:LARS:one}. 
		
		\item The one-at-a-time assumption does not hold in general\footnote{To see a counterexample leading to failure of Assumption \ref{assump:one_at_a_time}, one can imagine a sensing matrix $\Am\coloneqq \begin{bmatrix}
			\bar{\Am} & \bar{\Am}
		\end{bmatrix}$ with arbitrary matrix $\bar{\Am}$. In this case, since every column in $\Am$ has its replica, the E-LARS iteration cannot distinguish their corresponding components in $\wstar(t)$, whereby the deletion-insertion operations in Procedure \ref{procedure:LARS-sGMC} always involve more than two indices.}. In contrast, whether the sticky-equality assumption always holds is still an open problem. 
	\end{enumerate}
\end{remark}

Finally, we present a complexity analysis of the E-LARS algorithm, followed by a discussion on its efficiency.

\begin{corollary}
	The complexity of E-LARS is as follows:
	\begin{enumerate}[label=(\alph*)]
		\item The time complexity of the E-LARS iteration is \[\Ocal\parens{mn+m\abs{\supp(\sv)}^2+\abs{\supp(\sv)}^3},\]
		where $\sv$ is the input zone indicator, and the space complexity is $\Ocal(m+n)$. See Table \ref{tab:complexity_of_LARS_sGMC} in Appx. \ref{app:ELARS_deri} for a step-by-step analysis.
		
		\item Let $\sv_0=\mathbf{0}_{2n}$ be the initial zone indicator of the E-LARS algorithm, and set $R_y=+\infty$ in the stopping rule in Remark \ref{remark:termi_ELARS}. If the assumptions used in Thm. \ref{thm:E_LARS} hold for every E-LARS iteration, then it requires at most $3^{2n}$ iterations for E-LARS to obtain the whole min-norm extended solution map $\wstar(\yv,\lambda)$ of sGMC.
	\end{enumerate}
\end{corollary}

\begin{remark}[Efficiency of E-LARS]
	\begin{enumerate}[label=\arabic*)]
		\item For the E-LARS iteration, since $\abs{\supp(\sv)}\leq 2n$, its time complexity can be up to $\Ocal(n^3)$, which does not seem so efficient when $n$ is large. However, notice that $|\supp(\sv)|=\norm{\wstar(\bv,\lambda)}_0$ within the current linear zone, and in practice we are only interested in sparse solutions of sGMC. Hence one can expect that E-LARS is only applied in cases where $\abs{\supp(\sv)}\ll 2n$, i.e., the time complexity should be much lower than $\Ocal\parens{n^3}$.

		\item If the sensing matrix $\Am$ is generated randomly, the from Thm. \ref{thm:shape_uniqueness_sparseness} \ref{it:shape:unique}, $\abs{\supp(\sv)}\leq 2m$ holds with probability one. In this case the time complexity of the E-LARS iteration is at most $\Ocal(m^3)$, which is acceptable in the standard setting of sparse recovery problems where $m\ll n$.
		
		\item If we want to use E-LARS to compute sGMC solutions on a grid of problem parameters, then the required number of iterations depends on the number of linear zones the target parameters lie on. Especially, if the target parameters lie on a low-dimensional manifold in $\Real^{2m}\times\Real_{++}$, then one can imagine that with high probability, the required number of E-LARS iterations should be much smaller than the worst-case complexity $3^{2n}$.
	\end{enumerate}
\end{remark} 

\subsection{Comments on Proof Techniques of \cite{tibshirani2013}}
\label{sec:comments_on_tib13}

In the end of this section, we review the proof techniques of \cite{tibshirani2013} and the assumptions implicitly used there, and we analyze how one can adapt these techniques to obtain partial results of the current study.

We first introduce the proof idea of \cite[Appx. A.1]{tibshirani2013}. Suppose that we have already obtained the iterative step of the LARS algorithm (which can be derived from a homotopy method \cite{osborne2000LARS,efron2004}), then to study the geometric properties of the min-norm LASSO regularization path $\xstar(\lambda)$, a natural idea is to prove correctness of LARS by induction, whereby continuity, piecewise simple behvior and continuous-discrete correspondence of $\xstar(\lambda)$ can all be implied by properties of the LARS iteration. Here we run the LARS algorithm in the following setting:
\begin{enumerate}[label=\arabic*)]
	\item we select $\sv_0\coloneqq\mathbf{0}$ as the initial input indicator, such that the correctness of the first iteration can be easily verified,
	
	\item the parameter $(\yv,\lambda)$ moves as follows: $\yv\equiv\yv_0$ is fixed all the time, and $\lambda$ decreases from $+\infty$ to $0$.
\end{enumerate}

Let $\sv$, $\sv_+$ respectively be the input and output of the current LARS iteration, and let $\lambda_+$ be linear-piece-switching time estimated by the LARS iteration. The induction step of \cite[Appx. A.1]{tibshirani2013} is composed of two parts:
\begin{enumerate}[label=\arabic*)]
	\item demonstrating continuity of the regularization path generated by LARS at $\lambda_+$. This is achieved by the ``insertion-deletion" lemma \cite[Lemma 17]{tibshirani2013} which proves
	\begin{align}\label{eq:LARS_wEQ_continuity}
		\wEQ(\sv;\yv_0,\lambda_+)=\wEQ(\sv_+;\yv_0,\lambda_+)
	\end{align}
	via analysis of a block-structured linear system.
	
	\item demonstrating that $\wEQ(\sv_+;\yv_0,\lambda)$, i.e., the new linear piece generated by LARS, satisfies the optimality condition of LASSO (which is a variant of (\ref{eq:NQ}), see \cite[last paragraph in p. 1483]{tibshirani2013}) for every $\lambda$ within the new linear piece. This is shown by combining (\ref{eq:LARS_wEQ_continuity}) with the induction hypothesis, i.e., $\wEQ(\sv;\yv_0,\lambda)$ satisfies the optimality condition at $\lambda=\lambda_+$.
\end{enumerate}
The induction analysis above does not require the solution-uniqueness of LASSO, however, it implicitly uses the following two assumptions:
\begin{enumerate}[label=\arabic*)]
	\item the proof of the ``insertion-deletion" lemma implicitly uses the ``one-at-a-time" assumption. We should note that it is not difficult to generalize the proof of the ``deletion" case \cite[Lemma 17, Case 2]{tibshirani2013} to a multiple-deletion scenario, but if one wants to generalize the proof of the ``insertion" case, one needs to prove that 
	\[\brackets{X^\Tr_{i^{\text{join}}_{k+1}}(I-P)X^\Tr_{i^{\text{join}}_{k+1}}}\]
	appearing in the RHS of \cite[the last equation in p. 1485]{tibshirani2013} is invertible. As far as the authors see, this does not necessarily hold in the multiple-deletion scenario.
	
	\item the second part of the induction step indeed only proves that $\wEQ(\sv_+;\yv_0,\lambda)$ satisfies the optimality condition at $\lambda=\lambda_+$. To prove that the optimality condition holds for every $\lambda$ within the new linear piece, we need to assume that there exists a neighborhood $U$ of $\lambda_+$, such that $\wEQ(\sv_+;\yv_0,\lambda)$ satisfies the optimality condition for every $\lambda\in U$. This leads to a similar assumption to the ``sticky-equality" condition (where $\wstar(t)$ in Assumption \ref{assump:sticky_equality} is replaced by $\wEQ(\sv_+;\yv_0,\lambda)$ here).
\end{enumerate}
We should note that the ``one-at-a-time" assumption is intuitively a mild condition (see Remark \ref{remark:compare_assumptions}), but the mildness of the second assumption above is not obvious from \cite{tibshirani2013}.

Finally, we analyze how one can adapt the analysis of \cite{tibshirani2013} to obtain partial results of the current study. Since the optimality condition (\ref{eq:OPT}) of the sGMC model is structurally similar to the KKT condition of LASSO, it is possible to apply the induction analysis above to the sGMC model. However, due to the presence of the matrix $\Dm$ in (\ref{eq:OPT}), some efforts need to be made to facilitate the adaptation, including
\begin{enumerate}[label=\arabic*)]
	\item before adapting the induction analysis of LARS, we need to prove Thm. \ref{thm:shape_uniqueness_sparseness} by Appx. \ref{proof:geometry_and_uniqueness} to show that the sGMC model enjoys similar geometric properties to LASSO,
	
	\item when adapting the proof of the ``insertion-deletion" lemma, we need to prove Lemma \ref{lemma:blk_matrices} in Appx. \ref{app:technical_lemma} to simplify the analysis of the block-structured linear system therein.
\end{enumerate}
By the aforementioned efforts, we can conduct an induction analysis to a special implementation of E-LARS for computing the min-norm sGMC regularization path $\wstar(\lambda)$, whereby we can prove its continuity, piecewise simple behavior and continuous-discrete correspondence under the two assumptions introduced above.

\section{Proof of Theorem \ref{thm:continuity}}
\label{sec:proof:continuity}

We first introduce Lemma \ref{lemma:continuity_of_composition}, which is the core of our proof.

\begin{lemma}[Composition lemma]\label{lemma:continuity_of_composition}
	For $D_T\subset\Real^{p}$, $D_M\subset\Real^{m}$, suppose that the following holds for single-valued mapping $T:D_T\to \Real^{m}$ and set-valued mapping $\Mcal:D_M\rightrightarrows\Real^n$:
	\begin{enumerate}[label=(\alph*)]
		\item $T$ is continuous relative to $D_T$,
		\item For every $\thetav\in D_T$, $T(\thetav)\in D_M$.
		\item $\Mcal$ is osc (resp. isc) relative to $D_M$.
	\end{enumerate}
	Then $\Mcal\circ T(\cdot)\coloneqq \Mcal(T(\cdot))$ is osc (resp. isc) relative to $D_T$.
\end{lemma}
\begin{IEEEproof}
	For $\bar{\thetav}\in D_T$, let $(\thetav_k)_{k\in\Natural}\subset D_T$ be an arbitrary sequence convergent to $\bar{\thetav}$. Then from the assumptions (a) and (b), we have $T(\thetav_k)\in D_M$ for every $k$ and $\lim_{k\to\infty}T(\thetav_k)=T(\bar{\thetav})\in D_M$. If $\Mcal$ is osc relative to $D_M$, we have
	\begin{align*}
		\limsup_{k\to\infty}\Mcal(T(\thetav_k))\subset \Mcal(T(\bar{\thetav}))
	\end{align*}
	 from Def. \ref{def:continuity_of_setvalued_mapping}. Hence $\Mcal\circ T(\cdot)$ is osc relative to $D_T$. Similarly, if $\Mcal$ is isc/continuous relative to $D_T$, then $\Mcal\circ T(\cdot)$ is isc/continuous relative to $D_T$.
\end{IEEEproof}

Lemma \ref{lemma:continuity_of_composition} indicates that composition with a continuous parameterization function preserves semicontinuity of a set-valued mapping. This inspires a ``repeated-decomposition" technique for studying set-valued continuity, i.e., one can decompose complicated set-valued mappings into simpler ones, and transfer continuity properties from them. Using this technique, we can prove Theorem \ref{thm:continuity} by three steps:
\begin{enumerate}[label=\arabic*)]
	\item \textbf{Step 1:} we verify outer-semicontinuity of $\Se$ via the optimality condition (\ref{eq:OPT}). Exploiting the definition of outer-semicontinuity, we can further prove continuity of $\betae$ and $\gammae$, whereby Thm. \ref{thm:continuity} \ref{it:continuity:para} is obtained.
	
	\item \textbf{Step 2:} we decompose $\Se$ as a paramemterized LASSO solution set $\SLA$ via (\ref{eq:Se_variant_expression}), and further decompose $\SLA$ as a parameterized solution set of a linear program. The parameterization functions therein involve $\betae$ and $\gammae$, whose continuity has been established in the step 1.
	
	\item \textbf{Step 3:} we prove inner-semicontinuity of the solution set of linear programs using results from a previous work, then transfer inner-semicontinuity backward to $\Se$ via Lemma \ref{lemma:continuity_of_composition}. Continuity of $P_{\Se}(\bar{\wv})$ follows from continuity of $\Se$, whereby Thm. \ref{thm:continuity} \ref{it:continuity:lmd} is obtained.
\end{enumerate}
Combining the three steps above completes the proof.

\subsection{Step 1: Verify Outer-Semicontinuity of $\Se$ and Identify Continuous Parameterization Functions Therein}

In the step 1, we consider $\Se$ as a set-valued mapping, and consider $\betae$, $\gammae$ as single-valued mappings of $\paragroup$.

The outer semicontinuity of $\Se$ with respect to $\paragroup$ can be derived directly from (\ref{eq:OPT}) and Thm. \ref{thm:shape_uniqueness_sparseness}; see the following proposition.

\begin{proposition}\label{prop:outer_semicontinuity_of_Se}
	For an arbitrary sequence of parameters $(\Am_k,\rho_k,\bv_k,\lambda_k)_{k\in\Natural}\subset\Pscr$ convergent to $\paragroup\in\Pscr$, let $(\wv_k)_{k\in\Natural}\subset\Real^{2n}$ be a sequence of points satisfying
	\[\wv_k\in\Se(\Am_k,\rho_k,\bv_k,\lambda_k),\]
	then the following holds:
	\begin{enumerate}[label=(\alph*)]
		\item $(\wv_k)_k$ is bounded.
		
		\item For every cluster point of $(\wv_k)_k$, say $\wv_{\infty}$, we have $\wv_{\infty}\in\Se(\Am,\rho,\bv,\lambda)$. In other words, $\Se$ is osc relative to $\Pscr$.
	\end{enumerate}
	
\end{proposition}
\begin{IEEEproof}
	(a) From Thm. \ref{thm:shape_uniqueness_sparseness} \ref{it:shape:bound}, we have
	\[(\forall k\in\Natural)\;\;\norm{\wv_k}_1\leq \frac{1}{\lambda_k(1-\rho_k)}\norm{\yv_k}^2_2+\frac{\norm{\rv_k}^2_2}{\lambda_k}.\]
	Since $(\Am_k,\rho_k,\bv_k,\lambda_k)$ converges to $\paragroup$, for large enough $k$ we have
	\[\norm{\wv_k}_1\leq \frac{1}{\lambda(1-\rho)}\norm{\yv}^2_2+\frac{\norm{\rv}^2_2}{\lambda}+1,\]
	which verifies the boundedness of $(\wv_k)_{k\in\Natural}$. 
	
	\vspace{0.5em}
	\noindent
	(b) The boundedness of $(\wv_k)_{k\in\Natural}$ guarantees the existence of a cluster point $\wv_{\infty}$. Hence by passing to a subsequence, we can assume $\lim_{k\to\infty}\wv_k=\wv_{\infty}$ without loss of generality. From (\ref{eq:OPT}), for every $k$ and $i\in\{1,\dots,2n\}$, we have
	\begin{equation}\label{eq:convergence_of_extended_solution_temp1}
		\cv_{k,i}(\bv_k-\Dm_k\Cm_k\wv_k)\in[-\lambda_k,\lambda_k],
	\end{equation}
	where $(\Dm_k,\Cm_k)$ are the counterpart parts of $(\Dm,\Cm)$ in (\ref{eq:OPT}) for $\Se(\Am_k,\rho_k,\bv_k,\lambda_k)$, and $\cv_{k,i}$ is the $i$th column vector of $\Cm_k$. Taking $k\to\infty$ on both sides of (\ref{eq:convergence_of_extended_solution_temp1}) yields
	\begin{equation}\label{eq:convergence_of_extended_solution_temp2}
		\cv_i^{\Tr}(\bv-\Dm\Cm\wv_{\infty})\in[-\lambda,\lambda].
	\end{equation}
	For every $i\in\{1,\dots,2n\}$ such that $[\wv_{\infty}]_i\neq 0$, we have $\sign([\wv_k]_i)=\sign([\wv_{\infty}]_i)$ for large enough $k$, hence
	\begin{equation*}
		\cv_{k,i}^{\Tr}(\bv_k-\Dm_k\Cm_k\wv_k)=\lambda_k\,\sign([\wv_k]_i)=\lambda_k\,\sign([\wv_\infty]_i)
	\end{equation*}
	from (\ref{eq:OPT}). Taking $k\to\infty$ on both sides of the equation above yields $\cv_i^{\Tr}(\bv-\Dm\Cm\wv_{\infty})=\lambda\,\sign([\wv_{\infty}]_i)$. Combining this with (\ref{eq:convergence_of_extended_solution_temp2}), one can verify that $\wv_{\infty}$ satisfies (\ref{eq:OPT}) with parameters $\paragroup$, hence $\wv_{\infty}\in\Se\paragroup$.
\end{IEEEproof}

Exploiting Prop. \ref{prop:outer_semicontinuity_of_Se}, we can further establish single-valued continuity of $\betae$ and $\gammae$, which will serve as parameterization functions for the decomposition of $\Se$ in the step 2.

\begin{corollary}\label{cor:coninuity_of_betae_gammae}
	Single-valued mappings $\betae$ and $\gammae$ are both continuous with respect to $\paragroup$ relative to $\Pscr$.
\end{corollary}
\begin{IEEEproof}
	For $\paragroup\in\Pscr$, let $(\Am_k,\rho_k,\bv_k,\lambda_k)_{k\in\Natural}$ be an arbitrary sequence in $\Pscr$ convergent to $\paragroup$. Note that for every $k\in\Natural$, we can find some $\wv_k\in\Se(\Am_k,\rho_k,\bv_k,\lambda_k)$, hence from Thm. \ref{thm:shape_uniqueness_sparseness} \ref{it:shape:linear_fit} we have
	\begin{align*}
		\betav_k &\coloneqq\betae(\Am_k,\rho_k,\bv_k,\lambda_k)=\Cm_k\wv_k,\\
		\gamma_k &\coloneqq\gammae(\Am_k,\rho_k,\bv_k,\lambda_k)=\norm{\wv_k}_1.
	\end{align*}
	where $\Cm_k\coloneqq\blkdiag(\Am_k,\sqrt{\rho}\Am_k)$. We further define
	\begin{align*}
		\betav_{\infty} &\coloneqq  \betae\paragroup, &\gamma_{\infty} &\coloneqq \gammae\paragroup,\\
		p_k&\coloneqq \norm{\betav_k-\betav_{\infty}}_2, &q_k &\coloneqq \left\lvert\gamma_k-\gamma_{\infty}\right\rvert
	\end{align*}
	Then to prove Cor. \ref{cor:coninuity_of_betae_gammae}, we only need to prove
	\[\lim_{k\to\infty} p_k=\lim_{k\to\infty} q_k=0.\] Since $p_k\geq 0,q_k\geq 0$, it is sufficient to prove
	\begin{equation*}
		\limsup_{k\to\infty}p_k=\limsup_{k\to\infty}q_k=0.
	\end{equation*}
	
	From Fact \ref{fact:limsup_liminf} in Appx. \ref{app:other_useful_facts}, there exists a subsequence $(p_{k_i})_{i\in\Natural}$ of $(p_k)_{k\in\Natural}$ such that $\lim_{i\to\infty}p_{k_i}=\limsup_{k\to\infty}p_k$. Since $(\wv_{k_i})_{i\in\Natural}$ is bounded from Prop. \ref{prop:outer_semicontinuity_of_Se} (a), by further passing to a subsequence (if necessary), we can assume that $(\wv_{k_i})_{i\in\Natural}$ is convergent and $\lim_{i\to\infty}\wv_{k_i}=\wv_{p,\infty}$. Then from Prop. \ref{prop:outer_semicontinuity_of_Se} (b) we have $\wv_{p,\infty}\in\Se\paragroup$, which implies $\betav_\infty=\Cm\wv_{p,\infty}$ with $\Cm\coloneqq \blkdiag(\Am,\sqrt{\rho}\Am)$. Accordingly, we have
	\begin{align*}
		\limsup_{k\to\infty} p_k=\lim_{i\to\infty}p_{k_i}=\lim_{i\to\infty}\norm{\Cm_{k_i}\wv_{k_i}-\Cm\wv_{p,\infty}}_2=0.
	\end{align*}
	
	Similarly, there exists a subsequence $(q_{k_j})_{j\in\Natural}$ of $(q_k)_{k\in\Natural}$ such that $\lim_{j\to\infty}q_{k_j}=\limsup_{k\to\infty}q_k$. By further passing to a subsequence (if necessary), we can assume that $\lim_{j\to\infty}\wv_{k_j}=\wv_{q,\infty}$ for some $\wv_{q,\infty}\in\Se\paragroup$. Hence we have $\gamma_{\infty}=\norm{\wv_{q,\infty}}_1$ and
	\begin{align*}
		\limsup_{k\to\infty} q_k=\lim_{j\to\infty}q_{k_j}=\lim_{k\to\infty}\left\lvert\norm{\wv_{k_j}}_1-\norm{\wv_{q,\infty}}_1\right\rvert=0.
	\end{align*}
	Combining the discussion above completes the proof.
\end{IEEEproof}

Combining Prop. \ref{prop:outer_semicontinuity_of_Se} and Cor. \ref{cor:coninuity_of_betae_gammae} proves Thm. \ref{thm:continuity} \ref{it:continuity:para}.

\subsection{Step 2: Repeatedly Decompose $\Se$ into Parameterized Versions of Simpler Set-Valued Mappings}

In the step 2 and 3, for fixed $(\Am,\rho)\in\Real^{m\times n}\times[0,1)$, we consider $\Se$ as a set-valued mapping and consider $\betae$, $\gammae$, $P_{\Se}(\bar{\wv})$ as single-valued mappings of $(\bv,\lambda)\in\Real^{2m}\times\Real_{++}$.

We first decompose $\Se(\bv,\lambda)$ via the LASSO solution set mapping $\SLA(\cdot,\cdot,\cdot)$ (defined in the paragraph after (\ref{eq:LASSO})). From (\ref{eq:Se_variant_expression}) in Thm. \ref{thm:shape_uniqueness_sparseness}, for every $\lambda\in\Real_{++}$ we have
\begin{align*}
	\Se(\lambda)=\SLA\parens{\Am,\frac{\yv-\rho\betad(\bv,\lambda)}{1-\rho},\frac{\lambda}{1-\rho}}\times \SLA(\sqrt{\rho}\Am,\rv+\sqrt{\rho}\betap(\bv,\lambda),\lambda),
\end{align*}
where $\betap\coloneqq \Am\xv$ with $\xv\in\Sp$, $\betad\coloneqq \Am\zv$ with $\zv\in\Sd$ (Thm. \ref{thm:shape_uniqueness_sparseness} \ref{it:shape:linear_fit}). Let us define $D_T\coloneqq \Real^{2m}\times\Real_{++}$, $D_M\coloneqq \Real^{m}\times \Real_{++}$ and the following mappings:
\begin{align*}
	T_1&:D_T\to D_M: (\bv,\lambda)\mapsto \parens{\frac{\yv-\rho\betad(\bv,\lambda)}{1-\rho},\frac{\lambda}{1-\rho}},\\
	T_2&:D_T\to D_M:(\bv,\lambda)\mapsto \parens{\rv+\sqrt{\rho}\betap(\bv,\lambda),\lambda},\\
	\Mcal_1&:D_M\rightrightarrows \Real^n:(\zetav,\eta)\mapsto \SLA(\Am,\zetav,\eta),\\
	\Mcal_2&:D_M\rightrightarrows\Real^n:(\zetav,\eta)\mapsto\SLA(\sqrt{\rho}\Am,\zetav,\eta).
\end{align*}
Then one can verify that 
\begin{equation}\label{temp:continuity_of_Se}
	\Se(\yv,\lambda)=\Mcal_1(T_1(\yv,\lambda))\times\Mcal_2(T_2(\yv,\lambda)).
\end{equation}
Moreover, one can verify that the decomposition above satisfies the first two conditions required by Lemma \ref{lemma:continuity_of_composition}:
\begin{enumerate}[label=(\alph*)]
	\item $T_1$ and $T_2$ are both continuous relative to $D_T\coloneqq\Real^{2m}\times\Real_{++}$. This can be verified by continuity of $\betae$ (Cor. \ref{cor:coninuity_of_betae_gammae}) and recalling $\betae=\begin{bmatrix}
		\betap^\Tr & \sqrt{\rho}\betad^\Tr
	\end{bmatrix}^\Tr$ from (\ref{eq:Se_SpSd}).
	
	\item For every $(\bv,\lambda)\in\Real^{2m}\times\Real_{++}$, $T_1(\bv,\lambda),T_2(\bv,\lambda)\in D_M\coloneqq \Real^m\times\Real_{++}$ from the definitions of $T_1$, $T_2$.
\end{enumerate}
Thus from Lemma \ref{lemma:continuity_of_composition}, inner-semicontinuity of $\Se(\yv,\lambda)$ can be transferred from that of $\Mcal_1$ and $\Mcal_2$, i.e., we only need to prove that for a fixed sensing matrix $\Am\in\Real^{m\times n}$, $\SLA(\Am,\yv,\lambda)$ (abbreviated as $\SLA(\yv,\lambda)$) is isc with respect to $(\yv,\lambda)$.

Next we further decompose $\SLA(\yv,\lambda)$ into a parameterized solution set of a linear program. For $(\yv,\lambda)\in\Real^m\times\Real_{++}$, let $\xv_{\text{p}}\in\SLA(\yv,\lambda)$, then $\xv\in\Real^n$ is in $\SLA(\yv,\lambda)$ if and only if
\begin{equation}\label{eq:optimality_LASSO}
	\frac{1}{2}\norm{\yv-\Am\xv}^2_2+\lambda\norm{\xv}_1\leq\frac{1}{2}\norm{\yv-\Am\xv_{\text{p}}}^2_2+\lambda\norm{\xv_{\text{p}}}_1.
\end{equation}
Let $\betaLA$ and $\gammaLA$ respectively be the common linear fit and common $\ell_1$-norm of LASSO solutions (Thm. \ref{thm:shape_uniqueness_sparseness} \ref{it:shape:linear_fit}), i.e.,
\begin{equation*}
	\betaLA\coloneqq \Am\xv_{\text{p}},\;\; \gammaLA\coloneqq \norm{\xv_{\text{p}}}_1\text{ with }\xv_{\text{p}}\in\SLA(\yv,\lambda),
\end{equation*}
then the optimality condition (\ref{eq:optimality_LASSO}) is further equivalent to
\begin{equation}\label{temp:continuity_of_S_LA1}
	\Am\xv=\betaLA(\yv,\lambda),\; \norm{\xv}_1\leq\gammaLA(\yv,\lambda),
\end{equation}
Notice that the $\ell_1$-norm can be rewritten as
\begin{equation}\label{temp:continuity_of_S_LA2}
	\norm{\xv}_1=\sum_{i=1}^n \sign([\xv]_i)[\xv]_i=\max_{\qv\in\{\pm 1\}^n}\qv^{\Tr}\xv,
\end{equation}
by substituting (\ref{temp:continuity_of_S_LA2}) into (\ref{temp:continuity_of_S_LA1}), one can verify that $\xv\in\SLA(\yv,\lambda)$ is equivalent to
\begin{equation}\label{temp:continuity_of_S_LA3}
	\Am\xv=\betaLA(\yv,\lambda),\; \Qm\xv\leq \gammaLA(\yv,\lambda)\mathbf{1}_{2^n},
\end{equation}
where $\Qm\in\Real^{2^n\times n}$ is the matrix with its rows being all possible vectors in $\{\pm 1\}^n$. Let us define $\Lcal_{\Am,\Qm}:\Real^m\times\Real^{2^n}\rightrightarrows \Real^n$ as the solution set mapping of a linear program:
\begin{align}
	\Lcal_{\Am,\Qm}(\zetav,\etav)\coloneqq \braces{\xv\in\Real^n\;\middle\vert\; \Am\xv=\zetav,\Qm\xv\leq \etav}.\label{eq:LAQ_def}
\end{align}
Moreover, define the following mappings:
\begin{align*}
	T_3&:\Real^m\times\Real_{++}\to\Real^{m}\times\Real^{2^n}:(\yv,\lambda)\mapsto \parens{\betaLA(\yv,\lambda),\gammaLA(\yv,\lambda)\mathbf{1}_{2^n}},\\
	\Mcal_3 &: \dom\Lcal_{\Am,\Qm}\rightrightarrows \Real^n: (\zetav,\etav)\mapsto \Lcal_{\Am,\Qm}(\zetav,\etav), 
\end{align*}
Then (\ref{temp:continuity_of_S_LA3}) yields the following decomposition of $\SLA$:
\begin{equation}\label{eq:decompose_SLA}
	\SLA(\yv,\lambda)=\Mcal_3(T_3(\yv,\lambda)).
\end{equation}
One can verify that the decomposition above satisfies the first two conditions required by Lemma \ref{lemma:continuity_of_composition}:
\begin{enumerate}[label=(\alph*)]
	\item $T_3$ is continuous relative to $\Real^m\times\Real_{++}$. This can be verified by continuity of $\betae$ (Cor. \ref{cor:coninuity_of_betae_gammae}) in the LASSO case.
	
	\item For every $(\yv,\lambda)\in\Real^m\times\Real_{++}$, $T_3(\yv,\lambda)\in\dom\Lcal_{\Am,\Qm}$. This follows from (\ref{temp:continuity_of_S_LA3}) and nonemptiness of $\SLA(\yv,\lambda)$.
\end{enumerate}
Thus from Lemma \ref{lemma:continuity_of_composition}, inner-semicontinuity of $\SLA(\yv,\lambda)$ can be transferred from that of $\Mcal_3$, i.e., we only need to prove that $\Lcal_{\Am,\Qm}(\cdot,\cdot)$ is isc relative to $\dom\Lcal_{\Am,\Qm}$.

\subsection{Step 3: Establish Inner-Semicontinuity of the Simplest Set-Valued Mapping and Transfer the Continuity Backward}

Inner-semicontinuity of the solution sets of linear programs can be implied from a previous work (see Fact \ref{fact:Lipschitzian_continuity_of_linear_program} in Appx. \ref{app:other_useful_facts}); see the following lemma.
\begin{lemma}\label{lemma:continuity_of_linear_programs}
	Given $\Am\in\Real^{m\times n}$, $\Qm\in\Real^{k\times n}$, for $\zetav\in\Real^m$ and $\etav\in\Real^k$, define $\Lcal_{\Am,\Qm}:\Real^m\times\Real^k\rightrightarrows \Real^n$ as:
	\begin{align*}
		\Lcal_{\Am,\Qm}(\zetav,\etav)\coloneqq \braces{\xv\in\Real^n\;\middle\vert\; \Am\xv=\zetav,\Qm\xv\leq \etav}.
	\end{align*}
	Then $\Lcal_{\Am,\Qm}(\cdot,\cdot)$ is continuous relative to its domain 
	$\dom\Lcal_{\Am,\Qm}\coloneqq \braces{(\zetav,\etav)\in\Real^m\times\Real^k \;\middle\vert\; \Lcal_{\Am,\Qm}(\zetav,\etav)\neq\emptyset}$.
\end{lemma}
\begin{IEEEproof}
	For $\bar{\thetav}\coloneqq (\bar{\zetav},\bar{\etav})\in\dom\Lcal_{\Am,\Qm}$, let $(\thetav_k)_{k\in\Natural}\coloneqq(\zetav_k,\etav_k)_{k\in\Natural}\subset\dom\Lcal_{\Am,\Qm}$ be a sequence convergent to $\bar{\thetav}$.
	
	We first prove the outer semicontinuity of $\Lcal_{\Am,\Qm}(\cdot,\cdot)$. Let $\bar{\xv}$ be an arbitrary point in $\limsup_{k}\Lcal_{\Am,\Qm}(\thetav_k)$, then from Def. \ref{def:limit_of_sets}, $\bar{\xv}$ is a cluster point of some sequence $(\xv_k)_{k\in\Natural}$ satisfying $\xv_k\in\Lcal_{\Am,\Qm}(\thetav_k)$. Hence there exists a subsequence of $(\xv_k)_{k\in\Natural}$, say  $(\xv_{k_i})_{i\in\Natural}$, that converges to $\bar{\xv}$. Since $\xv_{k_i}\in\Lcal_{\Am,\Qm}(\thetav_{k_i})$, we have $\Am\xv_{k_i}=\zetav_{k_i}$ and $\Qm\xv_{k_i}\leq\etav_{k_i}$.
	Taking the limit $i\to\infty$ on the equation and inequation yields
	$\Am\bar{\xv}= \bar{\zetav}$ and $\Qm\bar{\xv}\leq \bar{\etav}$,
	which implies that $\bar{\xv}\in \Lcal_{\Am,\Qm}(\bar{\thetav})$. Thus we have $\limsup_k\Lcal_{\Am,\Qm}(\thetav_k)\subset\Lcal_{\Am,\Qm}(\bar{\thetav})$. From Def. \ref{def:continuity_of_setvalued_mapping}, we conclude that $\Lcal_{\Am,\Qm}(\cdot,\cdot)$ is osc relative to $\dom\Lcal_{\Am,\Qm}$.
	
	Next we prove the inner semicontinuity of $\Lcal_{\Am,\Qm}(\cdot,\cdot)$. Let $\bar{\xv}$ be an arbitrary point in $\Lcal_{\Am,\Qm}(\bar{\thetav})$. Then by Fact \ref{fact:Lipschitzian_continuity_of_linear_program} in Appx. \ref{app:other_useful_facts}, for every $k$, we can find $\xv_k\in\Lcal_{\Am,\Qm}(\thetav_k)$ such that
	\begin{equation*}
		\norm{\bar{\xv}-\xv_k}_{\infty}\leq \mu_{\Lcal}\norm{\begin{bmatrix}
				\bar{\zetav}-\zetav_k \\ \bar{\etav}-\etav_k
		\end{bmatrix}}_{\Lcal}.
	\end{equation*}
	Taking the limit $k\to\infty$ on the inequation above, one can verify that $\lim_{k\to\infty}\xv_k=\bar{\xv}$. Hence from Def. \ref{def:limit_of_sets}, we have $\bar{\xv}\in\liminf_k\Lcal_{\Am,\Qm}(\bar{\thetav})$, which further implies that $\Lcal_{\Am,\Qm}(\bar{\thetav})\subset \liminf_k\Lcal_{\Am,\Qm}(\thetav_k)$. From Def. \ref{def:continuity_of_setvalued_mapping}, we conclude that $\Lcal_{\Am,\Qm}(\cdot,\cdot)$ is isc relative to $\dom\Lcal_{\Am,\Qm}$.
\end{IEEEproof}

Thus by applying Lemma \ref{lemma:continuity_of_composition} to (\ref{temp:continuity_of_Se}) and (\ref{eq:decompose_SLA}), we can derive inner-semicontinuity of $\Se(\bv,\lambda)$, combining which with outer-semicontinuity of $\Se$ (Prop. \ref{prop:outer_semicontinuity_of_Se}) further implies that $\Se(\bv,\lambda)$ is continuous relative to $\Real^{2m}\times\Real_{++}$. Continuity of $P_{\Se}(\bar{\wv})$ with respect to $(\bv,\lambda)\in\Real^{2m}\times\Real_{++}$ follows directly from the set-valued continuity above (see Fact \ref{fact:continuity_of_projection} in Appx. \ref{app:other_useful_facts}). Combining the aforementioned results completes the proof of Thm. \ref{thm:continuity} \ref{it:continuity:lmd}.

\section{Proof of Theorem \ref{thm:correspondence}}
\label{sec:proof:correspondence}

We first introduce a discrete inequality (i.e., (\ref{eq:s_inequality}) below) that relates a candidate indicator $\sv$ to two discrete quantities $\sign(\wstar(\bv,\lambda))$ and $\se(\bv,\lambda)$, which is the core of our proof.
\begin{definition}\label{def:partial_order_of_s}
	For $\sv_1,\sv_2\in\braces{+1,0,-1}^{2n}$, let us define $\Ecal_1\coloneqq\supp(\sv_1)$ and $\Ecal_2\coloneqq\supp(\sv_2)$. We say $\sv_1\preceq\sv_2$ if
	\begin{equation*}
		\Ecal_1\subset\Ecal_2\text{ and }[\sv_1]_{\Ecal_1}=[\sv_2]_{\Ecal_1},
	\end{equation*}
	i.e., if the components of $\sv_1$ and $\sv_2$ coincide on the smaller support $\Ecal_1$ of $\sv_1$. Then $\preceq$ is a partial order on ${\Vcal}$.
\end{definition}
\begin{proposition}\label{prop:partial_order}
	For every $\sv\in\braces{+1,0,-1}^{2n}$ and $(\bv,\lambda)\in\Real^{2m}\times\Real_{++}$, the following holds:
	\begin{enumerate}[label=(\alph*)]	
		\item \label{it:order:inequality}
		Given $\wv\in\Se(\bv,\lambda)$, then we always have
		\begin{equation}\label{eq:inequality_in_Se}
			\sign(\wv)\preceq \se(\bv,\lambda),
		\end{equation}
		Moreover, $\sv$ satisfies the following inequality 
		\begin{equation}\label{eq:s_inequality}
			\sign(\wv)\preceq \sv\preceq \se(\bv,\lambda)
		\end{equation}
		if and only if $\wv\in\Spoly(\sv;\bv,\lambda)$.
		
		\item \label{it:order:equality} Partial equality of (\ref{eq:s_inequality}) leads to some favorable equalities:
		\begin{itemize}
			\item if $\sv=\sign(\wstar(\bv,\lambda))$, then we have 
			\begin{equation}\label{eq:wstar_equal_wEQ}
				\wstar(\bv,\lambda)=\wEQ(\sv;\bv,\lambda),
			\end{equation}			
			
			\item if $\sv=\se(\bv,\lambda)$, then we have
			\begin{equation}\label{eq:Se_equal_SEQNQ}
				\Se(\bv,\lambda)=\Spoly(\sv;\bv,\lambda).
			\end{equation}			 
		\end{itemize}
	\end{enumerate}
\end{proposition}

The proof of Prop. \ref{prop:partial_order} can be found in Appx. \ref{app:auxiliary_results}, along with the proofs of other two auxiliary results Lemma \ref{lemma:SEQNQ_property} and \ref{lemma:IEN_property} in Sec. \ref{sec:correspondence}. Prop. \ref{prop:partial_order} provides a powerful tool for establishing conditional equality between $\wstar(\yv,\lambda)$ and affine function $\wEQ(\sv;\yv,\lambda)$. Exploiting the conditions that imply partial equalities of (\ref{eq:s_inequality}), we can prove Theorem \ref{thm:correspondence}. Our proof is composed of three steps:
\begin{enumerate}[label=\arabic*)]
	\item \textbf{Step 1:} exploiting the condition that leads to equality between $\wstar(\bv,\lambda)$ and $\wEQ(\sv;\bv,\lambda)$ (Prop. \ref{prop:partial_order} \ref{it:order:equality}), we can introduce a partition of the parameter space $\Real^{2m}\times\Real_{++}$, combining which with linearity of $\wEQ(\sv;\cdot,\cdot)$ yields piecewise linearity of $\wstar(\cdot,\cdot)$, whereby proving partial result of Thm. \ref{thm:correspondence} \ref{it:corres:PL}.
	
	\item \textbf{Step 2:} we further study the behavior of $\wstar(\bv,\lambda)$ within a set $\bar{I}(\sv)$ from the partition of $\Real^{2m}\times\Real_{++}$ introduced in step 1. Exploiting the properties of $\wEQ(\sv;\bv,\lambda)$ and $\Iselect(\sv)$ (Lemma \ref{lemma:IEN_property} in Sec. \ref{sec:correspondence}), we prove that such $\bar{I}(\sv)$ is a convex cone, and is exactly a linear zone in $\wstar(\cdot,\cdot)$, whereby proving the result \ref{it:corres:zone} and unproven part of \ref{it:corres:zone}.
	
	\item \textbf{Step 3:} we show that if $\Am$ has columns in general position, then for every zone indicator $\sv$, $\sign(\wEQ(\sv;\bv,\lambda))$ equals to $\se(\bv,\lambda)$ for every $(\bv,\lambda)\in\interior(\Iselect(\sv))$, combining which with Prop. \ref{prop:partial_order} \ref{it:order:equality} yields the result \ref{it:corres:GP}.
\end{enumerate}

\subsection{Step 1: Establish Piecewise Linearity of the Min-Norm Solution Map by a Partition of Parameter Space}
We first prove piecewise linearity of $\wstar(\bv,\lambda)$, i.e., part of Thm. \ref{thm:correspondence} \ref{it:corres:PL}, by giving a finite partition of $\Real^{2m}\times\Real_{++}$ where $\wstar(\cdot,\cdot)$ is linear within each sets therein. Define
\begin{equation}\label{eq:domain}
	\breve{I}(\sv)\coloneqq\braces{(\bv,\lambda)\in\Real^{2m}\times\Real_{++}\;\middle\vert\; \sv=\sign(\wstar(\bv,\lambda))}
\end{equation}
for $\sv\in\braces{+1,0,-1}^{2n}$. Then from (\ref{eq:wstar_equal_wEQ}) in Prop. \ref{prop:partial_order} \ref{it:order:equality}, we have the following result for every $\sv\in\braces{+1,0,-1}^{2n}$:
\begin{align}\label{eq:breve_Is_PL}
	(\forall(\bv,\lambda)\in \breve{I}(\sv))\quad\wstar(\bv,\lambda)=\wEQ(\sv;\bv,\lambda),
\end{align}
which indicates that every nonempty $\breve{I}(\sv)$ is a region where $\wstar(\cdot,\cdot)$ is linear with $\sign(\wstar(\bv,\lambda))\equiv \sv$ being constant.

Moreover, from the definition (\ref{eq:domain}) of $\breve{I}(\sv)$, it is evident that 
\begin{equation}\label{eq:I_s_do_not_intersect}
	\sv_1\neq\sv_2\implies\breve{I}(\sv_1)\cap\breve{I}(\sv_2)=\emptyset,
\end{equation}
and the union of $\breve{I}(\sv)$ for all candidate indicators $\sv$ gives a cover of the whole parameter space $\Real^{2m}\times\Real_{++}$, i.e.,
\begin{equation}\label{eq:I_s_cover_space}
	\parens{\cup_{\sv\in\braces{+1,0,-1}^{2n}} \breve{I}(\sv)}\supset \parens{\Real^{2m}\times\Real_{++}}.
\end{equation}
Combining (\ref{eq:I_s_do_not_intersect}) and (\ref{eq:I_s_cover_space}), we can deduce that the union of all $\breve{I}(\sv)$ defines a partition of the parameter space $\Real^{2m}\times\Real_{++}$. Notice that there only exist a finite number of $\breve{I}(\sv)$, thus from the discussion just after (\ref{eq:breve_Is_PL}), we can directly deduce piecewise linearity of $\wstar(\cdot,\cdot)$.

In the sequel, we further simplify the aforementioned partition of the parameter space $\Real^{2m}\times\Real_{++}$ by $\breve{I}(\sv)$, which will be useful in subsequent analysis. Notice that for every candidate indicator $\sv$, either of the following two cases happens:
\begin{enumerate}[label=\arabic*)]
	\item $\interior(\breve{I}(\sv))\neq\emptyset$, in which case $\breve{I}(\sv)\subset\Real^{2m}\times\Real_{++}$ has positive Lebesgue measure, and must contain\footnote{We note that so far, the connectedness of $\breve{I}(\sv)$ has not been proven yet, thus it can be the union of interior volumes of multiple linear zones.} a part of interior of some linear zone in $\wstar(\cdot,\cdot)$,
	
	\item $\interior(\breve{I}(\sv))=\emptyset$, in which case $I(\sv)$ is \textit{nowhere dense}. One can imagine that such $\breve{I}(\sv)$ is either empty or the union of some vertices/edges/faces of linear zones in $\wstar(\cdot,\cdot)$.
\end{enumerate}
Since there are finite possible values for $\sv\in\braces{+1,0,-1}^{2n}$, the union of $\breve{I}(\sv)$ for all $\sv$ satisfying $\interior(\breve{I}(\sv))=\emptyset$ gives a \textit{meager} set, which implies that the collection of all $\breve{I}(\sv)$ satisfying $\interior(\breve{I}(\sv))\neq\emptyset$ defines a partition of a dense subset of $\Real^{2m}\times\Real_{++}$ (Baire category theorem). From this observation, in the following analysis, we only need to consider candidate indicators $\sv$ satisfying $\interior(\breve{I}(\sv))\neq\emptyset$.

\subsection{Step 2: Establish the One-to-One Correspondence between Linear Zones and Their Indicators}

Next we prove Thm. \ref{thm:correspondence} \ref{it:corres:zone} and unproven results of Thm. \ref{thm:correspondence} \ref{it:corres:PL} (i.e., piecewise constant sparsity pattern of $\wstar(\cdot,\cdot)$), by studying the behavior of $\wstar(\bv,\lambda)$ within some $\breve{I}(\sv)$ satisfying $\interior(\breve{I}(\sv))\neq\emptyset$. Let us define $\bar{I}(\sv)$ as the closure of $\breve{I}(\sv)$ in $\Real^{2m}\times\Real_{++}$, i.e.,
\begin{equation}
	\bar{I}(\sv)\coloneqq \closure(\breve{I}\parens{\sv)}\cap\parens{\Real^{2m}\times\Real_{++}}.
\end{equation}
Hereafter we assume $\interior(\bar{I}(\sv))\neq\emptyset$. Then from (\ref{eq:breve_Is_PL}) and continuity of $\wstar(\cdot,\cdot)$ (Thm. \ref{thm:continuity}), we have 
\begin{align}\label{eq:bar_I_PL}
	(\forall(\bv,\lambda)\in \bar{I}(\sv))\quad\wstar(\bv,\lambda)=\wEQ(\sv;\bv,\lambda).
\end{align}
In the sequel, we will further prove the following results:
\begin{enumerate}[label=\arabic*)]
	\item $\bar{I}(\sv)\subset\Iselect(\sv)$,
	
	\item $\bar{I}(\sv)$ is a convex cone, which implies that $\bar{I}(\sv)$ is a connected subset of some linear zone in $\wstar(\cdot,\cdot)$,
	
	\item there do not exist $\sv'\neq\sv$ satisfying $\interior(\breve{I}(\sv'))\neq\emptyset$, such that $\interior(\breve{I}(\sv))$ and $\interior(\breve{I}(\sv'))$ are contained in the interior of the same linear zone in $\wstar(\cdot,\cdot)$, which implies that $\bar{I}(\sv)$ is a single linear zone in $\wstar(\cdot,\cdot)$.
\end{enumerate}
Combining the aforementioned results yields Thm. \ref{thm:correspondence} \ref{it:corres:zone} and piecewise constant sparsity pattern of $\wstar(\cdot,\cdot)$. 

Now we prove the claimed results above in order.

\vspace{0.5em}
\noindent
(I) We first prove the inclusion $\bar{I}(\sv)\subset\Iselect(\sv)$. From Def. \ref{def:EQ_NQ_quantities}, one can verify that the closure of $\Iselect(\sv)$ in $\Real^{2m}\times\Real_{++}$ is exactly itself, thus to prove $\bar{I}(\sv)\subset\Iselect(\sv)$, it is sufficient to prove $\interior(\bar{I}(\sv))=\interior(\breve{I}(\sv))\subset\interior(\Iselect(\sv))$. This follows directly from (\ref{eq:breve_Is_PL}), i.e.,
\begin{align}
	(\forall (\bv,\lambda)\in \breve{I}(\sv))\quad \wEQ(\sv;\bv,\lambda)=\wstar(\bv,\lambda)\in\Se(\bv,\lambda), \label{temp:breve_Is_wEQ_in_Se}
\end{align}
as well as Lemma \ref{lemma:IEN_property} in Sec. \ref{sec:preparation_for_corres}.

\vspace{0.5em}
\noindent
(II) Next we prove that $\bar{I}(\sv)$ is a convex cone. It is sufficient to prove that $\breve{I}(\sv)$ is a convex cone.

From the result (I), we have $\breve{I}(\sv)\subset\Iselect(\sv)$. Then by Lemma \ref{lemma:IEN_property} in Sec. \ref{sec:preparation_for_corres}, for every $(\bv,\lambda)\in\interior(\Iselect(\sv))$,
\begin{align}
	\sign(\wEQ(\sv;\bv,\lambda))=\sign(\wstar(\bv,\lambda))\equiv \sv,\label{temp:sign_wEQ}\\
	\se(\bv,\lambda)=\sOPT(\wEQ(\sv;\bv,\lambda);\bv,\lambda)\equiv\sbar\succeq \sv,\label{temp:csign_wEQ}
\end{align}
where the first equality follows from (\ref{temp:breve_Is_wEQ_in_Se}), the inequality in (\ref{temp:csign_wEQ}) follows from (\ref{eq:inequality_in_Se}). For convenience of analysis, here after we define $\Ecal\coloneqq\supp(\sv)$ and $\Ebar\coloneqq\supp(\sbar)$. Combining (\ref{temp:csign_wEQ}) with (\ref{eq:Se_equal_SEQNQ}) in Prop. \ref{prop:partial_order}, we can deduce that
\begin{align*}
	(\forall(\bv,\lambda)\in\interior(\Iselect(\sv))) \quad \Se(\bv,\lambda)\equiv\Spoly(\sbar;\bv,\lambda).
\end{align*}
Thus for every $(\bv,\lambda)\in\interior(\Iselect(\sv))$, $\wstar(\bv,\lambda)$ is the unique solution to the following convex program:
\begin{equation*}
	\begin{aligned}
		\underset{\wv\in\Real^{2n}}{\text{minimize}}\quad & \frac{1}{2}\norm{\wv}^2_2\\
		\text{subject to}\quad & \wv\in\Spoly(\sbar;\bv,\lambda).
	\end{aligned}
\end{equation*}
By analyzing the KKT condition \cite[Sec. 5.5.3]{boyd2004} of the convex program above (the derivation is the same as (\ref{temp:KKT_analysis_first_line})-(\ref{temp:w_opt_necessary_conditions4b}) in Appx. \ref{proof:partial_order}), we can deduce that for every $(\bv,\lambda)\in\interior(\Iselect(\sv))$, $\wEQ(\sv;\bv,\lambda)=\wstar(\bv,\lambda)$ if and only if there exist dual variables $\muv$ and $\nuv\coloneqq(\nuv_1,\nuv_2,\nuv_3)$ such that
\begin{align}
	\nuv_1\geq \mathbf{0},\;\; \nuv_2\geq\mathbf{0},\;\; \nuv_3\geq \mathbf{0},\label{temp:KKT_bar_beg}\\
	\nuv_1\odot([\sbar]_{\Ebar}\odot\zv)=\mathbf{0},\label{temp:KKT_comple1}\\
	\nuv_2\odot\parens{{\Cm_{\neg\Ebar}}^{\Tr}(\Dm\Cm_{\Ebar}\zv-\bv)- \lambda\mathbf{1}}=\mathbf{0},\label{temp:KKT_comple2}\\
	\nuv_3\odot\parens{{\Cm_{\neg\Ebar}}^{\Tr}(\Dm\Cm_{\Ebar}\zv-\bv)+ \lambda\mathbf{1}}=\mathbf{0},\label{temp:KKT_comple3}\\
	\nabla L(\zv,\muv,\nuv)=\mathbf{0},\label{temp:KKT_bar_end}
\end{align}
where $\zv\coloneqq [\wEQ(\sv;\bv,\lambda)]_{\Ebar}$, $L(\zv,\mu,\nuv)$ is a Lagrangian similar to (\ref{temp:Lagrangian}) in Appx. \ref{proof:partial_order} (one only needs to replace $(\Estar,\sstar)$ therein by $(\Ebar,\sbar)$ to obtain the expression used here). We note that the primal feasibility conditions do not appear in the KKT conditions above because $\wEQ(\sv;\bv,\lambda)$ naturally satisfies them.

From (\ref{temp:csign_wEQ}) and Def. \ref{def:Es_mapping}, one can verify that for every $i\in\neg\Ebar$ and $(\bv,\lambda)\in\interior(\Iselect(\sv))$, we have
\begin{align*}
	\abs{\cv_i^\Tr(\bv-\Dm\Cm\wEQ(\sv;\bv,\lambda))}<\lambda,
\end{align*}
combining which with (\ref{temp:KKT_comple2}) and (\ref{temp:KKT_comple3}) yields $\nuv_2=\nuv_3=\mathbf{0}$. Moreover, from (\ref{temp:sign_wEQ}), we can derive that $[\sv]_{\Ecal}\odot [\zv]_{\Ecal}>\mathbf{0}$, combining which with (\ref{temp:KKT_comple1}) yields $[\nuv_1]_{\Ecal}=\mathbf{0}$ (here we assume that the index systems of components in $\zv$ and $\nuv_1$ are compatible with that of $\wEQ(\sv;\bv,\lambda)$). Applying the results above to the KKT system (\ref{temp:KKT_bar_beg})-(\ref{temp:KKT_bar_end}), we can deduce that for every $(\bv,\lambda)\in\interior(\Iselect(\sv))$, $\wEQ(\sv;\bv,\lambda)=\wstar(\bv,\lambda)$ if and only if there exist dual variables $\muv$ and $\nuv_1$ such that
\begin{align}
	\nuv_1 &\geq \mathbf{0},\label{temp:KKT_Is_1}\\ 
	\begin{bmatrix}
		[\wEQ(\sv;\bv,\lambda)]_{\Ecal}\\
		-[\sbar]_{\Ical}\odot[\nuv_1]_{\Ical}
	\end{bmatrix}
	&={\Cm_{\Ebar}}^{\Tr}\Dm^{\Tr}{\Cm_{\Ebar}\muv},\label{temp:KKT_Is_2}
\end{align}
where $\Ical\coloneqq\Ebar\setminus\Ecal$.

For every $(\bv',\lambda'),(\bv'',\lambda'')\in\breve{I}(\sv)$, from (\ref{temp:breve_Is_wEQ_in_Se}), there respectively exists $(\muv',\nuv'_1)$ (resp. $(\muv'',\nuv''_1)$) such that (\ref{temp:KKT_Is_1})-(\ref{temp:KKT_Is_2}) holds for $(\bv',\lambda')$ (resp. $(\bv'',\lambda'')$). For arbitrary $\theta_1,\theta_2\geq 0$, define
\begin{align*}
	(\bv_{\theta},\lambda_{\theta}) &\coloneqq \theta_1(\bv',\lambda')+\theta_2(\bv'',\lambda''),\\
	\muv_{\theta} &\coloneqq \theta_1\muv'+\theta_2\muv'',\\
	\nuv^{\theta}_1 &\coloneqq \theta_1\nuv'_1+\theta_2\nuv''_1,
\end{align*}
then from linearity of $\wEQ(\sv;\cdot,\cdot)$, one can verify that $(\muv_{\theta},\nuv^{\theta}_1)$ satisfies (\ref{temp:KKT_Is_1})-(\ref{temp:KKT_Is_2}) for $(\bv_{\theta},\lambda_{\theta})$, thus \[\sign(\wstar(\bv_{\theta},\lambda_{\theta}))=\sign(\wEQ(\sv;\bv_{\theta},\lambda_{\theta}))=\sv,\]
which implies that $(\bv_{\theta},\lambda_{\theta})\in\breve{I}(\sv)$, i.e., $\breve{I}(\sv)$ is a convex cone.

\vspace{0.5em}
\noindent
(III) Finally, we show that there do not exist $\sv'\neq\sv$ satisfying $\interior(\breve{I}(\sv'))\neq\emptyset$, such that $\interior(\breve{I}(\sv))$ and $\interior(\breve{I}(\sv'))$ are contained in the interior of the same linear zone in $\wstar(\cdot,\cdot)$.

Suppose that the contrary of the aforementioned statement holds, then from the result (I), we have
\begin{align*}
	\interior(\breve{I}(\sv))\subset\interior(\Iselect(\sv)), && \interior(\breve{I}(\sv'))\subset\interior(\Iselect(\sv')).
\end{align*}
Since $\interior(\breve{I}(\sv))$ and $\interior(\breve{I}(\sv'))$ are contained in the interior of the same linear zone of $\wstar(\cdot,\cdot)$, the slope of $\wEQ(\sv;\cdot,\cdot)$ and $\wEQ(\sv';\cdot,\cdot)$ should be the same,
which further implies $\Iselect(\sv)=\Iselect(\sv')$ from Lemma \ref{lemma:computation_wEQ_IEN} in Appx. \ref{app:computation}. Thus from Lemma \ref{lemma:IEN_property} in Sec. \ref{sec:preparation_for_corres}, we can deduce that for every $(\bv,\lambda)\in\interior(\Iselect(\sv))=\interior(\Iselect(\sv'))$,
\begin{align*}
	\sv\equiv\sign(\wEQ(\sv;\bv,\lambda))\equiv\sign(\wEQ(\sv';\bv,\lambda))\equiv\sv',
\end{align*}
which leads to a contradiction. We note that the first equality above follows from (\ref{temp:sign_wEQ}); by replacing $\sv$ in (\ref{temp:sign_wEQ}) by $\sv'$, we can derive the third equality above. Accordingly, $\bar{I}(\sv)$ is a single linear zone in $\wstar(\cdot,\cdot)$.

\subsection{Step 3: Establish Equality between Signs of Candidate Solution Map and Equicorrelation Signs within a Linear Zone}

Finally, we prove Thm. \ref{thm:correspondence} \ref{it:corres:GP} under the assumption that $\Am$ in (\ref{eq:def_of_G}) has columns in general position. 

We will show that for every $\sv$ satisfying $\interior(\breve{I}(\sv))\neq\emptyset$ and every $(\bv,\lambda)\in\interior(\Iselect(\sv))$, the following equality holds
\begin{align}\label{temp:s_equiv_sign_csign}
	\sv\equiv\sign(\wEQ(\sv;\bv,\lambda))\equiv \se(\bv,\lambda).
\end{align}
The proof of (\ref{temp:s_equiv_sign_csign}) will be provided later. From Lemma \ref{lemma:IEN_property} \ref{it:IEN:wEQ_min_in_SEQNQ} in Sec. \ref{sec:preparation_for_corres}, for every $(\bv,\lambda)\in\interior(\Iselect(\sv))$, $\wEQ(\sv;\bv,\lambda)$ is the unique min-norm element in $\Spoly(\sv;\bv,\lambda)$. From (\ref{temp:s_equiv_sign_csign}), for every $(\bv,\lambda)\in\interior(\Iselect(\sv))$, we further have
\begin{equation}\label{temp:Spoly_equiv_Se}
	\Spoly(\sv;\bv,\lambda)=\Spoly(\se(\bv,\lambda);\bv,\lambda)=\Se(\bv,\lambda),
\end{equation}
where the latter equality follows from Prop. \ref{prop:partial_order} \ref{it:order:equality}. Hence $\wEQ(\sv;\bv,\lambda)$ is the min-norm element in $\Se(\bv,\lambda)$ for every $(\bv,\lambda)\in\interior(\Iselect(\sv))$, which is exactly $\wstar(\bv,\lambda)$. This further implies $\interior(\Iselect(\sv))\subset\breve{I}(\sv)$ from (\ref{eq:domain}) and (\ref{temp:s_equiv_sign_csign}), combining which with $\interior(\Iselect(\sv))\supset\interior(\breve{I}(\sv))$ in Step 2 yields $\Iselect(\sv)=\bar{I}(\sv)$, i.e., $\Iselect(\sv)$ equals to the linear zone $\bar{I}(\sv)$ corresponding to $\sv$. Combining this with (\ref{temp:breve_Is_wEQ_in_Se}), (\ref{temp:Spoly_equiv_Se}), and recall that $\wstar(\cdot,\cdot)$ and $\Se(\cdot,\cdot)$ are continuous (Thm. \ref{thm:continuity}), we can prove Thm. \ref{thm:correspondence} \ref{it:corres:GP}.

\vspace{0.5em}
Now we prove the claimed conditional equality (\ref{temp:s_equiv_sign_csign}). The left equality in (\ref{temp:s_equiv_sign_csign}) follows directly from (\ref{temp:sign_wEQ}). Hence we only need to prove the right equality. From (\ref{temp:csign_wEQ}), we have
\begin{align*}
	(\forall(\bv,\lambda)\in\interior(\Iselect(\sv))) \quad \se(\bv,\lambda)\equiv \sv'\succeq\sv.
\end{align*}
Thus to prove $\sv'=\sv$, by Def. \ref{def:Es_mapping}, it is sufficient to show that there does not exist $i_0\not\in\Ecal\coloneqq\supp(\sv)$ such that 
\begin{align*}
	(\forall(\bv,\lambda)\in\interior(\Iselect(\sv))) \quad \abs{\cv_{i_0}^\Tr(\bv-\Dm\Cm\wEQ(\sv;\bv,\lambda))}=\lambda.
\end{align*}
Suppose that the contrary of the statement above holds. Let us define the following function:
\begin{align*}
	g_{i_0}(\bv,\lambda)\coloneqq \cv_{i_0}^\Tr(\bv-\Dm\Cm\wEQ(\sv;\bv,\lambda)).
\end{align*}
Then since $g_{i_0}(\cdot,\cdot)$ is affine and $\interior(\Iselect(\sv))$ has nonempty interior, we can deduce that
\begin{align*}
	(\forall(\bv,\lambda)\in\Real^{2m}\times\Real)\quad g_{i_0}(\bv,\lambda)\equiv +\lambda\text{ or }{-\lambda}.
\end{align*}
Without loss of generality, let us assume $g_{i_0}(\bv,\lambda)\equiv +\lambda$. Then by expanding the expression of $\wEQ(\sv;\bv,\lambda)$ in $g_{i_0}(\cdot,\cdot)$, we can derive that
\begin{align}
	\cv_{i_0}^\Tr\brackets{\bv-\Dm\Cm_{\Ecal}\parens{{\Cm_{\Ecal}}^{\Tr}\Dm\Cm_{\Ecal}}^{\dagger}{\Cm_{\Ecal}}^\Tr\bv} \equiv 0 \label{temp:EN_b_constant_1}\\
	\cv_{i_0}^\Tr{\Dm\Cm_{\Ecal}\parens{{\Cm_{\Ecal}}^{\Tr}\Dm\Cm_{\Ecal}}^{\dagger}[\sv]_{\Ecal}}\lambda \equiv \lambda \label{temp:EN_b_constant_2}
\end{align}
holds for every $(\bv,\lambda)\in\Real^{2m}\times\Real$. From Lemma \ref{lemma:blk_matrices} \ref{it:blkmat:proj} in Appx. \ref{app:technical_lemma}, (\ref{temp:EN_b_constant_1}) implies $\cv_{i_0}^\Tr\parens{\Imat- P_{\Col(\Cm_{\Ecal})}}\bv=0$ for all $\bv\in\Real^{2m}$, where $P_{\Col(\Cm_{\Ecal})}$ is the orthogonal projector onto $\Col(\Cm_{\Ecal})$. Arbitrariness of $\bv$ implies $\cv_{i_0}\in\Col(\Cm_{\Ecal})$. 

On the other hand, since $\Iselect(\sv)\neq\emptyset$, we have $[\sv]_{\Ecal}={\Cm_{\Ecal}}^\Tr{\Cm_{\Ecal}}^{\Tr\dagger}[\sv]_{\Ecal}$ from Lemma \ref{lemma:computation_wEQ_IEN} \ref{it:comp:IEN} in Appx. \ref{app:computation}, substituting which into (\ref{temp:EN_b_constant_2}) and use Lemma \ref{lemma:blk_matrices} \ref{it:blkmat:proj} in Appx. \ref{app:technical_lemma}, we have $\cv_{i_0}^\Tr{\Cm_{\Ecal}}^{\Tr\dagger}[\sv]_{\Ecal}=1$. Hence (\ref{temp:EN_b_constant_1}) and (\ref{temp:EN_b_constant_2}) are equivalent to 
\begin{align}
	\cv_{i_0}\in\Col\parens{\Cm_{\Ecal}} \text{ and }	\cv_{i_0}^\Tr{\Cm_{\Ecal}}^{\Tr\dagger}[\sv]_{\Ecal}=1.\label{temp:ci0_condition}
\end{align}
Let $\betav\coloneqq {\Cm_{\Ecal}}^{\dagger}\cv_{i_0}$, then the right equality in (\ref{temp:ci0_condition}) yields
\begin{align}
	\sum_{i\in\Ecal} [\betav]_i[\sv]_i=1,\label{temp:coef_sum_to_one}
\end{align}
where we assume that the index system of components in $\betav$ is compatible with that of $\sv$. Moreover, notice that \[\Cm_{\Ecal}\betav=\Cm_{\Ecal}{\Cm_{\Ecal}}^{\dagger}\cv_{i_0}=\cv_{i_0},\]
where the last equality follows from the left inclusion in (\ref{temp:ci0_condition}). Thus we have
\begin{align}
	\sum_{i\in\Ecal} ([\betav]_i[\sv]_i)([\sv]_i\cv_i)=\sum_{i\in\Ecal} [\betav]_i\cv_i=\cv_{i_0}.\label{temp:coef_aff_combination}
\end{align}
Combining (\ref{temp:coef_sum_to_one}) and (\ref{temp:coef_aff_combination}) yields
\begin{align}
	\cv_{i_0}\in\aff\braces{[\sv]_i\cv_i}_{i\in\Ecal}.\label{temp:ci0_in_aff_hull}
\end{align}
Recall that $\Cm=\blkdiag(\Am,\sqrt{\rho}\Am)$, if $i_0\in\braces{1,2,\dots,n}$, then
\begin{align*}
	\av_{i_0}\in\aff\braces{[\sv]_i\av_i}_{i\in\parens{\Ecal\cap\braces{1,2,\dots,n}}},
\end{align*}
which contradicts with the general-positioning condition of $\Am$. If $i_0\in\braces{n+1,n+2,\dots,2n}$, then (\ref{temp:ci0_in_aff_hull}) implies
\begin{align*}
	\sqrt{\rho}\av_{i_0-n}\in\aff\braces{\sqrt{\rho}[\sv]_{i-n}\av_{i-n}}_{i\in\parens{\Ecal\cap\braces{n+1,n+2,\dots,2n}}},
\end{align*}
Notice that in this case, we have $\rho\neq 0$ from $g_{i_0}(\bv,\lambda)\equiv +\lambda$, thus again we yield a contradiction with the general-positioning condition of $\Am$. Hence the equality (\ref{temp:s_equiv_sign_csign}) holds.

\section{Conclusion}\label{sec:conclusion}
In this paper, we studied solution-set geometry and min-norm solution map of the scaled generalized minimax concave (sGMC) model, which is an extension of the LASSO model with a nonconvex sparse regularizer. We have mainly proven the following favorable properties of the sGMC model:
\begin{enumerate}[label=(\alph*)]
	\item the sGMC solution set can be expressed as a nonlinearly parameterized LASSO solution set, and is almost surely a singleton with its unique solution being sparse,
	
	\item the sGMC solution set is a continuous set-valued mapping of $(\yv,\lambda)$, which implies vector-valued continuity of the min-norm solution map $\xstar(\yv,\lambda)$ and other quantities.
	
	\item $\xstar(\yv,\lambda)$ is piecewise linear with each linear zone being a convex cone, and $\sign(\xstar(\yv,\lambda))$ is constant within the interior of each linear zone, 
	
	\item the whole min-norm solution map of $\xstar(\yv,\lambda)$ can be computed by an extension of the LARS algorithm termed E-LARS within finite iterations.
\end{enumerate}
Among these results, the property (a) is a moderate extension of previous work from LASSO to sGMC, the property (b) is a novel result even for the LASSO model, and the properties (c) and (d) are significant extensions of previous work from the min-norm regularization path $\xstar(\lambda)$ of LASSO to the min-norm solution map $\xstar(\yv,\lambda)$ of sGMC.

\section*{Acknowledgment}
The authors sincerely appreciate two anonymous reviewers for insightful comments, which have greatly helped us improve the quality of this paper. They are also grateful to Professor R. J. Tibshirani for kindly answering their questions regarding the LASSO regularization path and LARS algorithm.

\appendices
\section{Known Facts}
\label{app:other_useful_facts}

\subsection{Real and Matrix Analysis}
\begin{fact}[{\cite[Prop. 6.4.12(e)]{tao2016}}]\label{fact:limsup_liminf}
	Let $(a_k)_{k\in\Natural}$ be a sequence of real numbers, if its limit superior $L^+$ (resp. limit inferior $L^-$) is finite, then it is a cluster point of $(a_k)_{k\in\Natural}$, i.e., there exists a subsequence of $(a_k)_{k\in\Natural}$ convergent to $L^+$ (resp. $L^-$).
\end{fact}

\begin{fact}[{\cite[Thm 2.1]{lu2002}}]\label{fact:invertibility}
	Consider the following square matrix 
	\begin{equation*}
		\Rm\coloneqq \begin{bmatrix}
			\Am & \Bm \\
			\Cm & \Dm
		\end{bmatrix},
	\end{equation*}
	where $\Am$ and $\Dm$ are square matrices, $\Bm$ and $\Cm$ are of compatible dimension. Suppose that $\Am$ is invertible. Then $\Rm$ is invertible if and only if $(\Dm-\Cm\Am^{-1}\Bm)$ is invertible.
\end{fact}

\subsection{Convex Analysis}
\begin{fact}[{\cite[Thm. 25.1]{rockafellar1997}}]\label{fact:differentiability_of_convex_functions}
Let $f:\Real^n\to\Real$ be a convex function, and let $\xv\in\Real^n$ be a point where $f$ is finite. If $f$ has a unique subgradient at $\xv$, then $f$ is differentiable at $\xv$.
\end{fact}

\begin{fact}[{\cite[Prop. 16.59]{bauschke2017}}]\label{fact:subgradient_of_inf_function}
Let $F\in\Gamma_0(\Real^m\times\Real^n)$, and set:
\begin{equation*}
f:\Real^m\to [-\infty,+\infty]:\uv\mapsto \inf_{\zv\in\Real^n}F(\uv,\zv).
\end{equation*}
Suppose that $f$ is proper and that $(\bar{\uv},\bar{\zv})\in\Real^m\times\Real^n$ satisfies $f(\bar{\uv})=F(\bar{\uv},\bar{\zv})$, and let $\vv\in\Real^m$. Then $\vv\in\partial f(\bar{\uv})\iff (\vv,\mathbf{0}_n)\in\partial F(\bar{\uv},\bar{\zv})$.
\end{fact}

\begin{fact}[{\cite[Prop. 9.3]{bauschke2017}}]\label{fact:sup_of_convex_functions}
Let $(f_i)_{i\in I}$ be a family of convex functions from $\Real^n$ to $[-\infty,+\infty]$. Then $\sup_{i\in I}f_i$ is convex.
\end{fact}

\subsection{Set-Valued Variational Analysis}

\begin{fact}[{\cite[Example 5.57]{rockafellar2009}}]\label{fact:continuity_of_projection}
	Let $\Mcal:\Real^n\rightrightarrows\Real^m$ be continuous relative to $\dom\Mcal\coloneqq \braces{\thetav\in\Real^n\;\middle\vert\; \Mcal(\thetav)\neq\emptyset}$ and convex-valued. For each $\uv\in\Real^m$, define $\sv_{\uv}:\Real^n\to\Real^m$ by taking $\sv_{\uv}(\thetav)$ to be the projection $P_{\Mcal(\thetav)}(\uv)$ of $\uv$ on $\Mcal(\thetav)$. Then $\sv_{\uv}$ is continuous relative to $\dom\Mcal$.
\end{fact}

\begin{fact}[{\cite[Thm. 2.2]{mangasarian1987}}]
\label{fact:Lipschitzian_continuity_of_linear_program}
Given $\Am\in\Real^{m\times n}$, $\Qm\in\Real^{k\times n}$, for $\yv\in\Real^m$ and $\zv\in\Real^k$, define $\Lcal_{\Am,\Qm}:\Real^m\times\Real^k\rightrightarrows \Real^n$ as:
\begin{align*}
\Lcal_{\Am,\Qm}(\yv,\zv)\coloneqq \braces{\xv\in\Real^n\;\middle\vert\; \Am\xv=\yv,\Qm\xv\leq \zv}.
\end{align*}
Then $\Lcal_{\Am,\Qm}(\cdot,\cdot)$ is ``Lipschitz continuous" on its domain 
\[\dom\Lcal_{\Am,\Qm}\coloneqq \braces{(\yv,\zv)\in\Real^m\times\Real^k \;\middle\vert\; \Lcal_{\Am,\Qm}(\yv,\zv)\neq\emptyset}.\]
More precisely, there exists $\mu_{\Lcal}>0$ and a well-defined norm $\norm{\cdot}_{\Lcal}$ on $\Real^{m+k}$ satisfying: for every $(\yv_1,\zv_1),(\yv_2,\zv_2)\in\dom\Lcal_{\Am,\Qm}$, for every $\xv_1\in\Lcal_{\Am,\Qm}(\yv_1,\zv_1)$, there exists an $\xv_2\in\Lcal(\yv_2,\zv_2)$ closest to $\xv_1$ in the $\ell_\infty$-norm such that
\begin{align*}
\norm{\xv_1-\xv_2}_{\infty}\leq \mu_{\Lcal}\norm{\begin{bmatrix}
\yv_1-\yv_2 \\ \zv_1-\zv_2
\end{bmatrix}}_{\Lcal}.
\end{align*}
\end{fact}

\subsection{Solution Set of General $\ell_1$-Regularization Models}\label{sec:general_l1_model}

Consider the following $\ell_1$-regularization model which embraces LASSO as a special instance:
\begin{equation}\label{eq:general_l1_model}
	\underset{\xv\in\Real^n}{\text{minimize}}\;\;J_{\textrm{GL1}}(\xv)\coloneqq f(\Am\xv)+\lambda \norm{\xv}_1,
\end{equation}
where the loss function $f:\Real^m\to\Real$ is differentiable, strictly convex and bounded from below, $\Am\in\Real^{m\times n}$ is the sensing matrix and $\lambda>0$ is the regularization parameter. Then the solution set of (\ref{eq:general_l1_model}) has the following geometric properties.

\begin{fact}[A slightly improved version of {\cite[Para. 2 in Sec. 2.3]{tibshirani2013}}]\label{fact:same_linear_fit_and_l1_norm}
	The solution set of (\ref{eq:general_l1_model}) is nonempty, closed, bounded and convex. Moreover, every solution of (\ref{eq:general_l1_model}), say $\hat{\xv}$, gives the same linear fit $\Am\hat{\xv}$ and has the same $\ell_1$-norm $\norm{\hat{\xv}}_1$, i.e., if $\hat{\xv}_1,\hat{\xv}_2$ are two solutions of (\ref{eq:general_l1_model}), then \[\Am\hat{\xv}_1=\Am\hat{\xv}_2,\;\;\norm{\hat{\xv}_1}_1=\norm{
		\hat{\xv}_2}_1.\]
\end{fact}
\begin{IEEEproof}
	Since $\norm{\cdot}_1$ is coercive and $f(\Am\cdot)$ is convex and bounded from below, $J_{\textrm{GL1}}$ is a continuous coercive convex function. Hence the solution set of (\ref{eq:general_l1_model}) is nonempty and bounded from \cite[Prop. 11.12 and 11.15]{bauschke2017}. The closedness and convexity of the solution set follow respectively from the continuity of $J_{\textrm{GL1}}$ and \cite[Prop. 11.6]{bauschke2017}. The uniqueness of $\Am\hat{\xv}$ and $\norm{\hat{\xv}}_1$ follows from \cite[Para. 2 in Sec. 2.3]{tibshirani2013}.
\end{IEEEproof}

Moreover, if the sensing matrix $\Am$ has columns in {general position} (cf. Def. \ref{def:general_position} in Sec. \ref{sec:preliminaries}), then the solution uniqueness and sparseness of (\ref{eq:general_l1_model}) can be guaranteed.

\begin{fact}[{\cite[Lemma 5 and the paragraph before it]{tibshirani2013}}]\label{fact:general_l1_unique}
	If $\Am\in\Real^{m\times n}$ has columns in general position, then (\ref{eq:general_l1_model}) has a unique solution with at most $\min\braces{m,n}$ nonzero components.
\end{fact}

\section{A Technical Lemma}\label{app:technical_lemma}

\begin{lemma}\label{lemma:blk_matrices}
	For integers $m>0$, $r_1\geq 0$ and $r_2\geq 0$, let $\Pm$ be a block diagonal matrix such that $\Pm\coloneqq\blkdiag(\Pm_{1},\Pm_{2})$ with $\Pm_{1}\in\Real^{m\times r_1}$, $\Pm_{2}\in\Real^{m\times r_2}$. Let
	\begin{align*}
		\Lmdm_1\coloneqq\begin{bmatrix}
			a\Imat_m & \\
			& d\Imat_m
		\end{bmatrix},\;\;
		\Lmdm_2\coloneqq\begin{bmatrix}
			& b \Imat_m\\
			c\Imat_m & 
		\end{bmatrix}
	\end{align*}
	and $\Lmdm\coloneqq \Lmdm_1+\Lmdm_2$ with $a>0$, $d>0$ and $bc\leq 0$. 
	
	Then the following holds:
	\begin{enumerate}[label=(\alph*)]
		\item \label{it:blkmat:col} $\Col\parens{\Pm^{\Tr}\Lmdm\Pm}=\Col\parens{\Pm^{\Tr}}$.
		
		\item \label{it:blkmat:proj} $\Lmdm\Pm\parens{\Pm^\Tr\Lmdm\Pm^\Tr}^{\dagger}\Pm^\Tr=\Pm\Pm^{\dagger}$.
	\end{enumerate}
\end{lemma}
\begin{IEEEproof}
	\ref{it:blkmat:col} To prove the result \ref{it:blkmat:col}, we only need to prove 
	\[\Null\parens{\Pm^{\Tr}\Lmdm^{\Tr}\Pm}=\Null\parens{\Pm}.\]
	If $\Pm\xv=\mathbf{0}$, then we immediately have $\Pm^{\Tr}\Lmdm^{\Tr}\Pm\xv=\mathbf{0}$, thus $\Null(\Pm)\subset\Null\parens{\Pm^{\Tr}\Lmdm^{\Tr}\Pm}$. Reversely, if $\Pm^{\Tr}\Lmdm^{\Tr}\Pm\xv=\mathbf{0}$, then we can derive the following equality:
	\begin{align*}
		&\Pm^{\Tr}\parens{\Lmdm_1^{\Tr}+\Lmdm_2^{\Tr}}\Pm\xv=\mathbf{0} \\
		\implies & \Pm^{\Tr}\Lmdm^{\Tr}_1\Pm\xv=-\Pm^{\Tr}\Lmdm_2^{\Tr}\Pm\xv \\
		\implies & \Pm^{\Tr}\Pm\xv=-\Pm^{\Tr}\Lmdm_1^{-1}\Lmdm_2^{\Tr}\Pm\xv,
	\end{align*}
	multiply both sides of the equation above by $\Pm^{\Tr\dagger}$ yields
	\begin{align*}
		&\Pm^{\Tr\dagger}\Pm^{\Tr}\Pm\xv=-\Pm^{\Tr\dagger}\Pm^{\Tr}\Lmdm_1^{-1}\Lmdm_2^{\Tr}\Pm\xv \\
		\implies & \Pm\Pm^{\dagger}\Pm\xv=-\Pm\Pm^{\dagger}\Lmdm_1^{-1}\Lmdm_2^{\Tr}\Pm\xv \\
		\implies & \Pm\xv=-\Pm\Pm^{\dagger}\Lmdm^{-1}_1\Lmdm_2^{\Tr}\Pm\xv,
	\end{align*}
	where the second implication follows from the symmetricity of $\Pm\Pm^{\dagger}$ and the third follows from the definition of $\Pm^{\dagger}$. The equation above further implies
	\begin{equation}\label{eq:blk_matrices:temp}
		\parens{\Imat_{2m}+\Pm\Pm^{\dagger}\Lmdm_1^{-1}\Lmdm_2^{\Tr}}\Pm\xv=\mathbf{0}.
	\end{equation}
	If $\parens{\Imat_{2m}+\Pm\Pm^{\dagger}\Lmdm_1^{-1}\Lmdm_2^{\Tr}}$ is invertible, then we have $\Pm\xv=\mathbf{0}$, which implies $\Null\parens{\Pm^{\Tr}\Lmdm^{\Tr}\Pm}\subset\Null\parens{\Pm}$ and completes the proof. Next we prove that $\parens{\Imat_{2m}+\Pm\Pm^{\dagger}\Lmdm_1^{-1}\Lmdm_2^{\Tr}}$ is invertible.
	
	One can verify that $\parens{\Imat_{2m}+\Pm\Pm^{\dagger}\Lmdm_1^{-1}\Lmdm_2^{\Tr}}$ has the following block structure
	\begin{equation*}
		\parens{\Imat_{2m}+\Pm\Pm^{\dagger}\Lmdm_1^{-1}\Lmdm_2^{\Tr}}=\begin{bmatrix}
			\Imat_m & \frac{c}{a}\Pm_{1}\Pm_{1}^{\dagger}\\
			\frac{b}{d}\Pm_{2}\Pm_{2}^{\dagger} & \Imat_{m}
		\end{bmatrix}.
	\end{equation*}
	Since $\Imat_m$ is invertible, by Fact \ref{fact:invertibility}, $\parens{\Imat_{2m}+\Pm\Pm^{\dagger}\Lmdm_1^{-1}\Lmdm_2^{\Tr}}$ is invertible if and only if 
	\[{\Qm}\coloneqq \Imat_{m}+\frac{(-bc)}{ad}\Pm_{2}\Pm_{2}^{\dagger}\Pm_{1}\Pm_{1}^{\dagger}\]
	is invertible. We can prove the invertibility of $\Qm$ by contradiction. Suppose that $\Qm$ is singular, then there exists a nonzero vector $\zv\in\Null({\Qm})$ such that ${\Qm}\zv=\mathbf{0}$, which yields
	\begin{equation}\label{eq:fact_invert_1}
		\zv=\frac{bc}{ad}\Pm_{2}\Pm_{2}^{\dagger}\Pm_{1}\Pm_{1}^{\dagger}\zv\in\Col(\Pm_{2}).
	\end{equation}
	Thus we have $\zv=P_{\Col(\Pm_2)}(\zv)=\Pm_{2}\Pm_{2}^{\dagger}\zv$ from $\zv\in\Col(\Pm_{2})$, substituting which into the RHS of (\ref{eq:fact_invert_1}) yields
	
	\begin{align*}
		&\zv=\frac{bc}{ad}\Pm_{2}\Pm_{2}^{\dagger}\Pm_{1}\Pm_{1}^{\dagger}\Pm_{2}\Pm_{2}^{\dagger}\zv \\ 
		\implies & \left(\Imat_m+\frac{(-bc)}{ad}\Pm_{2}\Pm_{2}^{\dagger}\Pm_{1}\Pm_{1}^{\dagger}\Pm_{2}\Pm_{2}^{\dagger}\right)\zv=\mathbf{0}.
	\end{align*}
	Since $\Pm_1\Pm_1^{\dagger}$ and $\Pm_2\Pm_2^{\dagger}$ are positive semidefinite matrices, the coefficient matrix of the linear equation above is positive definite, hence we have $\zv=\mathbf{0}$, which leads to a contradiction. Thus ${\Qm}$ is invertible, which further implies that $\parens{\Imat_{2m}+\Pm\Pm^{\dagger}\Lmdm_1^{-1}\Lmdm_2^{\Tr}}$ is invertible. 
	
	Combining the discussion above proves the result \ref{it:blkmat:col}.
	
	\vspace{0.5em}
	\noindent
	\ref{it:blkmat:proj} Let us consider the following two linear equations with respect to $\xv\in\Real^{r_1+r_2}$:
	\begin{align}
		\Pm^\Tr\Pm\xv&= \Pm^\Tr\parens{\Lmdm^{\Tr}}^{-1}\yv,\label{temp:blkeq_1}\\
		\Pm^\Tr\Lmdm^\Tr\Pm\xv&=\Pm^\Tr\yv.\label{temp:blkeq_2}
	\end{align}
	The equation (\ref{temp:blkeq_1}) is evidently compatible, and every solution $\xv$ of (\ref{temp:blkeq_1}) satisfies
	\begin{equation}\label{temp:Px_blkeq_1}
		\Pm\xv=\Pm\Pm^{\dagger}\parens{\Lmdm^{\Tr}}^{-1}\yv.
	\end{equation}
	On the other hand, from $\Col\parens{\Pm^\Tr}=\Col\parens{\Pm^\Tr\Lmdm^\Tr\Pm}$ and $\Null\parens{\Pm}=\Null(\Pm^\Tr\Lmdm^\Tr\Pm)$, the equation (\ref{temp:blkeq_2}) is also compatible, and every solution $\xv$ of (\ref{temp:blkeq_2}) satisfies
	\begin{equation}\label{temp:Px_blkeq_2}
		\Pm\xv=\Pm\parens{\Pm^\Tr\Lmdm^\Tr\Pm}^{\dagger}\Pm^\Tr\yv.
	\end{equation}
	
	Let us substitute (\ref{temp:Px_blkeq_1}) into the LHS of (\ref{temp:blkeq_2}), we obtain
	\begin{align*}
		\text{LHS of (\ref{temp:blkeq_2})}=\Pm^\Tr\Lmdm^\Tr\parens{\Pm\Pm^{\dagger}}\parens{\Lmdm^{\Tr}}^{-1}\yv=\Pm^\Tr\yv,
	\end{align*}
	where the latter equality is obtained by expanding the expression of $\Lmdm^\Tr\parens{\Pm\Pm^{\dagger}}\parens{\Lmdm^{\Tr}}^{-1}$ via its block structure. Hence every solution of (\ref{temp:blkeq_1}) satisfies (\ref{temp:blkeq_2}), which implies
	\begin{equation}
		\Pm\Pm^{\dagger}\parens{\Lmdm^{\Tr}}^{-1}\yv=\Pm\parens{\Pm^\Tr\Lmdm^\Tr\Pm}^{\dagger}\Pm^\Tr\yv.
	\end{equation}
	Arbitrariness of $\yv$ yields $\Pm\Pm^{\dagger}=\Pm\parens{\Pm^\Tr\Lmdm^\Tr\Pm}^{\dagger}\Pm^\Tr\Lmdm^\Tr$. Taking transpose of both sides proves the result \ref{it:blkmat:proj}.
\end{IEEEproof}

\section{Derivation of the Optimality Condition (\ref{eq:OPT})}\label{app:derive_OPT}

Recall that $\wv_{\text{e}}\coloneqq\begin{bmatrix}
	\xv_{\text{p}}^\Tr & \zv_{\text{d}}^{\Tr}
\end{bmatrix}^\Tr$ is an extended solution if and only if $({\xv_{\text{p}}},{\zv_{\text{d}}})$ is a saddle point of the objective function $G(\xv,\zv)$ of (\ref{eq:primal_sGMC}), i.e.,
\begin{align}
	{\xv_{\text{p}}}&\in \arg\min_{\xv\in\Real^n}G(\xv,{\zv_{\text{d}}}),\label{eq:saddle_x}\\
	{\zv_{\text{d}}}&\in\arg\max_{\zv\in\Real^n}G({\xv_{\text{p}}},\zv).\label{eq:saddle_z}
\end{align}
Expanding the expression of $G(\xv,{\zv_{\text{d}}})$ in (\ref{eq:saddle_x}), $G({\xv_{\text{p}}},\zv)$ in (\ref{eq:saddle_z}), and omitting constant terms yields
\begin{align}
	{\xv_{\text{p}}}\in\arg\min_{\xv\in\Real^n}& \frac{1}{2}\norm{\yv-\Am\xv}^2_2-\frac{\rho}{2}\norm{\Am\xv-\Am\zv_{\text{d}}}^2_2+\lambda\norm{\xv}_1,\label{eq:simp_saddle_x}\\
	{\zv_{\text{d}}}\in\arg\min_{\zv\in\Real^n}&\frac{\rho}{2}\norm{\Am\xv_{\text{p}}-\Am\zv}^2_2-\sqrt{\rho}\rv^\Tr\Am\zv+\lambda\norm{\zv}_1,\label{eq:simp_saddle_z}
\end{align}
which is further equivalent to
	\begin{equation*}
	\begin{cases}
		\mathbf{0} \in\Am^{\Tr}(\Am((1-\rho){\xv_{\text{p}}}+\rho{\zv_{\text{d}}})-\yv)+\lambda\partial (\norm{\cdot}_1)({\xv_{\text{p}}}),\\
		\mathbf{0}\in\sqrt{\rho}\Am^{\Tr}\parens{\sqrt{\rho}\Am({\zv_{\text{d}}}-{\xv_{\text{p}}})-\rv}+\lambda\partial (\norm{\cdot}_1)({\zv_{\text{d}}}).
	\end{cases}
\end{equation*}
The inclusion problem above can be rewritten as
\begin{equation}\label{eq:KKT_simp}
	\Cm^{\Tr}(\bv-\Dm\Cm\wv_{\text{e}})\in\lambda\partial (\norm{\cdot}_1)(\wv_{\text{e}}),
\end{equation}
Expanding the expression of $\partial (\norm{\cdot}_1)(\cdot)$ in (\ref{eq:KKT_simp}) yields (\ref{eq:OPT}).

\section{Proof of Theorem \ref{thm:shape_uniqueness_sparseness}}
\label{proof:geometry_and_uniqueness}

The proof of Theorem \ref{thm:shape_uniqueness_sparseness} is composed of three steps:
\begin{enumerate}[label=\arabic*)]
	\item \textbf{Step 1:} we show that the primal sGMC model (\ref{eq:primal_sGMC}) can be reformulated as an instance of the general $\ell_1$-regularization model in Appx. \ref{sec:general_l1_model}, whereby we can obtain Thm. \ref{thm:shape_uniqueness_sparseness} \ref{it:shape:linear_fit} \ref{it:shape:unique} for $\Sp$ from Fact \ref{fact:same_linear_fit_and_l1_norm} and \ref{fact:general_l1_unique}.
	
	\item \textbf{Step 2:} by exploiting the LASSO solution mapping as a bridge between $\Sp$ and $\Sd$, we can transfer the properties stated in Thm. \ref{thm:shape_uniqueness_sparseness} \ref{it:shape:linear_fit} \ref{it:shape:unique} from $\Sp$ to $\Sd$.
	
	\item \textbf{Step 3:} we prove Thm. \ref{thm:shape_uniqueness_sparseness} \ref{it:shape:bound} \ref{it:shape:Se} using some intermediate results obtained from the first two steps.
\end{enumerate}

Combining the three steps above completes the proof.

\subsection{Step 1: Reformulate the sGMC Model as a General $\ell_1$-Regularization Model}\label{sec:solution_set_of_sGMC}

We can rewrite the cost function $J_{\text{sGMC}}$ of the primal sGMC model (\ref{eq:primal_sGMC}) as follows:
\begin{align}
	\max_{\zv\in\Real^n}G(\xv,\zv) =\frac{1}{2}\norm{\yv-\Am\xv}^2_2+\lambda\norm{\xv}_1-\min_{\zv\in\Real^n}\parens{\lambda\norm{\zv}_1+\frac{\rho}{2}\norm{\Am(\xv-\zv)}^2_2-\sqrt{\rho}\rv^\Tr\Am\zv}\label{eq:variant_expression_of_sGMC}
\end{align}

Define functions $\bar{g},\bar{f}:\Real^m\to\Real$ as follows:
\begin{align}
	\bar{g}(\uv)&\coloneqq-\min_{\zv\in\Real^n}\parens{\lambda\norm{\zv}_1+\frac{\rho}{2}\norm{\uv-\Am\zv}^2_2-\sqrt{\rho}\rv^\Tr\Am\zv},\label{eq:def_of_bar_g}\\
	\bar{f}(\uv) &\coloneqq  \frac{1}{2}\norm{\yv-\uv}^2_2 -\bar{g}(\uv),\label{eq:def_of_bar_f}
\end{align}
then substituting the definition of $\bar{f}$ into (\ref{eq:variant_expression_of_sGMC}) yields that
\begin{equation}\label{eq:decomposition_of_J_sGMC}
	J_{\textrm{sGMC}}(\xv)=\bar{f}(\Am\xv)+\lambda\norm{\xv}_1.
\end{equation}
Hence $\sGMCpara$ can be rewritten in the form of (\ref{eq:general_l1_model}). In order to apply results introduced in Appx. \ref{sec:general_l1_model} to the sGMC model, we still need to prove that $\bar{f}$ is differentiable, strictly convex and bounded from below. To show this, we first present three important properties of $\bar{g}$.

\begin{lemma}\label{lemma:properties_of_bar_g}
	For every $\lambda>0$, $\rho\in[0,1)$, $\Am\in\Real^{m\times n}$ and $\rv\in\Real^m$, the following holds for $\bar{g}$ defined in (\ref{eq:def_of_bar_g}):
	\begin{enumerate}[label=(\alph*)]
		\item \label{it:properties_of_bar_g:differentiable} $\bar{g}$ is differentiable.
		\item \label{it:properties_of_bar_g:convexity_preserving} $\parens{\frac{\rho}{2}\norm{\yv-\cdot}^2_2-\bar{g}(\cdot)}$ is convex.
		\item \label{it:properties_of_bar_g:bounded} For every $\uv\in\Real^m$, $\bar{g}(\uv)\leq \frac{\rho}{2}\norm{\uv}^2_2$.
	\end{enumerate}
\end{lemma}
\begin{IEEEproof}
	(a) We first prove the differentiability of $\bar{g}$. Define
	\[\bar{G}(\uv,\zv)\coloneqq \lambda\norm{\zv}_1+\frac{\rho}{2}\norm{\uv-\Am\zv}^2_2-\sqrt{\rho}\rv^\Tr\Am\zv,\]
	then one can verify that $\bar{G}\in\Gamma_0(\Real^m\times\Real^n)$ and $\bar{g}(\uv)=\min_{\zv\in\Real^n} \bar{G}(\uv,\zv)$. For $\bar{\uv}\in\Real^m$, let $\bar{\zv}\in\Real^n$ be a global minimizer of $\bar{G}(\bar{\uv},\cdot)$, then from Fact \ref{fact:subgradient_of_inf_function} in Appendix \ref{app:other_useful_facts}, $\vv\in\partial\bar{g}(\bar{\uv})$ if and only if $(\vv,\mathbf{0}_n)\in\partial \bar{G}(\bar{\uv},\bar{\zv})$, i.e.,
	\begin{align*}
		\vv &\in \partial (\bar{G}(\cdot,\bar{\zv}))(\bar{\uv})= \braces{\rho(\bar{\uv}-\Am\bar{\zv})},\\
		\mathbf{0}_n &\in \partial (\bar{G}(\bar{\uv},\cdot))(\bar{\zv}).
	\end{align*}
	The monotone inclusions above are further equivalent to
	\begin{gather*}
		\vv =\rho(\bar{\uv}-\Am\bar{\zv}),\\
		\bar{\zv}\in\arg\min_{\zv\in\Real^n} \lambda\norm{\zv}_1+\frac{1}{2}\norm{(\rv+\sqrt{\rho}\bar{\uv})-\sqrt{\rho}\Am\zv}^2_2,
	\end{gather*}
	which implies that 
	\begin{align*}
		\partial \bar{g}(\bar{\uv})=&\{\rho(\bar{\uv}-\Am\bar{\zv})\,\vert\, \bar{\zv}\in\SLA(\sqrt{\rho}\Am,\rv+\sqrt{\rho}\bar{\uv},\lambda)\},
	\end{align*}
	where $\SLA$ is the LASSO solution set mapping (see the paragraph after (\ref{eq:LASSO})). According to Fact \ref{fact:same_linear_fit_and_l1_norm} in Appx. \ref{sec:general_l1_model}, every solution $\bar{\zv}$ of $\text{LASSO}(\sqrt{\rho}\Am,\rv+\sqrt{\rho}\bar{\uv},\lambda)$ gives the same value of $\Am\bar{\zv}$, hence $\partial \bar{g}(\bar{\uv})$ is single-valued, which implies that $\bar{g}$ is differentiable at $\bar{\uv}$ from Fact \ref{fact:differentiability_of_convex_functions} in Appx. \ref{app:other_useful_facts}. The arbitrariness of $\bar{\uv}$ proves the differentiability of $\bar{g}$ on $\Real^m$. 
	
	\vspace{0.5em}
	\noindent
	(b) Define $\bar{h}(\uv)\coloneqq \frac{\rho}{2}\norm{\yv-\uv}^2_2-\bar{g}(\uv)$. Then we have
	\begin{align*}
		\bar{h}(\uv)=&\max_{\zv\in\Real^n}\frac{\rho}{2}\parens{\norm{\yv-\uv}^2_2-\norm{\uv-\Am\zv}^2_2}-\lambda\norm{\zv}_1+\sqrt{\rho}\rv^\Tr\Am\zv\\
		=&\max_{\zv\in\Real^n}\rho(\Am\zv-\yv)^{\Tr}\uv+\frac{\rho}{2}\parens{\norm{\yv}^2_2-\norm{\Am\zv}^2_2}-\lambda\norm{\zv}_1+\sqrt{\rho}\rv^\Tr\Am\zv.
	\end{align*}
	For given $\zv\in\Real^n$, let us define
	\begin{align*}
		h_{\zv}(\uv)\coloneqq & \rho(\Am\zv-\yv)^{\Tr}\uv+\frac{\rho}{2}\parens{\norm{\yv}^2_2-\norm{\Am\zv}^2_2}-\lambda\norm{\zv}_1+\sqrt{\rho}\rv^\Tr\Am\zv,
	\end{align*}
	then one can verify that $h_{\zv}(\cdot)$ is convex and \[\bar{h}(\uv)=\max_{\zv\in\Real^n}h_{\zv}(\uv).\]
	Hence $\bar{h}$ is convex from Fact \ref{fact:sup_of_convex_functions} in Appendix \ref{app:other_useful_facts}.
	
	\vspace{0.5em}
	\noindent
	(c) Since $\bar{g}(\uv)=\min_{\zv\in\Real^n}\bar{G}(\uv,\zv)$, we have
	\begin{equation*}
		\bar{g}(\uv)\leq \bar{G}(\uv,\mathbf{0}_n)=\frac{\rho}{2}\norm{\uv}^2_2.
	\end{equation*}
	Combining the discussion above completes the proof.
\end{IEEEproof}

Exploiting the properties of $\bar{g}$ described in Lemma \ref{lemma:properties_of_bar_g}, we can prove the differentiability, strict convexity and boundedness from below of $\bar{f}$. See the following proposition.
\begin{proposition}\label{prop:properties_of_bar_f}
	For every $\paragroup\in\Pscr$, $\bar{f}$ defined in (\ref{eq:def_of_bar_f}) is differentiable, strictly convex and bounded from below with the following lower bound:
	\[(\forall \uv\in\Real^m)\;\; \bar{f}(\uv)\geq \frac{-\rho}{2(1-\rho)}\norm{\yv}^2_2.\]
\end{proposition}
\begin{IEEEproof}
	Since $\bar{g}$ is differentiable from Lemma \ref{lemma:properties_of_bar_g} \ref{it:properties_of_bar_g:differentiable}, $\bar{f}(\uv)=\frac{1}{2}\norm{\yv-\uv}^2_2-\bar{g}(\uv)$ is also differentiable. On the other hand, we can rewrite $\bar{f}$ as
	\begin{align*}
		\bar{f}(\uv)=&\frac{1}{2}\norm{\yv-\uv}^2_2-\bar{g}(\uv)\\
		=& \frac{(1-\rho)}{2}\norm{\yv-\uv}^2_2+\parens{\frac{\rho}{2}\norm{\yv-\uv}^2_2-\bar{g}(
			\uv)}.
	\end{align*}
	Since $\rho<1$, $\frac{(1-\rho)}{2}\norm{\yv-\cdot}$ is strictly convex, combining which with Lemma \ref{lemma:properties_of_bar_g} \ref{it:properties_of_bar_g:convexity_preserving} yields the strict convexity of $\bar{f}$. Moreover, from Lemma \ref{lemma:properties_of_bar_g} \ref{it:properties_of_bar_g:bounded}, for every $\uv\in\Real^n$, we have
	\begin{align*}
		\bar{f}(\uv)= &\frac{1}{2}\norm{\yv-\uv}^2_2-\bar{g}(\uv)\\
		\geq & \frac{1}{2}\norm{\yv-\uv}^2_2-\frac{\rho}{2}\norm{\uv}^2_2\\
		= & \frac{(1-\rho)}{2}\norm{\uv-\frac{\yv}{1-\rho}}^2_2-\frac{\rho}{2(1-\rho)}\norm{\yv}^2_2 \\
		\geq & -\frac{\rho}{2(1-\rho)}\norm{\yv}^2_2,
	\end{align*}
	thus $\bar{f}$ is bounded from below.
\end{IEEEproof}

Therefore, Fact \ref{fact:same_linear_fit_and_l1_norm} and \ref{fact:general_l1_unique} introduced in Appx. \ref{sec:general_l1_model} apply to the primal sGMC solution set $\Sp$, thus we have
\begin{itemize}
	\item $\Sp$ is nonempty, closed, bounded and convex,
	
	\item every $\xv\in\Sp$ gives the same linear fit $\Am\xv$ and the same $\ell_1$-norm $\norm{\xv}_1$,
	
	\item if $\Am$ has columns in general position, then $\xv\in\Sp$ is unique with at most $\min\braces{m,n}$ nonzero components,
\end{itemize}
which proves Thm. \ref{thm:shape_uniqueness_sparseness} \ref{it:shape:linear_fit} \ref{it:shape:unique} for $\Sp$.

\subsection{Step 2: Transfer the Properties of $\Sp$ to $\Sd$ via the LASSO Solution Set Mapping}
\label{sec:dual_sGMC_solution_set}

We first show that the primal sGMC solution set $\Sp$ and the dual one $\Sd$ are connected by the LASSO solution set mapping $\SLA(\cdot,\cdot,\cdot)$ (defined in the paragraph after (\ref{eq:LASSO})). This can be derived from the saddle point condition of the min-max form sGMC model (\ref{eq:sGMC_as_minmax}): by properly scaling and rearranging the terms in (\ref{eq:simp_saddle_x}) and (\ref{eq:simp_saddle_z}) in Appx. \ref{app:derive_OPT}, one can veriy that $(\xv_{\text{p}},\zv_{\text{d}})\in\Real^n\times\Real^n$ is a saddle point of (\ref{eq:sGMC_as_minmax}) if and only if
\begin{align}
	\xv_{\text{p}} &\in\SLA\parens{\Am,\frac{\yv-\rho\Am\zv_{\text{d}}}{(1-\rho)},\frac{\lambda}{(1-\rho)}}, \label{eq:xp_in_SLA}\\
	\zv_{\text{d}}& \in\SLA\parens{\sqrt{\rho}\Am,\rv+\sqrt{\rho}\Am\xv_{\text{p}},\lambda}.\label{eq:zd_in_SLA}
\end{align}

Recall that we have proven in the step 1 that every $\xv_{\text{p}}\in\Sp$ gives the same linear fit $\Am\xv_{\text{p}}$, hence $\betap\coloneqq \Am\xv_{\text{p}}$ only depend on the problem parameter $\paragroup$, and is independent of the selection of $\zv_{\text{p}}$. Applying this into (\ref{eq:zd_in_SLA}), and notice that (\ref{eq:zd_in_SLA}) holds for and only for arbitrary $\zv_{\text{d}}\in\Sd$, we have
\begin{align}
	\Sd=\SLA(\sqrt{\rho}\Am,\rv+\sqrt{\rho}\betap,\lambda),\label{temp:Sd_as_LASSO}
\end{align}

Applying Fact \ref{fact:same_linear_fit_and_l1_norm} in Appx. \ref{sec:general_l1_model} to the LASSO solution set in (\ref{temp:Sd_as_LASSO}), one can verify that 
\begin{itemize}
	\item $\Sd$ is nonempty, closed, bounded and convex,
	\item every $\zv\in\Sd$ gives the same linear fit $\Am\zv$ and the same $\ell_1$-norm $\norm{\zv}_1$,
	
	\item if $\Am$ has columns in general position and $\rho\in(0,1)$, then $\sqrt{\rho}\Am$ also has columns in general position, thus $\zv\in\Sd$ is unique with at most $\min\braces{m,n}$ nonzero components.
\end{itemize}
If $\rho=0$, then $\sqrt{\rho}\Am=\mathbf{O}_{m\times n}$ does not have columns in general position regardless of $\Am$. However, in this case, the LASSO problem in the RHS of (\ref{temp:Sd_as_LASSO}) reduces to the minimization of $\norm{\cdot}_1$, hence $\zv=\mathbf{0}_n$ is the unique dual solution in $\Sd$. Accordingly, Thm. \ref{thm:shape_uniqueness_sparseness} \ref{it:shape:linear_fit} \ref{it:shape:unique} holds for $\Sd$.

\subsection{Step 3: Prove Theorem \ref{thm:shape_uniqueness_sparseness} \ref{it:shape:bound} \ref{it:shape:Se} by Intermediate Results}

It has been shown in the step 1 and 2 that for every $\xv_{\text{p}}\in\Sp$, $\zv_{\text{d}}\in\Sd$, the value of $\betap\coloneqq \Am\xv_{\text{p}}$ and $\betad\coloneqq \Am\zv_{\text{d}}$ only depend on the problem parameter $\paragroup$. Applying this to (\ref{eq:xp_in_SLA}) and (\ref{eq:zd_in_SLA}) yields the two equalities (\ref{eq:Se_SpSd}) and (\ref{eq:Se_variant_expression}) in Thm. \ref{thm:shape_uniqueness_sparseness} \ref{it:shape:Se}.

In the sequel, we prove Thm. \ref{thm:shape_uniqueness_sparseness} \ref{it:shape:bound}. We first show that the $\ell_1$-norm bound (\ref{eq:l1_bound}) holds for every $\xv_{\text{p}}\in\Sp$. From (\ref{eq:decomposition_of_J_sGMC}) and Prop. \ref{prop:properties_of_bar_f}, the following holds
\begin{equation}\label{eq:bound_of_l1_norm_of_x}
	J_{\textrm{sGMC}}(\xv_{\text{p}})\geq \lambda\norm{\xv_{\text{p}}}_1-\frac{\rho}{2(1-\rho)}\norm{\yv}^2_2.
\end{equation}
Since $\xv_{\text{p}}\in\Sp$ is a global minimizer of the cost function $J_{\textrm{sGMC}}$ of the primal sGMC model (\ref{eq:primal_sGMC}), we have
\begin{align*}
	J_{\textrm{sGMC}}(\xv_{\text{p}})=\min_{{\xv}\in\Real^n}J_{\textrm{sGMC}}({\xv})\leq J_{\textrm{sGMC}}(\mathbf{0}).
\end{align*}
Expanding the expression of $J_{\textrm{sGMC}}(\mathbf{0})$ by (\ref{eq:variant_expression_of_sGMC}) yields 
\begin{align}
	J_{\textrm{sGMC}}(\xv_{\text{p}})\leq&\frac{1}{2}\norm{\yv}^2_2-\min_{{\zv}\in\Real^n}\left( \lambda\norm{{\zv}}_1+\frac{\rho}{2}\norm{\Am{\zv}}^2_2-\sqrt{\rho}\rv^\Tr\Am\zv \right)\nonumber\\
	\leq &\frac{1}{2}\parens{\norm{\yv}^2_2+\norm{\rv}^2_2},\label{temp:JsGMC_xp_bound}
\end{align}
where the second inequality follows from
\begin{align}
	 &\lambda\norm{{\zv}}_1+\frac{\rho}{2}\norm{\Am{\zv}}^2_2-\sqrt{\rho}\rv^\Tr\Am\zv\nonumber\\
	 =&\lambda\norm{{\zv}}_1+\frac{1}{2}\norm{\sqrt{\rho}\Am{\zv}-\rv}^2_2-\frac{1}{2}\norm{\rv}^2_2 \label{temp:inequality_auxliary_term}\\
	 \geq & -\frac{1}{2}\norm{\rv}^2_2.\nonumber
\end{align}
Combining (\ref{temp:JsGMC_xp_bound}) with (\ref{eq:bound_of_l1_norm_of_x}) yields
\begin{align*}
	\lambda\norm{\xv_{\text{p}}}_1-\frac{\rho}{2(1-\rho)}\norm{\yv}^2_2\leq \frac{1}{2}\norm{\yv}^2_2+\frac{1}{2}\norm{\rv}^2_2,
\end{align*}
which implies \begin{equation}\label{eq:xp_l1_bound}
	(\forall \xv_{\text{p}}\in\Sp)\;\;\norm{\xv_{\text{p}}}_1\leq \frac{\norm{\yv}_2^2}{2\lambda(1-\rho)}+\frac{\norm{\rv}^2_2}{2\lambda}.
\end{equation}

Next we show that the $\ell_1$-norm bound (\ref{eq:l1_bound}) also holds for every $\zv_{\text{d}}\in\Sd$. Let us consider the cost function $\bar{J}_{\text{sGMC}}$ of the dual sGMC model (\ref{eq:dual_sGMC}):
\begin{align*}
	\bar{J}_{\textrm{sGMC}}(\zv)\coloneqq & -\min_{\xv\in\Real^n}G(\xv,\zv).
\end{align*}
Since every $\zv_{\text{d}}\in\Sd$ is a global minimizer of $\bar{J}_{\text{sGMC}}$, we have
\begin{align}
	\bar{J}_{\textrm{sGMC}}(\zv_{\text{d}})&=\min_{\zv\in\Real^n}\bar{J}_{\textrm{sGMC}}(\zv)\leq \bar{J}_{\textrm{sGMC}}(\mathbf{0})\nonumber\\
	&=-\min_{\xv\in\Real^n}\left[ \frac{1}{2}\norm{\yv-\Am\xv}^2_2-\frac{\rho}{2}\norm{\Am\xv}^2_2+\lambda\norm{\xv}_1 \right]\nonumber\\
	&\leq -\min_{\xv\in\Real^n}\left[ \frac{1}{2}\norm{\yv-\Am\xv}^2_2-\frac{\rho}{2}\norm{\Am\xv}^2_2\right]\nonumber\\
	&\leq -\min_{\uv\in\Real^m}\left[ \frac{1}{2}\norm{\yv-\uv}^2_2-\frac{\rho}{2}\norm{\uv}^2_2\right]\label{eq:thm_of_dual2}\\
	&=\frac{\rho}{2(1-\rho)}\norm{\yv}^2_2,\label{eq:thm_of_dual3}
\end{align}
where (\ref{eq:thm_of_dual3}) is obtained by solving the quadratic program in (\ref{eq:thm_of_dual2}) in closed-form. On the other hand, for every $\zv\in\Real^n$, the following holds:
\begin{align*}
	\bar{J}_{\textrm{sGMC}}(\zv) =&\max_{\xv\in\Real^n}-G(\xv,\zv)\geq -G(\mathbf{0},\zv)\\
	=&\lambda\norm{\zv}_1-{\frac{1}{2}\norm{\yv}^2_2+\frac{\rho}{2}\norm{\Am\zv}^2_2}-\sqrt{\rho}\rv^\Tr\Am\zv\\
	\geq&\lambda\norm{\zv}_1-\frac{1}{2}\norm{\yv}^2_2-\frac{1}{2}\norm{\rv}^2_2,
\end{align*}
where the last inequality follows from (\ref{temp:inequality_auxliary_term}).
Combining the inequality above with (\ref{eq:thm_of_dual3}) yields
\begin{align*}
	\lambda\norm{\zv_{\text{d}}}_1-\frac{1}{2}\norm{\yv}^2_2-\frac{1}{2}\norm{\rv}^2_2\leq \frac{\rho}{2(1-\rho)}\norm{\yv}^2_2,
\end{align*}
which implies \begin{equation}\label{eq:zd_l1_bound}
	(\forall \zv_{\text{d}}\in\Sd)\;\;\norm{\zv_{\text{d}}}_1\leq \frac{\norm{\yv}_2^2}{2\lambda(1-\rho)}+\frac{\norm{\rv}^2_2}{2\lambda}.
\end{equation}
Combining (\ref{eq:xp_l1_bound}) and (\ref{eq:zd_l1_bound}) completes the proof of the $\ell_1$-norm bound (\ref{eq:l1_bound}) in Thm. \ref{thm:shape_uniqueness_sparseness} \ref{it:shape:bound}.

\section{Computation of Candidate Solution Maps and Candidate Zones}
\label{app:computation}

In this appendix, we discuss about computation of the candidate solution map $\wEQ(\sv;\bv,\lambda)$ and the candidate zone $\Iselect(\sv)$ for a given $\sv\in\braces{+1,0,-1}^{2n}$. In the sequel, we present the main results. The derivations of them will be provided later in the next two subsections. 

The following Lemma \ref{lemma:computation_wEQ_IEN} indicates that given $\sv$, $\wEQ(\sv;\bv,\lambda)$ can be computed in closed-form, whilst $\Iselect(\sv;\bv,\lambda)$ can be obtained by solving a linear program.

\begin{lemma}\label{lemma:computation_wEQ_IEN}
	Given $\sv\in\braces{+1,0,-1}^{2n}$, let $\Ecal\coloneqq \supp(\sv)$. Then the following holds:
	\begin{enumerate}[label=(\alph*)]
		\item \label{it:comp:wEQ} The candidate solution map $\wEQ(\sv;\bv,\lambda)$ has the following closed-form expression
		\begin{equation}
			\begin{cases}
				[\wEQ(\sv;\bv,\lambda)]_{\Ecal}=\Rm(\sv)\begin{bmatrix}
					\bv\\ \lambda
				\end{bmatrix},\\
				[\wEQ(\sv;\bv,\lambda)]_{\neg\Ecal}=\mathbf{0},
			\end{cases}		
		\end{equation}
		where $\Rm(\sv)$ is defined as
		\begin{equation}\label{eq:Qs}
			\Rm(\sv)\coloneqq \parens{{\Cm_{\Ecal}}^{\Tr}\Dm\Cm_{\Ecal}}^{\dagger}
			\begin{bmatrix}
				\Cm^{\Tr}_{\Ecal} & 
				-[\sv]_{\Ecal}
			\end{bmatrix}.
		\end{equation}

		\item \label{it:comp:IEN} Computation of the candidate zone $\Iselect(\sv)$ can be divided into two cases: if we have		
		\[[\sv]_{\Ecal}\not\in\Col\parens{{\Cm_{\Ecal}}^\Tr},\]
		then $\Iselect(\sv)=\emptyset$; otherwise, $\Iselect(\sv)$ equals to the set of $(\bv,\lambda)\in\Real^{2m}\times\Real$ satisfying the following conditions,
		\begin{subequations}
			\makeatletter
			\def\@currentlabel{EN}
			\makeatother
			\label{eq:EN}
			\renewcommand{\theequation}{EN-\alph{equation}}
			\begin{empheq}[left=\empheqlbrace]{align}
				&(\forall i\in\Ecal) && [\sv]_i[\wEQ(\sv;\bv,\lambda)]_i\geq 0, \label{eq:EN-E}&\\
				&(\forall i\in\neg\Ecal) && \abs{\xi_i(\wEQ(\sv;\bv,\lambda))}\leq \lambda,\label{eq:EN-nE}&\\
				& && \lambda>0, \nonumber &
			\end{empheq}
		\end{subequations}
		where $\xi_i(\wv)\coloneqq \cv_i^\Tr(\bv-\Dm\Cm\wv)$. We note that (\ref{eq:EN}) is further equivalent to the following linear program
		\begin{subequations}\label{eq:LP_zone}
			\begin{align}
				\begin{bmatrix}
					\Tm_1(\sv)\\
					\Tm_2(\sv)\\
					\Tm_3(\sv)
				\end{bmatrix}\begin{bmatrix}
					\bv \\ \lambda
				\end{bmatrix}\leq \mathbf{0},\\
				\lambda>0,
			\end{align}			
		\end{subequations}
		where $\Tm_1(\sv)$, $\Tm_2(\sv)$ and $\Tm_3(\sv)$ are defined as
		\begin{align*}
			\Tm_1(\sv) &\coloneqq {{{\parens{-[\sv]_{\Ecal}\mathbf{1}^\Tr_{2m+1}}}\odot \Rm(\sv)}},\\
			\Tm_2(\sv) &\coloneqq (\Cm_{\neg\Ecal})^\Tr\Dm\Cm_{\Ecal}\Rm(\sv)-\begin{bmatrix}
				(\Cm_{\neg\Ecal})^\Tr & \mathbf{1}_{\abs{\neg\Ecal}}
			\end{bmatrix},\\
			\Tm_3(\sv) &\coloneqq \begin{bmatrix}
				(\Cm_{\neg\Ecal})^\Tr & -\mathbf{1}_{\abs{\neg\Ecal}}
			\end{bmatrix}-(\Cm_{\neg\Ecal})^\Tr\Dm\Cm_{\Ecal}\Rm(\sv).
		\end{align*}
		From the expression of (\ref{eq:LP_zone}), one can verify that $\Iselect(\sv)$ is a convex cone in $\Real^{2m}\times\Real_{++}$.
	\end{enumerate}
\end{lemma}

The next Lemma \ref{lemma:compuation_IEN_line} indicates that if the movement of $(\bv,\lambda)$ is limited to a straight line in $\Real^{2m+1}$, then the intersection of $\Iselect(\sv)$ and this straight line can be computed in closed-form.

\begin{lemma}\label{lemma:compuation_IEN_line}
	Given $\sv\in\braces{+1,0,-1}^{2n}$, suppose that $(\bv,\lambda)\in\Real^{2m}\times\Real$ is only allowed to move in a straight line:
	\begin{align*}
		\bv(t) & = \bv_0+\Deltab t\\
		\lambda(t) & = \lambda_0 +\Deltalmd t
	\end{align*}
	where $t\in\Real$ is a time index. For $k,b\in\Real$, let us define
	\begin{align*}
		f_{\mathrm{tmax}}(k,b)\coloneqq & \sup\braces{t\in\Real\;\middle\vert \; kt\leq b}\\
		=&\begin{cases}
			b/k, &\text{if }k>0,\\
			-\infty, & \text{if }k=0\text{ and }b< 0,\\
			+\infty, & \text{otherwise}.
		\end{cases},
	\end{align*}
	
	Then the following holds:
	\begin{enumerate}[label=(\alph*)]
		\item \label{it:tmax:ta} Let $t_i^a$ be the supremum value of $t\in\Real$ such that $(\bv(t),
		\lambda(t))$ satisfies the $i$th $\sv$-(\ref{eq:EN-E}) condition, then
		\begin{equation*}
			t_i^a=f_{\mathrm{tmax}}([\sv]_i[\pv(\sv)]_i,[\sv]_i[\qv(\sv)]_i),
		\end{equation*}
		where $\pv(\sv),\qv(\sv)\in\Real^{2n}$ are
		\begin{align*}
			[\pv(\sv)]_{\Ecal}&\coloneqq -\Rm(\sv)\begin{bmatrix}
				\Deltab \\ \Deltalmd
			\end{bmatrix},& [\pv(\sv)]_{\neg\Ecal}=\mathbf{0},\\
			[\qv(\sv)]_{\Ecal} &\coloneqq \Rm(\sv)\begin{bmatrix}
				\bv_0 \\ \lambda_0
			\end{bmatrix}, & [\qv(\sv)]_{\neg\Ecal}=\mathbf{0},
		\end{align*}
		with $\Rm(\sv)$ defined in Lemma \ref{lemma:computation_wEQ_IEN} \ref{it:comp:wEQ}. 
		
		\item \label{it:tmax:tb} Let $t_i^b$ be the supremum value of $t\in\Real$ such that $(\yv(t),
		\lambda(t))$ satisfies the $i$th $\sv$-(\ref{eq:EN-nE}) condition, then
		\begin{align*}
			t^b_i=&\min\Big\{f_{\mathrm{tmax}}\parens{-\cv_i^\Tr\uv(\sv)-\Deltalmd,\lambda_0+\cv_i^\Tr\vv(\sv)}, f_{\mathrm{tmax}}\parens{\cv_i^\Tr\uv(\sv)-\Deltalmd,\lambda_0-\cv_i^\Tr\vv(\sv)}\Big\},
		\end{align*}
		where $\uv(\sv)$ and $\vv(\sv)$ are defined as
		\begin{align*}
			\uv(\sv) & \coloneqq  \Deltab+\Dm\Cm_{\Ecal}[\pv(\sv)]_{\Ecal},\\
			\vv(\sv) &\coloneqq \bv_0-\Dm\Cm_{\Ecal}[\qv(\sv)]_{\Ecal}.
		\end{align*}
		
		\item \label{it:tmax:tc} Define 		
		$\Tzone(\sv)\coloneqq \braces{t\in\Real\;\middle\vert\; (\bv(t),\lambda(t))\in \Iselect(\sv)}$. Suppose that $\interior(\Tzone(\sv))\neq\emptyset$, then  \[\sup(\Tzone)=\min\braces{\min_{i\in\Ecal} t^a_i, \min_{i\in\neg\Ecal} t^b_i, t^c},\] 
		where $t^c\coloneqq f_{\mathrm{tmax}}(-\Deltalmd,-\lambda_0)$. Notice that $\inf (\Tzone)$ can be computed similarly by reversing the velocity vector $(\Deltab,\Deltalmd)$ in the definition of $(\bv(t),\lambda(t))$. Moreover, if the computed value of $\inf(\Tzone(\sv))$ is strictly smaller than that of $\sup(\Tzone(\sv))$, then one can confirm that $\interior(\Tzone(\sv))\neq\emptyset$ and the computed values are correct.
	\end{enumerate}
\end{lemma}

\subsection{Proof of Lemma \ref{lemma:computation_wEQ_IEN}}

\noindent
\ref{it:comp:wEQ} One can verify that $\wv\in\Real^{2n}$ satisfies $\sv$-(\ref{eq:EQ}) at $(\bv,\lambda)\in\Real^{2m}\times\Real_{++}$ if and only if
\begin{equation}\label{eq:EQ-simp}
	\begin{cases}
		{\Cm_{\Ecal}}^{\Tr}\Dm\Cm_{\Ecal}[\wv]_{\Ecal}={\Cm_{\Ecal}}^{\Tr}\bv-\lambda[\sv]_{\Ecal},\\
		[\wv]_{\neg\Ecal}=\mathbf{0}.
	\end{cases}
\end{equation}
From Def. \ref{def:EQ_NQ_quantities}, $\wEQ(\sv;\bv,\lambda)$ is the least-squares fitting of (\ref{eq:EQ-simp}), which has the following closed-form expression
\begin{equation*}
	\begin{cases}
		[\wEQ(\sv;\bv,\lambda)]_{\Ecal}=\parens{{\Cm_{\Ecal}}^{\Tr}\Dm\Cm_{\Ecal}}^{\dagger}\parens{\Cm^{\Tr}_{\Ecal}\bv-\lambda[\sv]_{\Ecal}},\\
		[\wEQ(\sv;\bv,\lambda)]_{\neg\Ecal}=\mathbf{0}.
	\end{cases}
\end{equation*}
This expression is equivalent to the one given in the result \ref{it:comp:wEQ}.

\vspace{0.5em}
\noindent
\ref{it:comp:IEN} From Def. \ref{def:EQ_NQ_quantities}, $\Iselect(\sv)$ is the set of $(\bv,\lambda)\in\Real^{2m}\times\Real_{++}$ such that $\wEQ(\sv;\bv,\lambda)$ simultaneously satisfies $\sv$-(\ref{eq:EQ}) and $\sv$-(\ref{eq:NQ}) at $(\bv,\lambda)$.

Since $\wEQ(\sv;\bv,\lambda)$ is the least-squares fitting of (\ref{eq:EQ}) at $(\bv,\lambda)$, it satisfies $\sv$-(\ref{eq:EQ}) at $(\bv,\lambda)\in\Real^{2m}\times\Real_{++}$ if and only if this linear equation has a solution. Let $\Ecal\coloneqq \supp(\sv)$, then from (\ref{eq:EQ-simp}), the solution set of $\sv$-(\ref{eq:EQ}) at $(\bv,\lambda)$ is nonempty if and only if the following equation with respect to $[\wv]_{\Ecal}$:
\begin{align}
	{\Cm_{\Ecal}}^{\Tr}\Dm\Cm_{\Ecal}[\wv]_{\Ecal}={\Cm_{\Ecal}}^{\Tr}\bv-\lambda[\sv]_{\Ecal}\label{eq:EQa-simp}
\end{align}
is compatible, which is further equivalent to
\begin{equation}\label{eq:compatible_condition}
	{\Cm_{\Ecal}}^{\Tr}\bv-\lambda[\sv]_{\Ecal}\in\Col\parens{{\Cm_{\Ecal}}^{\Tr}\Dm\Cm_{\Ecal}}.
\end{equation}
Notice that from Lemma \ref{lemma:blk_matrices} in Appx. \ref{app:technical_lemma}, we have
\begin{equation*}
	\Col\parens{{\Cm_{\Ecal}}^{\Tr}}=\Col\parens{{\Cm_{\Ecal}}^{\Tr}\Dm\Cm_{\Ecal}}.
\end{equation*}
Hence we always have ${\Cm_{\Ecal}}^\Tr\bv\in\Col\parens{{\Cm_{\Ecal}}^{\Tr}\Dm\Cm_{\Ecal}}$. Substituting this into (\ref{eq:compatible_condition}) and recall $\lambda>0$, we can yield that $\sv$-(\ref{eq:EQ}) at $(\bv,\lambda)$ has a solution if and only if
\begin{equation}\label{eq:compat_condition_simp}
	[\sv]_{\Ecal}\in\Col\parens{{\Cm_{\Ecal}}^{\Tr}\Dm\Cm_{\Ecal}}=\Col\parens{{\Cm_{\Ecal}}^{\Tr}}.
\end{equation}
Notice that this condition does not depend on the value of $(\bv,\lambda)$. Accordingly, if (\ref{eq:compat_condition_simp}) does not hold, then $\Iselect(\sv)=\emptyset$.

If the compatibility condition (\ref{eq:compat_condition_simp}) holds, then $(\bv,\lambda)\in\Iselect(\sv)$ if and only if $\wEQ(\sv;\bv,\lambda)$ satisfies $\sv$-(\ref{eq:NQ}) at $(\bv,\lambda)\in\Real^{2m}\times\Real_{++}$, i.e., $(\bv,\lambda)$ satisfies (\ref{eq:EN}) condition. Using the result \ref{it:comp:wEQ}, one can verify that (\ref{eq:EN}) is equivalent to (\ref{eq:LP_zone}).

\subsection{Derivation of Lemma \ref{lemma:compuation_IEN_line}}

One can verify from Lemma \ref{lemma:computation_wEQ_IEN} \ref{it:comp:wEQ} that the following holds:
\begin{align*}
	\wEQ(\sv;\bv(t),\lambda(t))&=\qv(\sv)-\pv(\sv)t,\\
	\bv-\Dm\Cm\wEQ(\sv;\bv(t),\lambda(t))&=\vv(\sv)+\uv(\sv)t.
\end{align*}

For $i\in\Ecal$ and $(\bv,\lambda)\coloneqq (\bv(t),\lambda(t))$, the $i$th $\sv$-(\ref{eq:EN-E}) condition can be rewritten as
\begin{equation*}
	[\sv]_i[\pv(\sv)]_i t\leq  [\sv]_i[\qv(\sv)]_i,
\end{equation*}
which gives the result \ref{it:tmax:ta}; the $i$th $\sv$-(\ref{eq:EN-nE}) condition can be rewritten as
\begin{align*}
	-(\cv_i^\Tr\uv(\sv)+\Deltalmd)t &\leq  \lambda_0+\cv_i^\Tr\vv(\sv),\\
	(\cv_i^\Tr\uv(\sv)-\Deltalmd)t &\leq \lambda_0-\cv_i^\Tr\vv(\sv),
\end{align*}			
which gives the result \ref{it:tmax:tb}.

Since $\interior(\Tzone(\sv))\neq\emptyset$, then from Lemma \ref{lemma:computation_wEQ_IEN} \ref{it:comp:IEN}, $\Tzone(\sv)$ is the set of $t\in\Real$ such that $(\bv(t),\lambda(t))$ satisfies (\ref{eq:EN}) conditions, which gives the result \ref{it:tmax:tc}.

\section{Proofs for Auxiliary Results in Section \ref{sec:proof:correspondence}}
\label{app:auxiliary_results}

This appendix addresses the proofs of three auxiliary results used in Sec. \ref{sec:proof:correspondence}: Lemma \ref{lemma:SEQNQ_property}, \ref{lemma:IEN_property} and Prop. \ref{prop:partial_order}. Due to the dependencies among these results, we will prove them in the following order: Lemma \ref{lemma:SEQNQ_property}, Prop. \ref{prop:partial_order}, Lemma \ref{lemma:IEN_property}.

\subsection{Proof of Lemma \ref{lemma:SEQNQ_property}}\label{proof:SEQNQ_property}
\begin{IEEEproof}
	For $\sv\in\braces{+1,0,-1}^{2n}$ and $(\bv,\lambda)\in\Real^{2m}\times\Real_{++}$, define $\Ecal\coloneqq \supp(\sv)$. Let $\wv\in\Spoly(\sv;\bv,\lambda)$, then ${\wv}$ satisfies $\sv$-(\ref{eq:EQ}) and $\sv$-(\ref{eq:NQ}) simultaneously at $(\bv,\lambda)$, i.e.,
	\begin{empheq}[left=\empheqlbrace]{align}
		&(\forall i\in\Ecal) &&\cv_{i}^{\Tr}(\bv-\Dm\Cm\wv)=\lambda[\sv]_{i}, \tag{\ref{eq:EQ-E}}&\\
		&(\forall i\in\neg\Ecal)&& [\wv]_{i}=0, \tag{\ref{eq:EQ-nE}}&
	\end{empheq}
	\begin{empheq}[left=\empheqlbrace]{align}
		&(\forall i\in\Ecal) && [\sv]_i[\wv]_i\geq 0, \tag{\ref{eq:NQ-E}}&\\
		&(\forall i\in\neg\Ecal) &&  |\cv_i^{\Tr}(\bv-\Dm\Cm\wv)|\leq \lambda. \tag{\ref{eq:NQ-nE}}&
	\end{empheq}
	To prove $\wv\in\Se(\bv,\lambda)$, we only need to show that $\wv$ satisfies the $i$th (\ref{eq:OPT}) condition (Lemma \ref{lemma:OPT}) for every $i\in\braces{1,2,\dots,2n}$, i.e.,
	\begin{equation}
		\cv_i^{\Tr}(\bv-\Dm\Cm\wv)\in\begin{cases}
			\braces{\lambda\,\sign([\wv]_i)}, & \textrm{if }[\wv]_i\neq 0,\\
			[-\lambda,\lambda], & \textrm{if }[\wv]_i=0.
		\end{cases}\tag{\ref{eq:OPT}}
	\end{equation}
	Let us define $\xi_i(\wv)\coloneqq \cv_i^\Tr(\bv-\Dm\Cm\wv)$ for abbreviation.
	
	For every $i\in\neg{\Ecal}$, (\ref{eq:EQ-nE}) and (\ref{eq:NQ-nE}) imply that $[\wv]_i=0$ and $\abs{\xi_i(\wv)}\leq \lambda$, hence $\wv$ satisfies the $i$th (\ref{eq:OPT}) condition. 
	
	On the other hand, for every $i\in{\Ecal}$, since $\Ecal\coloneqq\supp(\sv)$, we have $[\sv]_i\in\braces{\pm 1}$, which further implies from (\ref{eq:EQ-E}) that $\abs{\xi_i(\wv)}=\lambda$. Hence we can derive that
	\begin{itemize}
		\item if $[{\wv}]_i=0$, then since $\abs{\xi_i(\wv)}=\lambda\in[-\lambda,\lambda]$, $\wv$ satisfies the $i$th (\ref{eq:OPT}) condition,
		
		\item if $[{\wv}]_i\neq 0$, then from (\ref{eq:NQ-E}) and $[\sv]_i\in\braces{\pm 1}$ we have $[{\sv}]_i=\sign([{\wv}]_i)$, combining which with (\ref{eq:EQ-E}) implies 
		\[\xi_i(\wv)=\lambda[\sv]_i=\lambda\,\sign([\wv]_i).\]
		Thus $\wv$ satisfies the $i$th (\ref{eq:OPT}) condition.
	\end{itemize}
	
	Accordingly, ${\wv}$ satisfies the $i$th (\ref{eq:OPT}) condition for every $i\in\{1,2,\dots,2n\}$, which yields ${\wv}\in\Se(\bv,{\lambda})$. Hence we conclude that $\Spoly(\sv;\bv,\lambda)\subset\Se(\bv,\lambda)$.
\end{IEEEproof}

\subsection{Proof of Proposition \ref{prop:partial_order}}\label{proof:partial_order}

\begin{IEEEproof}
	Due to dependencies among the results, we first prove (\ref{eq:Se_equal_SEQNQ}) in the result \ref{it:order:equality}, then prove the result \ref{it:order:inequality}, finally prove (\ref{eq:wstar_equal_wEQ}) in the result \ref{it:order:equality}.  
	
	\textbf{Part 1: we first prove (\ref{eq:Se_equal_SEQNQ}) in the result \ref{it:order:equality}}. Given $(\bv,\lambda)\in\Real^{2m}\times\Real_{++}$, for convenience of derivation, let us define $\bar{\sv}\coloneqq \se(\bv,\lambda)$ and $\bar{\Ecal}\coloneqq \supp(\bar{\sv})$. 
	
	From Lemma \ref{lemma:SEQNQ_property}, we have $\Spoly(\bar{\sv};\bv,\lambda)\subset\Se(\bv,\lambda)$, hence to prove (\ref{eq:Se_equal_SEQNQ}), it is sufficient to prove $\Spoly(\bar{\sv};\bv,\lambda)\supset\Se(\bv,\lambda)$. Let ${\wv}\in\Se(\bv,{\lambda})$, then we have
	\begin{equation}\label{temp:smap_equal_csign}
		\sOPT(\wv;\bv,\lambda)=\se(\bv,\lambda)=:\bar{\sv}
	\end{equation}
	from Def. \ref{def:Es_mapping}. To prove $\wv\in\Spoly(\bar{\sv};\bv,\lambda)$, we only need to show that $\wv$ satisfies $\bar{\sv}$-(\ref{eq:EQ}) and $\bar{\sv}$-(\ref{eq:NQ}) at $(\bv,\lambda)$.

	From (\ref{temp:smap_equal_csign}) and Def. \ref{def:Es_mapping}, we can derive
	\begin{align}
		&(\forall i\in\bar{\Ecal}) & \cv_i^\Tr(\bv-\Dm\Cm\wv)&=\lambda[\bar{\sv}]_i, \label{temp:sbar_EQa}\\
		&(\forall i\in\neg\bar{\Ecal}) & \abs{\cv_i^\Tr(\bv-\Dm\Cm\wv)}&\neq\lambda.\label{temp:xi_i_neq_lambda}
	\end{align}
	Combining (\ref{temp:xi_i_neq_lambda}) with the $i$th (\ref{eq:OPT}) condition implies
	\begin{equation}\label{temp:sbar_EQb}
		(\forall i\in\neg\bar{\Ecal})\;\;[\wv]_i=0.
	\end{equation}
	Putting (\ref{temp:sbar_EQa}) and (\ref{temp:sbar_EQb}) together, we can verify that $\wv$ satisfies $\bar{\sv}$-(\ref{eq:EQ}) at $(\bv,\lambda)$. 
	
	In addition, for every $i\in\bar{\Ecal}$, we have
	\begin{itemize}
		\item if $i\not\in \supp(\wv)$, then we naturally have $[\bar{\sv}]_i[\wv]_i=0$,
		
		\item if $i\in\supp(\wv)$, then from (\ref{eq:OPT}) we have
		\begin{equation*}
			\cv_i^\Tr(\bv-\Dm\Cm\wv)=\lambda\,\sign([\wv]_i),
		\end{equation*}
		combining which with (\ref{temp:sbar_EQa}) yields $[\bar{\sv}]_i=[\sign(\wv)]_i$, hence $[\bar{\sv}]_i[\wv]_i> 0$,
	\end{itemize}
	thus we verify that $\wv$ satisfies $\bar{\sv}$-(\ref{eq:NQ-E}) at $(\bv,\lambda)$. Moreover, since $\wv\in\Se(\bv,\lambda)$ satisfies (\ref{eq:OPT}), we have
	\begin{equation*}
		\cv_i^\Tr(\bv-\Dm\Cm\wv)\in[-\lambda,\lambda]
	\end{equation*}
	for every $i\in\braces{1,2,\dots,2n}$, hence $\wv$ satisfies $\bar{\sv}$-(\ref{eq:NQ-nE}) at $(\bv,\lambda)$. Accordingly, $\wv$ satisfies $\bar{\sv}$-(\ref{eq:NQ}) at $(\bv,\lambda)$.
	
	Combining the discussion above yields $\Spoly(\bar{\sv};\bv,\lambda)\supset\Se(\bv,\lambda)$. Hence we have $\Spoly(\bar{\sv};\bv,\lambda)=\Se(\bv,\lambda)$, which proves the equality (\ref{eq:Se_equal_SEQNQ}) in the result \ref{it:order:equality}.
	
	\vspace{0.5em}
	\noindent
	\textbf{Part 2: next we prove the result \ref{it:order:inequality}}.
	
	For $(\bv,\lambda)\in\Real^{2m}\times\Real_{++}$ and $\wv\in\Se(\bv,\lambda)$, let us define $\bar{\sv}\coloneqq \sOPT(\wv;\bv,\lambda)$, $\bar{\Ecal}\coloneqq\supp(\bar{\sv})$. Then from the definition of $\sOPT(\wv;\bv,\lambda)$ in Def. \ref{def:Es_mapping} and the fact that $\wv\in\Se(\bv,\lambda)$ satisfies (\ref{eq:OPT}), we have
	\begin{align*}
		\supp(\wv) &\subset \supp(\bar{\sv}),\\
		[\sign(\wv)]_{\supp(\wv)} &=[\bar{\sv}]_{\supp(\wv)}.
	\end{align*}
	Thus we have the inequality
	\begin{align}
		\sign(\wv)\preceq \bar{\sv}\coloneqq \sOPT(\wv;\bv,\lambda)=\se(\bv,\lambda), \label{temp:sign_w_preceq_sOPT_w}
	\end{align}	
	where the last equality follows $\wv\in\Se(\bv,\lambda)$ and Def. \ref{def:Es_mapping}.
	
	For arbitrary $\sv\in\braces{+1,0,-1}^{2n}$, let $\Ecal\coloneqq \supp(\sv)$.
	If $\wv\in\Spoly(\sv;\bv,\lambda)$, then since $\wv$ satisfies $\sv$-(\ref{eq:EQ-E}) at $(\bv,\lambda)$, we have $\Ecal\subset\bar{\Ecal}$ and the following equality
	\begin{align*}
		\lambda[\sv]_{\Ecal}={\Cm_{\Ecal}}^\Tr(\bv-\Dm\Cm\wv)=\lambda[\bar{\sv}]_{\Ecal},
	\end{align*}
	where the second equality follows from Def. \ref{def:Es_mapping}. Hence one can derive $\sv\preceq \bar{\sv}$. Since $\wv$ satisfies ${\sv}$-(\ref{eq:EQ-nE}) at $(\bv,\lambda)$, we have $\supp(\wv)\subset\Ecal$. Combining these results with (\ref{temp:sign_w_preceq_sOPT_w}) implies
	\begin{align*}
		\sign(\wv) &\preceq\bar{\sv},\\
		\sv &\preceq \bar{\sv},\\
		\supp(\wv)&\subset\Ecal\coloneqq\supp(\sv),
	\end{align*}
	thus from the definition of $\preceq$ (cf. Def. \ref{def:partial_order_of_s}), we have
	\begin{equation*}
		\sign(\wv)\preceq\sv\preceq \bar{\sv}\coloneqq\se(\bv,\lambda).
	\end{equation*}
	
	Reversely, if $(\Ecal,\sv)$ satisfies the generalized inequality above, then we can derive the following:
	\begin{itemize}
		\item since $\supp(\wv)\subset \Ecal$, we have $[\wv]_{\neg\Ecal}=\mathbf{0}$, hence $\wv$ satisfies $\sv$-(\ref{eq:EQ-nE}) at $(\bv,\lambda)$,
		
		\item since $\wv\in\Se(\bv,\lambda)$, from Def. \ref{def:Es_mapping} we have
		\begin{equation*}
			{\Cm_{\bar{\Ecal}}}^\Tr(\bv-\Dm\Cm\wv)=\lambda\brackets{\bar{\sv}}_{\bar{\Ecal}},
		\end{equation*}
		combining which with $\sv\preceq\bar{\sv}$ yields 
		\[{\Cm_{{\Ecal}}}^\Tr(\bv-\Dm\Cm\wv)=\lambda\brackets{\bar{\sv}}_{{\Ecal}}=\lambda[\sv]_{\Ecal},\] hence $\wv$ satisfies $\sv$-(\ref{eq:EQ-E}) at $(\bv,\lambda)$,
		
		\item since $\wv\in\Se(\bv,\lambda)$, for every $i\in\{1,\dots,2n\}$, we have the following from (\ref{eq:OPT}):
		\begin{equation*}
			\abs{\cv_i^\Tr(\bv-\Dm\Cm\wv)}\leq \lambda
		\end{equation*}
		which implies that $\wv$ satisfies $\sv$-(\ref{eq:NQ-nE}) at $(\bv,\lambda)$,
		
		\item since $\wv\in\Se(\bv,\lambda)$, from (\ref{eq:Se_equal_SEQNQ}) in the result \ref{it:order:equality} (which has been proven in the Part 1), we have $\wv\in\Spoly(\bar{\sv};\bv,\lambda)$, hence $\wv$ satisfies $\bar{\sv}$-(\ref{eq:NQ-E}) at $(\bv,\lambda)$, i.e.,
		\begin{align*}
			(\forall i\in\bar{\Ecal})\quad& \brackets{\bar{\sv}}_{i}[\wv]_i\geq 0,
		\end{align*}
		combining which with $\sv\preceq\bar{\sv}\coloneqq\se(\bv,\lambda)$ yields that $\wv$ satisfies $\sv$-(\ref{eq:NQ-E}) at $(\bv,\lambda)$.
	\end{itemize}
	Hence $\wv$ satisfies $\sv$-(\ref{eq:EQ}) and $\sv$-(\ref{eq:NQ}) simultaneously at $(\bv,\lambda)$, i.e., $\wv\in\Spoly(\sv;\bv,\lambda)$.
	
	Combining the discussion above proves the result \ref{it:order:inequality}.

	\vspace{0.5em}
	\noindent
	\textbf{Part 3: next we prove (\ref{eq:wstar_equal_wEQ}) in the result \ref{it:order:equality}}. Define $\wstar\coloneqq \wstar(\bv,\lambda)$, $\sstar\coloneqq\sign(\wstar)$ and $\Estar\coloneqq\supp(\sstar)$ for simplicity.
	
	Since $\wstar\in\Se(\bv,\lambda)$, from the result \ref{it:order:inequality} (which has been proven in the Part 2) we have 
	\begin{equation}\label{temp:wEQ_equal_wstar:1}
		\wstar\in\Spoly(\sstar;\bv,\lambda)\subset\Se(\bv,\lambda),
	\end{equation}
	where the second inclusion follows from Lemma \ref{lemma:SEQNQ_property}. Hence $\wstar$ satisfies $\sstar$-(\ref{eq:EQ}) at $(\bv,\lambda)$, which is equivalent to
	\begin{equation}\label{eq:wstar_in_EQ}
		\begin{cases}
			{\Cm_{\Estar}}^{\Tr}\Dm\Cm_{\Estar}[\wstar]_{\Estar}={\Cm_{\Estar}}^{\Tr}\bv-\lambda[\sstar]_{\Estar},\\
			[\wstar]_{\neg\Estar}=\mathbf{0}.
		\end{cases}
	\end{equation}
	One can verify that $\wstar$ is the least squares solution of (\ref{eq:wstar_in_EQ}) (i.e., $\wstar=\wEQ(\sstar;\bv,\lambda)$), if and only if $[\wstar]_{\Estar}$ satisfies
	\begin{equation}\label{eq:wstar_equal_wEQ_condition}
		[\wstar]_{\Estar}\in\Null\parens{{\Cm_{\Estar}}^{\Tr}\Dm\Cm_{\Estar}}^{\perp}.
	\end{equation}
	Hereafter, we will show that the inclusion (\ref{eq:wstar_equal_wEQ_condition}) holds.
	
	From Cor. \ref{cor:wstar}, $\wstar$ is the unique minimum norm element in $\Se(\bv,\lambda)\supset\Spoly(\sstar;\bv,\lambda)\ni\wstar$. Hence $\wstar$ is the unique solution of the following convex program:
	\begin{equation*}
		\begin{aligned}
			\underset{\wv\in\Real^{2n}}{\text{minimize}}\quad & \frac{1}{2}\norm{\wv}^2_2\\
			\text{subject to}\quad & \wv\in\Spoly(\sstar;\bv,\lambda).
		\end{aligned}
	\end{equation*}
	From Def. \ref{def:EQ_NQ_quantities}, the convex program above is equivalent to
	\begin{subequations}\label{temp:w_opt_necessary_conditions2}
		\begin{align}
			\underset{\wv\in\Real^{2n}}{\text{minimize}}\quad & \frac{1}{2}\norm{\wv}^2_2\label{temp:KKT_analysis_first_line}\\
			\text{subject to}\quad & {\Cm_{\Estar}}^{\Tr}(\bv-\Dm\Cm\wv)=\lambda[\sstar]_{\Estar},\\
			&[\wv]_{\neg\Estar}=\mathbf{0},\label{temp:w_opt_necessary_conditions2b}\\
			& [\sstar]_{\Estar}\odot [\wv]_{\Estar}\geq \mathbf{0},\\
			&-\lambda\mathbf{1}\leq{\Cm_{\neg\Estar}}^{\Tr}(\bv-\Dm\Cm\wv) \leq \lambda\mathbf{1}.
		\end{align}
	\end{subequations}
	From (\ref{temp:w_opt_necessary_conditions2b}), we know that every solution $\wv$ of (\ref{temp:w_opt_necessary_conditions2}) satisfies $[\wv]_{\neg\Estar}=\mathbf{0}$. Define $\zv\coloneqq [\wv]_{\Estar}$ and notice that $\Cm\wv=\Cm_{\Estar}\zv$, then (\ref{temp:w_opt_necessary_conditions2}) is further equivalent to $[\wv]_{\neg\Estar}=\mathbf{0}$ and
	\begin{subequations}\label{temp:w_opt_necessary_conditions3}
		\begin{align}
			\underset{\zv\in\Real^{|\Estar|}}{\text{minimize}}\quad & \frac{1}{2}\norm{\zv}^2_2\\
			\text{subject to}\quad & {\Cm_{\Estar}}^{\Tr}(\bv-\Dm\Cm_{\Estar}\zv)=\lambda[\sstar]_{\Estar},\label{temp:w_opt_necessary_conditions3b}\\
			& [\sstar]_{\Estar}\odot \zv\geq \mathbf{0},\\
			&{\Cm_{\neg\Estar}}^{\Tr}(\Dm\Cm_{\Estar}\zv-\bv)\leq \lambda\mathbf{1}\\
			&{\Cm_{\neg\Estar}}^{\Tr}(\bv-\Dm\Cm_{\Estar}\zv) \leq \lambda\mathbf{1}.
		\end{align}
	\end{subequations}
	The Lagrangian function of (\ref{temp:w_opt_necessary_conditions3}) is
	\begin{align}
		L\parens{\zv,\muv,\nuv}=&\frac{1}{2}\norm{\zv}^2_2+\muv^{\Tr}\parens{{\Cm_{\Estar}}^{\Tr}(\bv-\Dm\Cm_{\Estar}\zv)-\lambda[\sstar]_{\Estar}}\nonumber\\
		&-\nuv_1^{\Tr}\parens{[\sstar]_{\Estar}\odot\zv}\nonumber\\
		&+\nuv_2^{\Tr}\parens{{\Cm_{\neg\Estar}}^{\Tr}(\Dm\Cm_{\Estar}\zv-\bv)-\lambda\mathbf{1}}\nonumber\\
		&+\nuv_3^{\Tr}\parens{{\Cm_{\neg\Estar}}^{\Tr}(\bv-\Dm\Cm_{\Estar}\zv)-\lambda\mathbf{1}}\label{temp:Lagrangian}
	\end{align}
	where $\muv$ and $\nuv\coloneqq (\nuv_1,\nuv_2,\nuv_3)$ are the Lagrange multipliers. Hence from the KKT condition, we can yield the following necessary conditions of optimality (cf. \cite[Sec. 5.5.3]{boyd2004}):
	\begin{subequations}
		\begin{align}
			\nabla_{\zv} L(\zv,\muv,\nuv)&=\mathbf{0}.\label{temp:w_opt_necessary_conditions4a}\\
			\nuv_1\odot([\sstar]_{\Estar}\odot\zv)&=\mathbf{0}.\label{temp:w_opt_necessary_conditions4b}
		\end{align}
	\end{subequations}
	Expanding the expression of $\nabla_{\zv}L(\zv,\muv,\nuv)$ in (\ref{temp:w_opt_necessary_conditions4a}) yields
	\begin{equation}\label{temp:w_opt_necessary_conditions5}
		\zv={\Cm_{\Estar}}^{\Tr}\Dm^{\Tr}\parens{\Cm_{\Estar}\muv+\Cm_{\neg\Estar}(\nuv_3-\nuv_2)}+[\sstar]_{\Estar}\odot\nuv_1.
	\end{equation}
	Notice that if $\zv$ is the solution of (\ref{temp:w_opt_necessary_conditions3}), then since $\Estar\coloneqq\supp(\wstar)$, we have
	\[\zv=[\wstar]_{\Estar}
	=[\wstar]_{\supp(\wstar)},\]
	which means that every component in $\zv$ is nonzero. Since every component in $[\sstar]_{\Estar}$ is also nonzero, from (\ref{temp:w_opt_necessary_conditions4b}) we can derive that $\nuv_1=\mathbf{0}$, substituting which into (\ref{temp:w_opt_necessary_conditions5}) yields
	\begin{equation*}
		\zv={\Cm_{\Estar}}^{\Tr}\Dm^{\Tr}(\Cm_{\Estar}\muv+\Cm_{\neg\Estar}(\nuv_3-\nuv_2))\in\Col\parens{{\Cm_{\Estar}}^{\Tr}}.
	\end{equation*}
	Combining the inclusion above with $\zv\coloneqq[\wstar]_{\Estar}$ and Lemma \ref{lemma:blk_matrices} in Appx. \ref{app:technical_lemma} yields
	\begin{align*}
		[\wstar]_{\Estar}\in\Col\parens{{\Cm_{\Estar}}^{\Tr}\Dm^{\Tr}\Cm_{\Estar}}=\Null\parens{{\Cm_{\Estar}}^{\Tr}\Dm\Cm_{\Estar}}^{\perp}.
	\end{align*}
	Hence the inclusion (\ref{eq:wstar_equal_wEQ_condition}) holds, which implies that $\wstar(\bv,\lambda)$ is the least squares solution of $\sstar$-(\ref{eq:EQ}) at $(\bv,\lambda)$. Hence we have $\wstar(\bv,\lambda)=\wEQ(\sstar;\bv,\lambda)$, which proves (\ref{eq:wstar_equal_wEQ}) in the result \ref{it:order:equality}.
\end{IEEEproof}

\subsection{Proof of Lemma \ref{lemma:IEN_property}}
\label{proof:IEN_property}

\begin{IEEEproof}
	For simplicity, we define $\wEQ\coloneqq \wEQ(\sv;\bv,\lambda)$. From the definition of $\Iselect(\sv)$, $(\bv,\lambda)\in\Iselect(\sv)$ is equivalent to $\wEQ\in\Spoly(\sv;\bv,\lambda)$.
	
	\vspace{0.5em}
	\noindent
	\ref{it:IEN:wEQ_in_Se} We first prove that $\wEQ\in\Spoly(\sv;\bv,\lambda)$ is equivalent to the statement \ref{it:IEN:wEQ_in_Se}.
	
	If $\wEQ\in\Spoly(\sv;\bv,\lambda)$, then we have $\wEQ\in\Se(\bv,\lambda)$ directly from Lemma \ref{lemma:SEQNQ_property}. Reversely, if $\wEQ\in\Se(\bv,\lambda)$, then $\sOPT(\wEQ;\bv,\lambda)=\se(\bv,\lambda)$. Define $\bar{\sv}\coloneqq \se(\bv,\lambda)$, $\bar{\Ecal}\coloneqq\supp(\bar{\sv})$ and $\Ecal\coloneqq\supp(\sv)$ for abbreviation. Notice that $\wEQ$ satisfies $\sv$-(\ref{eq:EQ-E}) at $(\bv,\lambda)$, we have \[\cv_i^\Tr(\bv-\Dm\Cm\wEQ)=\lambda[\sv]_i\]
	for every $i\in\Ecal$, which implies that
	\begin{align*}
		\Ecal\subset\bar{\Ecal} \text{ and }
		[\sv]_{\Ecal}=[\bar{\sv}]_{\Ecal}.
	\end{align*}
	Hence we have $\sv\preceq\bar{\sv}\coloneqq \se(\bv,\lambda)$. Moreover, since $\wEQ$ satisfies $\sv$-(\ref{eq:EQ-nE}) at $(\bv,\lambda)$, we have $\supp(\wEQ)\subset \Ecal$. Combining the discussion above with (\ref{eq:inequality_in_Se}), we obtain the following relations:
	\begin{align*}
		\sv &\preceq \bar{\sv},\\
		\supp(\wEQ) &\subset\Ecal\coloneqq\supp(\sv),\\
		\sign(\wEQ) &\preceq \se(\bv,\lambda)=:\bar{\sv}.
	\end{align*}
	One can verify that the relations above yields
	\begin{equation*}
		\sign(\wEQ) \preceq \sv\preceq \bar{\sv}\coloneqq\se(\bv,\lambda).
	\end{equation*}
	Since $\wEQ\in\Se(\bv,\lambda)$, this implies $\wEQ\in\Spoly(\sv;\bv,\lambda)$ from Prop. \ref{prop:partial_order} \ref{it:order:inequality}. Hence $\wEQ\in\Spoly(\sv;\bv,\lambda)$ is equivalent to the statement \ref{it:IEN:wEQ_in_Se}.
	
	\vspace{0.5em}
	\noindent
	\ref{it:IEN:wEQ_min_in_SEQNQ} Next we show that $\wEQ\in\Spoly(\sv;\bv,\lambda)$ is equivalent to the statement \ref{it:IEN:wEQ_min_in_SEQNQ}. Notice that the latter statement evidently implies the former one, we only need to prove the reverse.
	
	If $\wEQ\in\Spoly(\sv;\bv,\lambda)$, then $\wEQ$ must satisfy $\sv$-(\ref{eq:EQ}) at $(\bv,\lambda)$. From the definition of $\wEQ$ in Def. \ref{def:EQ_NQ_quantities}, $\wEQ$ must be the minimum $\ell_2$-norm element in the solution set of $\sv$-(\ref{eq:EQ}) at $(\yv,\lambda)$, which is a super set of $\Spoly(\bv,\lambda;\sv)$. Hence $\wEQ$ is the minimum $\ell_2$-norm element in $\Spoly(\sv;\bv,\lambda)$.
	Hence $\wEQ\in\Spoly(\sv;\bv,\lambda)$ is equivalent to the statement \ref{it:IEN:wEQ_min_in_SEQNQ}.
	
	\vspace{0.5em}
	Finally, we show that the values of two discrete functions $\sign(\wEQ(\sv;\bv,\lambda))$ and $\sOPT(\wEQ(\sv;\bv,\lambda);\bv,\lambda)$ respectively remain constant for all $(\bv,\lambda)\in\interior(\Iselect(\sv))$. From Lemma \ref{lemma:computation_wEQ_IEN} \ref{it:comp:IEN} in Appx. \ref{app:computation}, $\Iselect(\sv)$ is the set of $(\bv,\lambda)\in\Real^{2m}\times\Real$ satisfying the following linear program:
	\begin{subequations}
		\begin{empheq}[left=\empheqlbrace]{align}
			&(\forall i\in\Ecal) && [\sv]_i[\wEQ(\sv;\bv,\lambda)]_i\geq 0, \tag{\ref{eq:EN-E}}&\\
			&(\forall i\in\neg\Ecal) && \abs{\xi_i(\wEQ(\sv;\bv,\lambda))}\leq \lambda,\tag{\ref{eq:EN-nE}}&\\
			& && \lambda>0, \nonumber &
		\end{empheq}
	\end{subequations}
	where $\xi_i(\wv)\coloneqq \cv_i^\Tr(\bv-\Dm\Cm\wv)$. Let us define 
	\begin{align*}
		f_i(\bv,\lambda)\coloneqq [\sv]_i[\wEQ(\sv;\bv,\lambda)]_i, && g_i(\bv,\lambda)\coloneqq \xi_i(\wEQ(\sv;\bv,\lambda)).
	\end{align*}
	Then $f(\cdot,\cdot)$ and $g(\cdot,\cdot)$ are both affine functions. One can verify that for every $i\in\Ecal$, either of the following two cases happens for the $i$th (\ref{eq:EN-E}) condition:
	\begin{enumerate}[label=\arabic*)]
		\item $f_i(\bv,\lambda)\equiv 0$, then since $[\sv]_i\neq 0$, we have
		\begin{align*}
			(\forall (\bv,\lambda)\in\Real^{2m}\times\Real)\quad\quad[\sign(\wEQ(\sv;\bv,\lambda))]_i\equiv 0.
		\end{align*}
		
		\item $f_i(\bv,\lambda)\not\equiv 0$, in which case the $i$th (\ref{eq:EN-E}) condition defines a half-space of $\Real^{2m}\times\Real$, and strict inequality of the $i$th (\ref{eq:EN-E}) condition holds for every $(\bv,\lambda)\in\interior(\Iselect(\sv))$, thus we have
		\begin{align*}
			(\forall(\bv,\lambda)\in\interior(\Iselect(\sv)))\quad\quad[\sign(\wEQ(\sv;\bv,\lambda))]_i\equiv [\sv]_i.
		\end{align*}
	\end{enumerate}
	Combining the discussion above with the fact that $\wEQ(\sv;\bv,\lambda)$ satisfies $\sv$-(\ref{eq:EQ-nE}) at every $(\bv,\lambda)\in\Iselect(\sv)$, one can verify that $\sign(\wEQ(\sv;\cdot,\cdot))$ keeps constant within $\interior(\Iselect(\sv))$.
	
	Similarly, for every $i\in\neg\Ecal$, either of the following two cases happens for the $i$th (\ref{eq:EN-nE}) condition:
	\begin{enumerate}[label=\arabic*)]
		\item $g_i(\bv,\lambda)\equiv +\lambda$ or $-\lambda$, then from Def. \ref{def:Es_mapping}, for every $(\bv,\lambda)\in\Real^{2m}\times\Real_{++}$, we always have
		\begin{align*}
			[\sOPT(\wEQ(\sv;\bv,\lambda);\bv,\lambda)]_i\equiv +1\text{ or }{-1}.
		\end{align*}
		
		\item $g_i(\bv,\lambda)\not\equiv +\lambda$ and $g_i(\bv,\lambda)\not\equiv-\lambda$, in which case the $i$th (\ref{eq:EN-nE}) condition defines a nonempty polytope in $\Real^{2m}\times\Real$, and strict inequality of the $i$th (\ref{eq:EN-nE}) condition holds for every $(\bv,\lambda)\in\interior(\Iselect(\sv))$. From Def. \ref{def:Es_mapping}, this implies that for all $(\bv,\lambda)\in\interior(\Iselect(\sv))$,
		\begin{align*}
			[\sOPT(\wEQ(\sv;\bv,\lambda);\bv,\lambda)]_i\equiv 0.
		\end{align*}
	\end{enumerate}
	Combining the discussion above with the fact that $\wEQ(\sv;\bv,\lambda)$ satisfies $\sv$-(\ref{eq:EQ-E}) at every $(\bv,\lambda)\in\Iselect(\sv)$, one can verify that $\sOPT(\wEQ(\sv;\bv,\lambda);\bv,\lambda)$ keeps constant within $\interior(\Iselect(\sv))$.
\end{IEEEproof}

\section{Derivation and Analysis of E-LARS}
\label{app:ELARS_deri_ana}

\subsection{Derivation of the E-LARS Iteration}
\label{app:ELARS_deri}
Suppose that Thm. \ref{thm:correspondence} \ref{it:corres:GP} holds for the linear zones of $\sv$ and $\sadj$. Then according to Lemma \ref{lemma:computation_wEQ_IEN} \ref{it:comp:IEN} in Appx. \ref{app:computation}, the zone-switching time $\tadj$ is exactly the time when
\begin{align*}
	\wEQ(\sv;t)\coloneqq\wEQ(\sv;\yv(t),\lambda(t))
\end{align*} 
leaves $\Iselect(\sv)$. Accordingly,
\begin{equation*}
	\begin{aligned}
		t_+\coloneqq& \sup_{t\in\Real} \braces{t\;\middle\vert\; \wstar(t)\in\interior(\Iselect(\sv))}
	\end{aligned}
\end{equation*}
gives a reasonable estimate of $\tadj$, which can be computed by Lemma \ref{lemma:compuation_IEN_line} in Appx. \ref{app:computation}. Moreover, let $\epsilon>0$ be an infinitesimally small number. Then one can imagine that as $t$ increases from $t_+-\epsilon$ to $t_++\epsilon$, the following happens:
\begin{enumerate}[label=\arabic*)]
	\item For $t=t_+-\epsilon$, $\wEQ(\sv;t)\coloneqq\wEQ(\sv;\yv(t),\lambda(t))$ satisfies $\sv$-(\ref{eq:NQ}) at $t$, i.e.,
	\begin{subequations}
		\makeatletter
		\def\@currentlabel{EN'}
		\makeatother
		\label{eq:ENp}
		\renewcommand{\theequation}{EN'-\alph{equation}}
		\begin{empheq}[left=\empheqlbrace]{align}
			&(\forall i\in\Ecal) && [\sv]_i[\wEQ(\sv;t)]_i\geq 0,  \label{eq:ENp-E}&\\
			&(\forall i\in\neg\Ecal) && |\xi_i\parens{\wEQ(\sv;t)}|\leq \lambda(t),\label{eq:ENp-nE}&
		\end{empheq}
	\end{subequations}
	where $\Ecal\coloneqq\supp(\sv)$, $\xi_i(\wv)\coloneqq \cv_i^\Tr(\bv(t)-\Dm\Cm\wv)$.
	
	\item As $t$ increases infinitesimally from $t_+-\epsilon$ to $t_++\epsilon$, there emerges some $i\in\Ical_{\text{break}}^{\text{a}}$ or $i\in\Ical_{\text{break}}^{\text{b}}$ with 
	\begin{align*}
		\Ical_{\text{break}}^{\text{a}} &\coloneqq\braces{i\in\Ecal\mid [\sv]_i[\wEQ(\sv;t_++\epsilon)]_i<0},\\
		\Ical_{\text{break}}^{\text{b}} &\coloneqq \braces{i\in\neg\Ecal\mid \abs{\xi_i(\wEQ(\sv;t_++\epsilon))}>\lambda(t_++\epsilon)},
	\end{align*}
	such that for every $i\in\Ical_{\text{break}}^{\text{a}}$ (resp. $i\in\Ical_{\text{break}}^{\text{b}}$), the $i$th (\ref{eq:ENp-E}) (resp. (\ref{eq:ENp-nE})) condition breaks at $t=t_++\epsilon$.
	
	\item For $t=t_++\epsilon$, the value of $\sv$ has to change to $\sadj$ to maintain the validity of (\ref{eq:ENp}) at $t$. 
\end{enumerate}

Based on this observation, given $\sv$, if we can modify the value of $\sv$ such that it satisfies (\ref{eq:ENp}) at $t_++\epsilon$, then the modified value of $\sv$ should give a reasonable estimate of $\sadj$. A greedy modification of $\sv$ is to remove every condition in (\ref{eq:ENp}) that breaks at $t_++\epsilon$ from the expression of (\ref{eq:ENp}), i.e., remove every $i\in\Ical_{\text{break}}^{\text{a}}$ (resp. $i\in\Ical_{\text{break}}^{\text{b}}$)  from $\Ecal$ (resp. $\neg\Ecal$). 

See the following procedure for this modification.

\begin{procedure}[deletion-insertion process]\label{proc:deletion_insertion}
	Given the input zone indicator $\sv$, let $\Ecal\coloneqq\supp(\sv)$. Define
	\begin{align*}
		\Idel^+ &\coloneqq\{i\in\Ecal \mid [\wEQ(\sv;t_+)]_i=0\},\\
		\Ijoin^+ &\coloneqq\{i\in\neg\Ecal \mid \abs{\xi_i(\wEQ(\sv;t_+))}=\lambda(t_+)\}.
	\end{align*}
	From previous discussion, one can imagine that $\Idel^+=\Ical_{\text{break}}^{\text{a}}$ and $\Ijoin^+=\Ical_{\text{break}}^{\text{b}}$.
	We conduct the following deletion and insertion operations on $\sv$ to fit $\sv$ into (\ref{eq:ENp}) at $t=t_++\epsilon$:
	\begin{itemize}
		\item \textbf{Deletion:} for every $i\in\Idel^+$, we set $[\sv]_i\gets 0$ to remove the $i$th (\ref{eq:ENp-E}) condition that breaks at $t_++\epsilon$.
		
		\item \textbf{Insertion:} for every $i\in\Ijoin^+$, we set
		\begin{align*}
			[\sv]_i \gets\sign(\xi_i(\wstar(t_+)))=\sign(\xi_i(\wEQ(\sv;t_+)))
		\end{align*}  
		to remove the $i$th (\ref{eq:ENp-nE}) condition that breaks at $t_++\epsilon$, and to ensure that $\wstar(t_+)$ satisfies the $i$th $\sv$-(\ref{eq:EQ-E}) condition at $(\yv(t_+),\lambda(t_+))$.
	\end{itemize}
	Let the modified value of $\sv$ be $\sv_+$, then we set $\sv_+$ as the output of the E-LARS iteration, i.e., an estimate of $\sadj$.
\end{procedure}

To implement the E-LARS iteration, we still need to address the computation of $\Idel^+$ and $\Ijoin^+$. Notice that for every $i\in\Ecal$ (resp. $i\in\neg\Ecal$), the interval of $t\in\Real$ satisfying the $i$th (\ref{eq:ENp-E}) (resp. (\ref{eq:ENp-nE})) condition can be computed in closed-form (see Lemma \ref{lemma:compuation_IEN_line} in Appx. \ref{app:computation}). Exploiting the results therein, we can obtain a closed-form expression of $\Idel^+$ and $\Ijoin^+$ in Procedure \ref{proc:deletion_insertion}, by the following steps:
\begin{enumerate}[label=\arabic*)]
	\item \textbf{Step 1:} for every $i\in\Ecal$, compute $t^a_i$ by Lemma \ref{lemma:compuation_IEN_line} \ref{it:tmax:ta}, which gives the supremum value of $t\in\Real$ such that $(\yv(t),\lambda(t))$ satisfies the $i$th (\ref{eq:ENp-E}) condition.
	
	\item \textbf{Step 2:} for every $i\in\neg\Ecal$, compute $t^b_i$ by Lemma \ref{lemma:compuation_IEN_line} \ref{it:tmax:tb}, which gives the supremum value of $t\in\Real$ such that $(\yv(t),\lambda(t))$ satisfies the $i$th (\ref{eq:ENp-nE}) condition.
	
	\item \textbf{Step 3:} compute the zone-switching time $t_+$, which equals to $\sup(\Tzone(\sv))$ in Lemma \ref{lemma:compuation_IEN_line} \ref{it:tmax:tc}:
	\begin{align*}
		t_+=\sup(\Tzone(\sv))= \min\braces{\min_{i\in\Ecal}t^a_i,\min_{i\in\neg\Ecal}t^b_i,t^c},
	\end{align*}
	where $t^c\coloneqq f_{\mathrm{tmax}}(-\Deltalmd,-\lambda_0)$.
	
	\item \textbf{Step 4}: compute the deleted and inserted index sets
	\begin{align*}
		\Idel^+ &= \braces{i\in\Ecal\mid t^a_i=t_+},\\
		\Ijoin^+ &= \braces{i\in\neg\Ecal\mid t^b_i=t_+}.
	\end{align*}
\end{enumerate}

\begin{figure*}
	\centering
	\begin{minipage}{\textwidth}
		\begin{algorithm}[H]
			\caption{A single iteration of the E-LARS algorithm}\label{alg:LARS_sGMC}
			\textbf{Problem Parameters:} $(\Am,\rho)\in\Real^{m\times n}\times[0,1)$, $(\bv_0,\lambda_0)\in\Real^{2m}\times\Real$, $(\Deltab,\Deltalmd)\in\Real^{2m}\times\Real$\\
			\textbf{Input:} an indicator $\sv$ of some linear zone \\
			\textbf{Output:} an estimate $\sv_+$ of the indicator $\sadj$ of an adjacent linear zone
			\begin{algorithmic}[1]
				\LineComment{\textbf{Compute $\Cm$ and $\Dm$ (Lemma. \ref{lemma:OPT} in Sec. \ref{sec:optimality_condition})}} 
				\State $\Cm\gets\blkdiag(\Am,\sqrt{\rho}\Am)$
				\State $\Dm\gets\begin{bmatrix}(1-\rho)\Imat_m & \sqrt{\rho}\Imat_m \\ -\sqrt{\rho}\Imat_m & \Imat_m \end{bmatrix}$
				\LineComment{\textbf{Compute $t^a_i$, $t^b_i$ and $t_{+}$ by Lemma \ref{lemma:compuation_IEN_line} in Appx. \ref{app:computation}}}
				\State $[\pv]_{\Ecal}\gets -\parens{{\Cm_{\Ecal}}^{\Tr}\Dm\Cm_{\Ecal}}^{\dagger}\parens{\Cm^{\Tr}_{\Ecal}\Deltab-[\sv]_{\Ecal}\Deltalmd }$
				\State $[\qv]_{\Ecal}\gets\parens{{\Cm_{\Ecal}}^{\Tr}\Dm\Cm_{\Ecal}}^{\dagger}\parens{\Cm^{\Tr}_{\Ecal}\bv_0-[\sv]_{\Ecal}\lambda_0 }$ 	
				\State $\uv\gets \Deltab+\Dm\Cm_{\Ecal}[\pv]_{\Ecal}$ 
				\State $\vv\gets \bv_0 -\Dm\Cm_{\Ecal}[\qv]_{\Ecal}$  
				\For{$i\in\Ecal$}
				\State $t_i^a\gets f_{\mathrm{tmax}}([\sv]_i[\pv]_i,[\sv]_i[\qv]_i),$
				\EndFor
				\For{$i\in\neg\Ecal$}
				\State $t^b_i\gets\min\Big\{f_{\mathrm{tmax}}\parens{-\cv_i^\Tr\uv-\Deltalmd,\lambda_0+\cv_i^\Tr\vv},\; f_{\mathrm{tmax}}\parens{\cv_i^\Tr\uv-\Deltalmd,\lambda_0-\cv_i^\Tr\vv}\Big\},$
				\EndFor
				
				\State $t_+\gets \min\braces{\min_{i\in\Ecal}t^a_i,\;\min_{i\in\neg\Ecal}t^b_i,\;f_{\mathrm{tmax}}(-\Deltalmd,-\lambda_0)}$
				
				\LineComment{\textbf{Compute the deleted and inserted index sets, and obtain $\sv_+$}}
				
				\State $\Idel^+\gets\{i\in\Ecal \mid t^a_i= t_+\}$
				
				\State $\Ijoin^+\gets\{i\in\neg\Ecal \mid t^b_i= t_+\}$
				
				\For{$i=1,2,\dots,2n$}
				\State $[\sv_+]_i\gets \begin{cases}
					0, & \text{if }i\in\Idel^+,\\
					\sign\parens{\cv_i^\Tr(\uv+\vv t_+)}, & \text{if }i\in\Ijoin^+,\\
					[\sv]_i, & \text{otherwise},
				\end{cases}$
				\EndFor
			\end{algorithmic}
		\end{algorithm}
	\end{minipage}
\end{figure*}

See Algorithm \ref{alg:LARS_sGMC} for the pseudocode of the E-LARS iteration, and see Table \ref{tab:complexity_of_LARS_sGMC} for its step-by-step complexity analysis.

\begin{table*}
	\centering
	\caption{Time and space complexity of a single E-LARS iteration (Algorithm \ref{alg:LARS_sGMC}).}
	\label{tab:complexity_of_LARS_sGMC}
	\begin{tabular}{ccc}
		\hline
		Object of computation & Time complexity & Space complexity \\
		\hline
		$(\pv,\qv)$ (line 5-6) & $\Ocal\parens{\abs{\Ecal}^2m+\abs{\Ecal}^3}$ & $\Ocal\parens{n}$ \\
		$(\uv,\vv)$ (line 7-8) & $\Ocal\parens{m\abs{\Ecal}}$ & $\Ocal(m)$\\
		$(t^a_i)_{i\in\Ecal}$ (line 9-11) & $\Ocal\parens{\abs{\Ecal}}$ & $\Ocal\parens{\abs{\Ecal}}$\\
		$(t^b_i)_{i\in\neg\Ecal}$ (line 12-14) & $\Ocal\parens{m\abs{\neg\Ecal}}$ & $\Ocal\parens{\abs{\neg\Ecal}}$\\
		$\Idel^+$ and $\Ijoin^+$ (line 17-18) & $\Ocal\parens{n}$ & $\Ocal\parens{n}$\\
		$\sv_+$ (line 19-21) & $\Ocal\parens{\abs{\Idel^+}+\abs{\Ijoin^+}}$ & $\Ocal\parens{n}$\\
		In total & $\Ocal\parens{mn+\abs{\Ecal}^2m+\abs{\Ecal}^3}$ & $\Ocal\parens{m+n}$ \\
		\hline
	\end{tabular}
\end{table*}

\subsection{Proof of Theorem \ref{thm:E_LARS}}\label{proof:E_LARS}
\begin{IEEEproof}
	 Define $\wstar(t)\coloneqq \wstar(\bv(t),\lambda(t))$ and $\se(t)\coloneqq \se(\bv(t),\lambda(t))$. Our proof is based on two key results:
	\begin{itemize}
		\item The min-norm solution path $\wstar(t)$ is continuous,
		
		\item there exists $\epsilon>0$ such that for every $t\in (t_+-\epsilon,t_+)$, $\sv$ satisfies the following equality (Thm. \ref{thm:correspondence} \ref{it:corres:zone})
		\begin{align}
			\sv&\equiv\sign(\wstar(t))\equiv\se(t),\label{eq:s_equality}
		\end{align}		
		whilst for every $t\in (t_+,t_++\epsilon)$, $\sadj$ satisfies
		\begin{align}
			\sadj&\equiv\sign(\wstar(t))\equiv\se(t).\label{eq:s_kp1_equality}
		\end{align}		
	\end{itemize}
	
	For $i\in\braces{1,2,\dots,2n}$, let us study the condition that guarantees $[\sadj]_i=[\sv_+]_i$.
	
	If $i\not\in\Idel^+$ and $i\not\in\Ijoin^+$ (i.e., $i$ is not involved in the deletion-insertion operations), then $[\sadj]_i=[{\sv}_+]_i$ always holds, which can be proven directly from the two results above:
	\begin{itemize}
		\item If $i\in\Ecal\coloneqq\supp(\sv)$ and $i\not\in\Idel^+$, then from the deletion-insertion process, $\wstar(t_+)$ satisfies but does not attain equality of the $i$th $\sv$-(\ref{eq:NQ-E}) condition at $t_+$, i.e.,
		\[[\sv]_i[\wstar(t_+)]_i>0.\]
		From the continuity of $\wstar(\cdot)$, we know that there exists some $t'\in(t_+,t_++\epsilon)$ such that
		\[[\sign(\wstar(t'))]_i=[\sv]_i=[{\sv}_+]_i,\]
		which implies $[\sadj]_i=[{\sv}_+]_i$ from (\ref{eq:s_kp1_equality}).
		
		\item If $i\in\neg\Ecal$ and $i\not\in\Ijoin^+$, then from the deletion-insertion process, $\wstar(t_+)$ satisfies but does not attain equality of the $i$th $\sv$-(\ref{eq:NQ-nE}) condition at $t_+$, i.e.,
		\[\abs{\cv_i^\Tr(\bv(t_+)-\Dm\Cm\wstar(t_+))}<\lambda(t_+).\]
		Since $f(t)\coloneqq \abs{\cv_i^\Tr(\bv(t)-\Dm\Cm\wstar(t))}-\lambda(t)$ is a continuous function, there exists $t'\in(t_+,t_++\epsilon)$ such that 
		\[\abs{\cv_i^\Tr(\bv(t')-\Dm\Cm\wstar(t'))}<\lambda(t').\]
		From the definition of $\se(t)$ (Def. \ref{def:Es_mapping}), we have
		\[[\se(t')]_i=0=[\sv]_i=[\sv_+]_i,\]
		which implies $[\sadj]_i=[\sv_+]_i$ from (\ref{eq:s_kp1_equality}).
	\end{itemize}
	Thus correctness of the LARS-sGMC iteration only depends on whether $[\sadj]_i=[\sv_+]_i$ holds for $i\in\parens{\Idel^+\cup\Ijoin^+}$.

	\vspace{0.5em}
	\noindent
	\ref{it:LARS:sticky} In the sequel, we show that $[\sadj]_i=[\sv_+]_i$ holds for every $i\in\parens{\Idel^+\cup\Ijoin^+}$ if and only if the sticky-equality assumption holds:
	\begin{itemize}
		\item If $i\in\Idel^+\subset\Ecal$, then $[\sv_+]_i=0$. Thus from (\ref{eq:s_kp1_equality}), we have $[\sadj]_i=[\sv_+]_i$ if and only if
		\begin{equation}\label{eq:NQ-E-equality}
			[\sign(\wstar(t))]_i=0
		\end{equation}
		holds for every $t\in (t_+,t_++\epsilon)$. Notice that $i\in\Ecal$, hence $[\sv]_i\neq 0$, which implies from (\ref{eq:s_equality}) that (\ref{eq:NQ-E-equality}) does not hold for every $t\in(t_+-\epsilon,t_+)$. In addition, $i\in\Idel^+$ implies that (\ref{eq:NQ-E-equality}) holds at $t\coloneqq t_+$, thus (\ref{eq:NQ-E-equality}) holds at every $t\in(t_+,t_++\epsilon)$ if and only if there exists a neighborhood $U$ of $t_+$ such that (\ref{eq:NQ-E-equality}) holds for every $t\in U$.
		
		\item If $i\in\Ijoin^+\subset\neg\Ecal$, then \[[\sv_+]_i=\sign\parens{\cv_i^\Tr(\bv(t_+)-\Dm\Cm\wstar(t_+))}\neq 0.\]
		From (\ref{eq:s_kp1_equality}), we have $[\sadj]_i=[\sv_+]_i$ if and only if\footnote{If $[\sadj]_i=[\sv_+]_i\neq 0$, then from (\ref{eq:s_kp1_equality}), $[\se(t)]_i\neq 0$ holds for every $\lambda\in(t_+,t_++\epsilon)$, thus (\ref{eq:NQ-nE-equality}) holds; conversely, if (\ref{eq:NQ-nE-equality}) holds for every $\lambda\in(t_+,t_++\epsilon)$, then since $f(t)=\cv_i^\Tr(\bv(t)-\Dm\Cm\wstar(t))$ is a continuous function satisfying $\abs{f(t_+)}\neq \lambda(t_+)>0$, one can verify that at the instant that $t$ increases across $t_+$, $f(t)$ cannot change its sign, which implies $[\se(t)]_i=\sign(f(t_+))=[\sv_+]_i$ for every $\lambda\in(t_+,t_++\epsilon)$, hence we can derive $[\sadj]_i=[\sv_+]_i$ from (\ref{eq:s_kp1_equality}).}
		\begin{equation}\label{eq:NQ-nE-equality}
			\abs{\cv_i^\Tr(\bv(t)-\Dm\Cm\wstar(t))}=\lambda(t)
		\end{equation}
		holds for every $\lambda\in(t_+,t_++\epsilon)$. Notice that $i\in\neg\Ecal$, hence $[\sv]_i=0$, which implies from (\ref{eq:s_equality}) that (\ref{eq:NQ-nE-equality}) does not hold for every $\lambda\in(t_+-\epsilon,t_+)$. In addition, $i\in\Ijoin^+$ implies that (\ref{eq:NQ-nE-equality}) holds at $t\coloneqq t_+$, thus (\ref{eq:NQ-nE-equality}) holds at every $\lambda\in(t_+,t_++\epsilon)$ if and only if there exists a neighborhood $V$ of $t_+$ such that (\ref{eq:NQ-nE-equality}) holds for every $t\in V$.		
	\end{itemize}
	
	Combining the discussion above, we can conclude that $\sadj=\sv_+$ if and only if the sticky-equality assumption holds.

	\vspace{0.5em}
	\noindent
	\ref{it:LARS:one} Finally, we show that the one-at-a-time assumption is a more restrictive condition than the sticky-equality assumption.
	
	Suppose that the one-at-a-time assumption holds, then there exists only one index $i$ that is involved in the deletion-insertion operation. From our previous analysis, the LARS-sGMC iteration is correct if and only if $[\sadj]_i=[\sv_+]_i$:
	\begin{itemize}
		\item If $i\in\Idel^+\subset\Ecal$, then since $\sv\neq\sadj$, we must have $[\sadj]_i\neq[\sv]_i$, thus $[\sadj]_i\in\braces{0,-[\sv]_i}$. If $[\sadj]_i=-[\sv]_i$, then from (\ref{eq:s_equality}) and Def. \ref{def:Es_mapping}, we have
		\[{\cv_i^\Tr(\bv(t)-\Dm\Cm\wstar(t))}=\lambda(t)[\sv]_i\]
		for every $\lambda\in(t_+-\epsilon,t_+)$. Without loss of generality, let us assume $[\sv]_i=+1$, then taking $t\uparrow t_+$ on both sides of the equation above yields
		\[{\cv_i^\Tr(\bv(t_+)-\Dm\Cm\wstar(t_+))}=\lambda(t_+)[\sv]_i>0.\]
		Similarly, from (\ref{eq:s_kp1_equality}), taking $t\downarrow t_+$ for $t\in(t_+,t_++\epsilon)$ implies
		\[{\cv_i^\Tr(\bv(t_+)-\Dm\Cm\wstar(t_+))}=\lambda(t_+)[\sadj]_i<0,\]
		which leads to a contradiction. Hence $[\sadj]_i=0=[\sv_+]_i$.
		
		\item If $i\in\Ijoin^+\subset\neg\Ecal$, then since $\sv\neq\sadj$, we must have $[\sadj]_i\neq[\sv]_i$, thus $[\sadj]_i\in\braces{+1,-1}$. Since $[\sadj]_i\neq 0$, from (\ref{eq:s_kp1_equality}) and Def. \ref{def:Es_mapping}, we have
		\begin{align*}
			\cv_i^\Tr(\bv(t)-\Dm\Cm\wstar(t))=[\sadj]_i \lambda(t),
		\end{align*}
		for every $\lambda\in (t_+,t_++\epsilon)$. Taking $t\downarrow t_+$ on both sides of the equation above yields
		\begin{align*}
			\cv_i^\Tr(\bv(t_+)-\Dm\Cm\wstar(t_+))=[\sadj]_i \lambda(t_+),
		\end{align*}
		since $\lambda(t_+)>0$, we have
		\[[\sadj]_i=\sign\parens{\cv_i^\Tr(\bv(t_+)-\Dm\Cm\wstar(t_+))}=[\sv_+]_i.\]
	\end{itemize}
Combining the discussion above, we can yield that the one-at-a-time assumption implies correctness of the LARS-sGMC iteration, which is equivalent to the sticky-equality assumption. 

To see that the reverse does not hold, it is sufficient to check the following counterexample: let us consider a LASSO problem (i.e., $\rho=0$, $\rv=\mathbf{0}_{m}$ in the sGMC model) with $\Am\coloneqq \begin{bmatrix}
	1 & 1
\end{bmatrix}$. Suppose that $(y,\lambda)\in\Real\times\Real$ moves in the following line:
\begin{align*}
	y(t)\equiv 1, \;\; \lambda(t)=2-t,
\end{align*}
Define $T(\sv)\coloneqq \braces{t\in\Real\;\middle\vert\; (y(t),\lambda(t))\in\Iselect(\sv)}$.
Using Lemma \ref{lemma:compuation_IEN_line} from Appx. \ref{app:computation}, one can verify that the min-norm solution path $\xstar(y(t),\lambda(t))$ of this LASSO problem has two linear pieces with indicators as follows (since $\rho=0$, we omit the dual solution part in the indicators):
\begin{align*}
	\sv_0&\coloneqq \mathbf{0}_{2}\text{ with }\Tzone(\sv_0)=(-\infty,1],\\
	\sv_1&\coloneqq \mathbf{1}_{2}\text{ with }\Tzone(\sv_0)=[1,2).
\end{align*}
Let us run the LARS-sGMC iteration with its input being $\sv_0$, then one can verify that the output is $\sv_1$ which is correct, whereas the one-at-a-time assumption does not hold.
\end{IEEEproof}

\bibliographystyle{citations/IEEEtran}
\bibliography{citations/IEEEabrv,citations/refs}

\begin{thebibliography}{10}
\providecommand{\url}[1]{#1}
\csname url@samestyle\endcsname
\providecommand{\newblock}{\relax}
\providecommand{\bibinfo}[2]{#2}
\providecommand{\BIBentrySTDinterwordspacing}{\spaceskip=0pt\relax}
\providecommand{\BIBentryALTinterwordstretchfactor}{4}
\providecommand{\BIBentryALTinterwordspacing}{\spaceskip=\fontdimen2\font plus
\BIBentryALTinterwordstretchfactor\fontdimen3\font minus
  \fontdimen4\font\relax}
\providecommand{\BIBforeignlanguage}[2]{{%
\expandafter\ifx\csname l@#1\endcsname\relax
\typeout{** WARNING: IEEEtran.bst: No hyphenation pattern has been}%
\typeout{** loaded for the language `#1'. Using the pattern for}%
\typeout{** the default language instead.}%
\else
\language=\csname l@#1\endcsname
\fi
#2}}
\providecommand{\BIBdecl}{\relax}
\BIBdecl

\bibitem{candes2006}
E.~Candes, J.~Romberg, and T.~Tao, ``Robust uncertainty principles: exact
  signal reconstruction from highly incomplete frequency information,''
  \emph{IEEE Trans. Inf. Theory}, vol.~52, no.~2, pp. 489--509, 2006.

\bibitem{donoho2006}
D.~Donoho, ``Compressed sensing,'' \emph{IEEE Trans. Inf. Theory}, vol.~52,
  no.~4, pp. 1289--1306, 2006.

\bibitem{eldar2012}
Y.~C. Eldar and G.~Kutyniok, Eds., \emph{Compressed {Sensing}: {Theory} and
  {Applications}}.\hskip 1em plus 0.5em minus 0.4em\relax Cambridge: Cambridge
  University Press, 2012.

\bibitem{qaisar2013}
S.~Qaisar, R.~M. Bilal, W.~Iqbal, M.~Naureen, and S.~Lee, ``Compressive
  sensing: {From} theory to applications, a survey,'' \emph{J. Commun. Netw.},
  vol.~15, no.~5, pp. 443--456, 2013.

\bibitem{hastie2015}
T.~Hastie, R.~Tibshirani, and M.~Wainwright, \emph{Statistical learning with
  sparsity: the lasso and generalizations}.\hskip 1em plus 0.5em minus
  0.4em\relax Boca Raton: CRC Press, Taylor \& Francis Group, 2015, no. 143.

\bibitem{tibshirani1996}
R.~Tibshirani, ``Regression {Shrinkage} and {Selection} {Via} the {Lasso},''
  \emph{J. R. Stat. Soc., B: Stat. Methodol.}, vol.~58, no.~1, pp. 267--288,
  1996.

\bibitem{chen2001}
S.~S. Chen, D.~L. Donoho, and M.~A. Saunders, ``Atomic {Decomposition} by
  {Basis} {Pursuit},'' \emph{SIAM Review}, vol.~43, no.~1, pp. 129--159, 2001.

\bibitem{combettes2005}
P.~L. Combettes and V.~R. Wajs, ``Signal recovery by proximal forward-backward
  splitting,'' \emph{Multiscale Model. Simul.}, vol.~4, no.~4, pp. 1168--1200,
  2005.

\bibitem{beck2009}
A.~Beck and M.~Teboulle, ``A {Fast} {Iterative} {Shrinkage}-{Thresholding}
  {Algorithm} for {Linear} {Inverse} {Problems},'' \emph{SIAM J. Imaging Sci.},
  vol.~2, no.~1, pp. 183--202, 2009.

\bibitem{osborne2000}
M.~R. Osborne, B.~Presnell, and B.~A. Turlach, ``On the {LASSO} and its
  {Dual},'' \emph{J. Comput. Graph. Stat.}, vol.~9, no.~2, pp. 319--337, 2000.

\bibitem{efron2004}
B.~Efron, T.~Hastie, I.~Johnstone, and R.~Tibshirani, ``Least angle
  regression,'' \emph{Ann. Stat.}, vol.~32, no.~2, pp. 407 -- 499, 2004.

\bibitem{hastie2007}
T.~Hastie, J.~Taylor, R.~Tibshirani, and G.~Walther, ``Forward stagewise
  regression and the monotone lasso,'' \emph{Electron. J. Stat.}, vol.~1, no.
  none, Jan. 2007.

\bibitem{zou2007}
H.~Zou, T.~Hastie, and R.~Tibshirani, ``On the “degrees of freedom” of the
  lasso,'' \emph{Ann. Stat.}, vol.~35, no.~5, pp. 2173--2192, 2007.

\bibitem{mairal2012}
J.~Mairal and B.~Yu, ``Complexity {Analysis} of the {Lasso} {Regularization}
  {Path},'' in \emph{Proceedings of the 29th {International} {Coference} on
  {International} {Conference} on {Machine} {Learning}}, ser. {ICML}'12.\hskip
  1em plus 0.5em minus 0.4em\relax Omnipress, 2012, pp. 1835--1842.

\bibitem{tibshirani2012}
R.~J. Tibshirani and J.~Taylor, ``Degrees of freedom in lasso problems,''
  \emph{Ann. Stat.}, vol.~40, no.~2, pp. 1198 -- 1232, 2012.

\bibitem{tibshirani2013}
R.~J. Tibshirani, ``The lasso problem and uniqueness,'' \emph{Electron. J.
  Stat.}, vol.~7, 2013.

\bibitem{fan2001}
J.~Fan and R.~Li, ``Variable {Selection} via {Nonconcave} {Penalized}
  {Likelihood} and its {Oracle} {Properties},'' \emph{J. Am. Stat. Assoc.},
  vol.~96, no. 456, pp. 1348--1360, 2001.

\bibitem{zhang2010}
C.-H. Zhang, ``Nearly unbiased variable selection under minimax concave
  penalty,'' \emph{Ann. Stat.}, vol.~38, no.~2, 2010.

\bibitem{selesnick2017}
I.~Selesnick, ``Sparse {Regularization} via {Convex} {Analysis},'' \emph{IEEE
  Trans. Signal Process.}, vol.~65, no.~17, pp. 4481--4494, 2017.

\bibitem{lanza2019}
A.~Lanza, S.~Morigi, I.~W. Selesnick, and F.~Sgallari, ``Sparsity-{Inducing}
  {Nonconvex} {Nonseparable} {Regularization} for {Convex} {Image}
  {Processing},'' \emph{SIAM J. Imaging Sci.}, vol.~12, no.~2, pp. 1099--1134,
  2019.

\bibitem{abe2020}
J.~Abe, M.~Yamagishi, and I.~Yamada, ``Linearly involved generalized {Moreau}
  enhanced models and their proximal splitting algorithm under overall
  convexity condition,'' \emph{Inverse Probl.}, vol.~36, no.~3, p. 035012,
  2020.

\bibitem{alshabili2021}
A.~H. Al-Shabili, Y.~Feng, and I.~Selesnick, ``Sharpening {Sparse}
  {Regularizers} via {Smoothing},'' \emph{IEEE Open J. Signal Process.},
  vol.~2, pp. 396--409, 2021.

\bibitem{lanza2021}
A.~Lanza, S.~Morigi, I.~W. Selesnick, and F.~Sgallari, ``Convex {Non}-convex
  {Variational} {Models},'' in \emph{Handbook of {Mathematical} {Models} and
  {Algorithms} in {Computer} {Vision} and {Imaging}}.\hskip 1em plus 0.5em
  minus 0.4em\relax Cham: Springer International Publishing, 2021, pp. 1--57.

\bibitem{yata2022}
W.~Yata, M.~Yamagishi, and I.~Yamada, ``A constrained {LiGME} {Model} and {Its}
  {Proximal} {Splitting} {Algorithm} under overall convexity condition,''
  \emph{J. Appl. Numer.}, 2022.

\bibitem{liu2023}
X.~Liu, A.~J. Molstad, and E.~C. Chi, ``A convex-nonconvex strategy for grouped
  variable selection,'' \emph{Electronic Journal of Statistics}, vol.~17,
  no.~2, Jan. 2023.

\bibitem{zhang2023}
Y.~Zhang and I.~Yamada, ``A {Unified} {Framework} for {Solving} a {General}
  {Class} of {Nonconvexly} {Regularized} {Convex} {Models},'' \emph{IEEE Trans.
  Signal Process.}, vol.~71, pp. 3518--3533, 2023.

\bibitem{zhang2021}
------, ``{DC}-{LiGME}: {An} {Efficient} {Algorithm} for {Improved} {Convex}
  {Sparse} {Regularization},'' in \emph{55th {Asilomar} {Conference} on
  {Signals}, {Systems}, and {Computers}}, 2021, pp. 1348--1354.

\bibitem{zhang2022}
------, ``A {Unified} {Class} of {DC}-type {Convexity}-{Preserving}
  {Regularizers} for {Improved} {Sparse} {Regularization},'' in \emph{30th
  {European} {Signal} {Processing} {Conference} ({EUSIPCO})}, 2022, pp.
  2051--2055.

\bibitem{gregor2010}
K.~Gregor and Y.~LeCun, ``Learning {Fast} {Approximations} of {Sparse}
  {Coding},'' in \emph{International {Conference} on {Machine} {Learning}},
  2010.

\bibitem{rockafellar2009}
R.~T. Rockafellar and R.~J.-B. Wets, \emph{Variational analysis}, corr. 3nd
  print~ed.\hskip 1em plus 0.5em minus 0.4em\relax Berlin: Springer, 2009.

\bibitem{osborne2000LARS}
M.~Osborne, ``A new approach to variable selection in least squares problems,''
  \emph{IMA J. Numer. Anal.}, vol.~20, no.~3, pp. 389--403, 2000.

\bibitem{dossal2012}
C.~Dossal, ``A necessary and sufficient condition for exact sparse recovery by
  minimization,'' \emph{C. R. Math.}, vol. 350, no. 1-2, pp. 117--120, 2012.

\bibitem{rosset2007}
S.~Rosset and J.~Zhu, ``Piecewise {Linear} {Regularized} {Solution} {Paths},''
  \emph{Ann. Stat.}, vol.~35, no.~3, pp. 1012--1030, 2007.

\bibitem{tibshirani2011}
R.~J. Tibshirani and J.~Taylor, ``{THE} {SOLUTION} {PATH} {OF} {THE}
  {GENERALIZED} {LASSO},'' \emph{Ann. Stat.}, vol.~39, no.~3, pp. 1335--1371,
  2011.

\bibitem{tibshirani2011phd}
R.~J. Tibshirani, ``The solution path of the generalized lasso,'' Ph.D.
  dissertation, Stanford University, 2011.

\bibitem{berk2023}
A.~Berk, S.~Brugiapaglia, and T.~Hoheisel, ``{LASSO} {Reloaded}: {A}
  {Variational} {Analysis} {Perspective} with {Applications} to {Compressed}
  {Sensing},'' \emph{SIAM Journal on Mathematics of Data Science}, vol.~5,
  no.~4, pp. 1102--1129, 2023.

\bibitem{xiao2015}
W.~Xiao, Y.~Wu, and H.~Zhou, ``{ConvexLAR}: {An} {Extension} of {Least} {Angle}
  {Regression},'' \emph{J. Comput. Graph. Stat.}, vol.~24, no.~3, pp. 603--626,
  2015.

\bibitem{yukawa2016}
M.~Yukawa and S.-I. Amari, ``$\ell_p$-{Regularized} {Least} {Squares} $(0<p<1)$
  and {Critical} {Path},'' \emph{IEEE Trans. Inf. Theory}, vol.~62, no.~1, pp.
  488--502, 2016.

\bibitem{mishkin2022}
A.~Mishkin and M.~Pilanci, ``The {Solution} {Path} of the {Group} {Lasso},'' in
  \emph{{OPT} 2022: {Optimization} for {Machine} {Learning} ({NeurIPS} 2022
  {Workshop})}, 2022.

\bibitem{mishkin2023}
------, ``Optimal sets and solution paths of {ReLU} networks,'' in
  \emph{Proceedings of the 40th {International} {Conference} on {Machine}
  {Learning}}, ser. {ICML}'23.\hskip 1em plus 0.5em minus 0.4em\relax JMLR.org,
  2023.

\bibitem{hartung1982}
J.~Hartung, ``An extension of {Sion}'s minimax theorem with an application to a
  method for constrained games.'' \emph{Pacific Journal of Mathematics}, vol.
  103, no.~2, pp. 401 -- 408, 1982.

\bibitem{bauschke2017}
H.~H. Bauschke and P.~L. Combettes, \emph{Convex {Analysis} and {Monotone}
  {Operator} {Theory} in {Hilbert} {Spaces}}, 2nd~ed.\hskip 1em plus 0.5em
  minus 0.4em\relax Cham: Springer, 2017.

\bibitem{hogan1973}
W.~W. Hogan, ``Point-to-{Set} {Maps} in {Mathematical} {Programming},''
  \emph{SIAM Review}, vol.~15, no.~3, pp. 591--603, Jul. 1973.

\bibitem{rudin1976}
W.~Rudin, \emph{Principles of mathematical analysis}, 3rd~ed.\hskip 1em plus
  0.5em minus 0.4em\relax New York: McGraw-Hill, 1976.

\bibitem{mangasarian1987}
O.~L. Mangasarian and T.-H. Shiau, ``Lipschitz {Continuity} of {Solutions} of
  {Linear} {Inequalities}, {Programs} and {Complementarity} {Problems},''
  \emph{SIAM J. Control Optim.}, vol.~25, no.~3, pp. 583--595, 1987.

\bibitem{fukuda1997}
K.~Fukuda, T.~M. Liebling, and F.~Margot, ``Analysis of backtrack algorithms
  for listing all vertices and all faces of a convex polyhedron,''
  \emph{Computational Geometry}, vol.~8, no.~1, pp. 1--12, Jun. 1997.

\bibitem{boyd2004}
\BIBentryALTinterwordspacing
S.~Boyd and L.~Vandenberghe, \emph{Convex {Optimization}}, 1st~ed.\hskip 1em
  plus 0.5em minus 0.4em\relax Cambridge University Press, 2004. [Online].
  Available: \url{https://web.stanford.edu/~boyd/cvxbook/bv_cvxbook.pdf}
\BIBentrySTDinterwordspacing

\bibitem{tao2016}
T.~Tao, \emph{Analysis {I}}, 3rd~ed.\hskip 1em plus 0.5em minus 0.4em\relax
  Singapore: Springer Singapore, 2016.

\bibitem{lu2002}
T.-T. Lu and S.-H. Shiou, ``Inverses of 2 × 2 block matrices,'' \emph{Comput.
  Math. Appl.}, vol.~43, no. 1-2, pp. 119--129, 2002.

\bibitem{rockafellar1997}
R.~T. Rockafellar, \emph{Convex {Analysis}}.\hskip 1em plus 0.5em minus
  0.4em\relax Princeton, NJ: Princeton University Press, 1970.

\end{thebibliography}

\end{document}